\documentclass[11pt]{book} 

\usepackage{fancyhdr}

\fancypagestyle{plain}{\fancyhead{}}

\usepackage{amssymb,amscd,amsfonts,mathrsfs,amsmath}
\usepackage[all]{xy}

\input amssym.def

\def\double{\mathbb}

\def\cc{{\double C}}
\def\nn{{\double N}}
\def\zz{{\double Z}}

\def\rr{{\double R}}

\newtheorem{theorem}{Th\'eor\`eme}[section]
\newtheorem{lemma}[theorem]{Lemme}
\newtheorem{corollary}[theorem]{Corollaire}
\newtheorem{definition}[theorem]{D\'efinition}
\newtheorem{proposition}[theorem]{Proposition}
\newtheorem{remark}[theorem]{Remarque}

\newtheorem{example}[theorem]{Exemple}

\def\Kt{K^{\mathrm{top}}}
\def\Ka{K^{\mathrm{alg}}}

\def\res{\mathop{\mathrm{Res}}\limits_{z=0}}
\def\Pf{\mathop{\mathrm{Pf}}}

\def\cp{\rtimes}
\def\si{\sigma}
\def\Si{\Sigma}
\def\cinf{C^{\infty}}
\def\cinfc{C^{\infty}_c}

\newcommand{\be}{\begin{equation}}
\newcommand{\ee}{\end{equation}}
\newcommand{\beq}{\begin{eqnarray}}
\newcommand{\eeq}{\end{eqnarray}}
\newcommand{\om}{\omega}
\newcommand{\Om}{\Omega}
\newcommand{\al}{\alpha}
\def\nat{\natural}
\def\id{{\mathop{\mathrm{id}}}}

\newcommand{\la}{\lambda}
\newcommand{\Ec}{{\mathscr E}}
\newcommand{\Vc}{{\mathscr V}}

\newcommand{\non}{\nonumber}
\newcommand{\eps}{\varepsilon}

\newcommand{\Wc}{{\mathscr W}}
\newcommand{\Rc}{{\mathscr R}}
\newcommand{\Mc}{{\mathscr M}}
\newcommand{\Nc}{{\mathscr N}}
\newcommand{\Ic}{{\mathscr I}}
\newcommand{\Jc}{{\mathscr J}}

\newcommand{\Ind}{{\mathop{\mathrm{Ind}}}}
\def\ch{\mathrm{ch}}

\def\Td{\mathrm{Td}}

\def\tch{\mathrm{\,\slash\!\!\!\! \ch}}

\def\re{\mathrm{Re}}

\newcommand{\Tr}{{\mathop{\mathrm{Tr}}}}

\newcommand{\Ac}{{\mathscr A}}

\newcommand{\te}{\theta}

\newcommand{\Te}{\Theta}
\newcommand{\Det}{\textup{Det}}
\newcommand{\cqfd}{\hfill\rule{1ex}{1ex}}
\def\deb{\overline{\partial}}
\def\Id{\mathrm{Id}}
\def\zb{\overline{z}}

\def\d{\partial}
\def\dd{\mathrm{\bf d}}

\def\Hc{{\mathscr H}}

\def\Bc{{\mathscr B}}
\def\Cc{{\mathscr C}}
\def\Jc{{\mathscr J}}
\def\Kc{{\mathscr K}}
\def\Fc{{\mathscr F}}
\def\Pc{{\mathscr P}}

\def\bb{\overline{b}}

\def\psib{\overline{\psi}}

\def\hom{{\mathop{\mathrm{Hom}}}}
\def\dom{{\mathop{\mathrm{Dom}}}}
\def\ran{{\mathop{\mathrm{Ran}}}}

\def\End{{\mathop{\mathrm{End}}}}

\def\hotimes{\hat{\otimes}}

\def\Psit{\widetilde{\Psi}}

\def\gt{\tilde{g}}
\def\ut{\tilde{u}}

\def\Gh{\widehat{G}}

\def\et{\tilde{e}}

\def\Omh{\widehat{\Omega}}
\def\Th{\widehat{T}}
\def\chih{\widehat{\chi}}
\def\etah{\widehat{\eta}}

\def\Dh{\widehat{D}}
\def\Fh{\widehat{F}}
\def\Rch{\widehat{\mathscr R}}
\def\Mch{\widehat{\mathscr M}}
\def\Jch{\widehat{\mathscr J}}
\def\Ome{\Omega_{\epsilon}}

\def\Omc{\Omega_c}

\def\Xh{\widehat{X}}
\def\Jh{\widehat{J}}
\def\Tc{{\mathscr T}}
\def\Dc{{\mathscr D}}

\def\Kc{{\mathscr K}}

\def\mod{\ \mathrm{mod}\ }

\def\supp{\mathrm{supp}\,}

\def\eh{\hat{e}}

\def\uh{\hat{u}}

\def\Omb{\underline{\Omega}}
\def\zzb{\underline{\zz}}
\def\Hdr{H_{\mathrm{dR}}}
\def\Zdr{Z_{\mathrm{dR}}}

\def\SL{SL}
\def\GL{GL}

\begin{document}

\begin{center}
{\bf \Large Universit\'e Claude Bernard - Lyon 1 }

\vskip 1cm
{\Large  HABILITATION A DIRIGER DES RECHERCHES}

\vskip 1cm

{\Large Discipline: Math\'ematiques}

\vskip 5cm

{\huge \bf Sur les formules locales de l'indice}

\vskip 3cm

{\Large\bf Denis PERROT}

\vskip 1cm

{\Large \bf Novembre 2008}

\end{center}

\thispagestyle{empty}

\newpage
\thispagestyle{empty}

\thispagestyle{empty}

\tableofcontents

\chapter*{Introduction}
\addcontentsline{toc}{chapter}{Introduction}
\fancyhead[LO]{\bfseries Introduction}
\fancyhead[RE]{\bfseries Introduction}

Ce m\'emoire d'habilitation est la synth\`ese de travaux effectu\'es en th\'eorie de l'indice non-commutative \cite{P1}-\cite{P10}. Notre objectif est de fournir des \emph{formules locales} et d'\'etudier leurs relations avec les anomalies de la th\'eorie quantique des champs. On se place dans le cadre de la g\'eom\'etrie diff\'erentielle non-commutative et de l'homologie cyclique d\'evelopp\'ees par Connes \cite{C0, C1}. Un espace non-commutatif y est repr\'esent\'e par une alg\`ebre associative. En pratique il s'agit d'une alg\`ebre de Banach, ou de Fr\'echet, ou m\^eme plus g\'en\'eralement d'une alg\`ebre bornologique \cite{Me}. La $K$-th\'eorie et l'homologie cyclique d\'ecrivent les invariants de topologie alg\'ebrique d'un tel espace. De fa\c{c}on g\'en\'erale on peut dire que la th\'eorie de l'indice non-commutative consiste \`a \'etudier l'image de ces invariants sous l'action d'un bimodule de Kasparov \cite{Bl} (ou d'un quasihomomorphisme \cite{Cu}) suffisamment ``lisse''. Plus pr\'ecis\'ement la situation qui nous int\'eresse est la suivante. D\'esignons par $\Ic=\ell^p$ l'id\'eal de Schatten des op\'erateurs $p$-sommables sur un espace de Hilbert. Un $\Ac$-$\Bc$-bimodule $p$-sommable entre deux alg\`ebres de Fr\'echet $\Ac$ et $\Bc$ induit alors un morphisme d'image directe en $K$-th\'eorie topologique $\Kt_*(\Ic\hotimes\Ac)\to \Kt_*(\Ic\hotimes\Bc)$, o\`u le produit tensoriel projectif compl\'et\'e $\Ic\hotimes\cdot$ est une version $p$-sommable de stabilisation. On cherche alors \`a construire un caract\`ere de Chern dans la cohomologie cyclique bivariante de $\Ac$ et $\Bc$, de sorte que la fl\^eche induite en homologie cyclique $HC_*(\Ac)\to HC_*(\Bc)$ s'ins\`ere dans un diagramme commutatif 
\be
\vcenter{\xymatrix{
\Kt_*(\Ic\hotimes\Ac) \ar[r] \ar[d] & \Kt_*(\Ic\hotimes\Bc) \ar[d] \\
HC_*(\Ac) \ar[r] & HC_*(\Bc) }} \label{grr}
\ee
Il existe deux approches compl\'ementaires, chacune apportant son lot d'avantages et d'inconv\'enients. La premi\`ere est bas\'ee sur les propri\'et\'es abstraites de la $K$-th\'eorie et de la cohomologie cyclique bivariantes qui garantissent l'existence d'un caract\`ere de Chern ``universel'' muni des propri\'et\'es voulues. C'est la voie suivie par Cuntz dans le cas des alg\`ebres localement convexes \cite{Cu1, Cu2} ou par Puschnigg pour les $C^*$-alg\`ebres \cite{Pu}. Cette m\'ethode est enti\`erement satisfaisante d'un point de vue th\'eorique car elle garantit un r\'esultat beaucoup plus g\'en\'eral que (\ref{grr}), \`a savoir la compatibilit\'e entre le produit de Kasparov en $K$-th\'eorie bivariante et le produit de composition en cohomologie cyclique bivariante. Par contre elle ne donne pas v\'eritablement de formule concr\`ete pour le caract\`ere de Chern.\\
La deuxi\`eme approche renonce \`a construire un caract\`ere de Chern universel et se concentre sur des bimodules de Kasparov $p$-sommables munis de propri\'et\'es ad\'equates, tels que ceux qui apparaissent dans les situations d'origine g\'eom\'etrique. On peut alors donner des formules relativement simples, mais il n'est pas garanti qu'elles repr\'esentent le caract\`ere de Chern universel. Dans ce cas la commutativit\'e du diagramme (\ref{grr}) doit \^etre v\'erifi\'ee \emph{a posteriori}. C'est la voie que nous adoptons ici. Notre objectif est double. D'abord on d\'egage les conditions qui permettent de construire un caract\`ere de Chern concret et assurent l'existence de diagrammes commutatifs comme ci-dessus. Ensuite on explique comment obtenir des formules locales en s'inspirant des techniques de renormalisation en th\'eorie quantique des champs. En fait la diagonale $\Delta:\Kt_*(\Ic\hotimes \Ac)\to HC_*(\Bc)$ du diagramme (\ref{grr}) est l'exacte g\'en\'eralisation du calcul de l'anomalie chirale associ\'ee \`a une th\'eorie de jauge non-commutative. L'\'evaluation de l'image de $\Delta$ sur une classe de cohomologie cyclique de $\Bc$ donne alors une formule locale de l'indice. \\

Le chapitre \ref{ccar} d\'ecrit la construction du caract\`ere de Chern bivariant au moyen de superconnexions de Quillen \cite{Q1}. En r\'ealit\'e on donne deux formules. La premi\`ere repose sur le noyau de la chaleur $\exp(-tD^2)$ associ\'e \`a un op\'erateur de Dirac \cite{P4}; elle est donc adapt\'ee aux bimodules non-born\'es  ``$\te$-sommables''. L'autre est obtenue par un proc\'ed\'e de r\'etraction et fonctionne pour les bimodules born\'es $p$-sommables \cite{P5} v\'erifiant certaines conditions d'admissibilit\'e. Vu que l'on ne cherche pas ici \`a parler de $K$-th\'eorie, il suffit de se placer dans la cat\'egorie la plus g\'en\'erale, celle des alg\`ebres bornologiques. Le caract\`ere de Chern prend alors ses valeurs dans la cohomologie cyclique bivariante enti\`ere $HE^*(\Ac,\Bc)$, qui contient des cocycles de dimension infinie bien adapt\'es aux formules bas\'ees sur l'utilisation du noyau de la chaleur. De plus, dans des circonstances favorables la limite $t\downarrow 0$ permet d'obtenir un repr\'esentant local du caract\`ere de Chern. A titre d'exemple, nous \'etablissons au chapitre \ref{ceq} un th\'eor\`eme de l'indice pour les actions propres et isom\'etriques d'un groupe localement compact $G$ sur une vari\'et\'e de Riemann \cite{P6}. On d\'emontre que l'indice d'un op\'erateur elliptique $G$-invariant de type Dirac, qui d\'etermine une classe d'homologie cyclique sur l'alg\`ebre du groupe, est donn\'e par une formule de localisation aux points fixes de l'action. \\

A partir du chapitre \ref{csec} on se restreint aux alg\`ebres de Fr\'echet multiplicativement convexes, c'est-\`a-dire les limites projectives de suites d'alg\`ebres de Banach. Elles poss\`edent deux types distincts d'invariants: les invariants primaires, stables par homotopie diff\'erentiable tels que la $K$-th\'eorie topologique \cite{Ph} et l'homologie cyclique p\'eriodique, et les invariants secondaires tels que la $K$-th\'eorie multiplicative \cite{K1, K2} et les versions instables de l'homologie cyclique. Ces diff\'erents types d'invariants sont reli\'es par des suites exactes longues. En utilisant les formules obtenues pr\'ec\'edemment pour le caract\`ere de Chern d'un bimodule born\'e, on \'etablit qu'un quasihomomorphisme $p$-sommable, de parit\'e $p\mod 2$, muni de certaines propri\'et\'es d'admissibilit\'e induit des morphismes d'image directe pour les invariants primaires et secondaires tout en respectant les suites exactes. Cela se traduit par un diagramme commutatif
$$
\xymatrix{
\ldots \Kt_{n+1}(\Ic\hotimes\Ac) \ar[r] \ar[d] & HC_{n-1}(\Ac) \ar[r] \ar[d] & MK^{\Ic}_n(\Ac)  \ar[r] \ar[d]  & \Kt_n(\Ic\hotimes\Ac)  \ar[d] \ldots  \\
\ldots \Kt_{n+1-p}(\Ic\hotimes\Bc) \ar[r]  & HC_{n-1-p}(\Bc) \ar[r]  & MK^{\Ic}_{n-p}(\Bc)  \ar[r]  & \Kt_{n-p}(\Ic\hotimes\Bc) \ldots }
$$
avec $\Ic$ une alg\`ebre $p$-sommable et $MK^{\Ic}_n$ la $K$-th\'eorie multiplicative introduite dans \cite{P8}. L'entier $p$ est la \emph{dimension relative} du quasihomomorphisme. En cons\'equence on obtient non seulement le diagramme (\ref{grr}) en $K$-th\'eorie topologique mais aussi des diagrammes analogues reliant $K$-th\'eorie multiplicative et les versions instables d'homologie cyclique. Notons que le produit de Kasparov entre quasihomomorphismes n'est pas d\'efini dans ce contexte.  \\

La m\'ethode g\'en\'erale permettant d'\'etablir des formules locales pour le caract\`ere de Chern bivariant est expos\'ee au chapitre \ref{cano}, ainsi que sa relation avec les anomalies en th\'eorie quantique des champs \cite{P9}. Le lien entre th\'eorie de l'indice et anomalies n'est pas nouveau. Citons par exemple Atiyah et Singer \cite{AS2,S} ou plus r\'ecemment Mickelsson et coauteurs \cite{AM, CMM, M, MS}. Cependant l'approche que nous pr\'esentons ici est diff\'erente. On introduit la \emph{cocha\^{\i}ne \^eta renormalis\'ee} comme s\'erie formelle dans le complexe cyclique bivariant, dont le bord fournit automatiquement un repr\'esentant local du caract\`ere de Chern. Cette s\'erie formelle est reli\'ee tr\`es explicitement \`a la fonctionnelle d'action quantique d'une th\'eorie de jauge non-commutative, ce qui explique le lien avec les anomalies. L'avantage de cette m\'ethode est qu'elle laisse une grande libert\'e dans le choix de renormalisation: dans chaque situation de nature ``g\'eom\'etrique'' il existe un choix assez naturel qui donne un repr\'esentant local particulier du caract\`ere de Chern. Par exemple, la renormalisation z\^eta est utilisable en pr\'esence d'un op\'erateur de Dirac. Le caract\`ere de Chern est alors donn\'e par une somme de r\'esidus de fonctions z\^eta et g\'en\'eralise la formule de Connes et Moscovici valable pour les triplets spectraux \cite{CM95}. Ce n'est \'evidemment pas le seul choix possible. On illustre dans le chapitre \ref{cgrou} un autre type de renormalisation dans le cas d'un groupe op\'erant sur le plan complexe par \emph{transformations conformes}. Ici aucune m\'etrique riemannienne n'est pr\'eserv\'ee et l'introduction d'un op\'erateur de Dirac n'est pas naturelle. On peut n\'eanmoins renormaliser sans briser la sym\'etrie conforme, ce qui m\`ene encore une fois \`a une formule de l'indice localis\'ee aux points fixes. Elle fait intervenir des nombres de Lefschetz g\'en\'eralis\'es ainsi qu'une classe de Todd non-commutative bas\'ee sur le groupe d'automorphismes modulaires \cite{P10}. \\

\chapter{Caract\`ere de Chern bivariant}\label{ccar}
\fancyhead[LO]{\bfseries\rightmark}
\fancyhead[RE]{\bfseries\leftmark}

Ce chapitre pr\'esente deux formules candidates pour le caract\`ere de Chern d'un $\Ac$-$\Bc$-bimodule muni des propri\'et\'es ad\'equates. La premi\`ere est bas\'ee sur le noyau de la chaleur associ\'e \`a un op\'erateur de Dirac. Elle est donc adapt\'ee au bimodules non-born\'es $\te$-sommables \cite{P4}. La deuxi\`eme n'utilise que la phase de l'op\'erateur de Dirac et par cons\'equent est applicable aux bimodules born\'es $p$-sommables \cite{P5}. En fait ces deux caract\`eres de Chern sont reli\'es, au moins formellement, par un processus de r\'etraction et d\'efinissent la m\^eme classe de cohomologie cyclique bivariante.  \\
Pour se placer dans le cadre le plus g\'en\'eral possible, on consid\`ere la cat\'egorie des alg\`ebres \emph{bornologiques}. Le caract\`ere de Chern d'un $\Ac$-$\Bc$-bimodule vit alors dans la cohomologie cyclique bivariante enti\`ere $HE^*(\Ac,\Bc)$ \cite{Me}. Par commodit\'e nous rappelons dans la premi\`ere section quelques rudiments de th\'eorie cyclique enti\`ere. Les deux sections suivantes d\'etaillent la construction du caract\`ere de Chern respectivement dans le cas d'un bimodule non-born\'e $\te$-sommable et d'un bimodule born\'e $p$-sommable. Le mat\'eriel expos\'e ici est adapt\'e des articles \\

\noindent \cite{P4} D. Perrot: A bivariant Chern character for families of spectral triples, {\it Comm. Math. Phys.} {\bf 231} (2002) 45-95.\\

\noindent \cite{P5} D. Perrot: Retraction of the bivariant Chern character, {\it K-Theory} {\bf 31} (2004) 233-287.

\section{Alg\`ebres bornologiques}\label{sbo}

Rappelons qu'une \emph{bornologie} sur un $\cc$-espace vectoriel $\Vc$ est la donn\'ee d'un ensemble de parties de $\Vc$, dites born\'ees, v\'erifiant certains axiomes \cite{Me}. Il existe une notion de \emph{compl\'etude} au sens bornologique. L'exemple standard d'espace vectoriel bornologique est un espace vectoriel localement convexe $\Vc$ muni de sa bornologie dite de von Neumann, constitu\'ee des parties born\'ees pour toutes les semi-normes d\'efinissant la topologie de $\Vc$. \\
Soient $\Vc$ et $\Wc$ deux espaces vectoriels bornologiques. Une application lin\'eaire $\Vc \to \Wc$ est born\'ee si elle envoie les parties born\'ees de $\Vc$ sur les parties born\'ees de $\Wc$. L'ensemble des applications lin\'eaires born\'ees $\hom(\Vc,\Wc)$ est lui-m\^eme un espace vectoriel bornologique, complet si $\Wc$ l'est. On d\'efinit aussi le \emph{produit tensoriel bornologique compl\'et\'e} $\Vc\hotimes\Wc$ par une propri\'et\'e universelle de factorisation. Ce produit tensoriel est associatif. Notons que dans le cas des espaces de Fr\'echet, $\hom(\Vc,\Wc)$ est exactement l'espace des applications lin\'eaires continues de $\Vc$ vers $\Wc$, et le produit tensoriel bornologique $\hotimes$ co\"incide essentiellement avec le produit tensoriel projectif (modulo quelques subtilit\'es, voir \cite{Me}). \\
Une alg\`ebre bornologique compl\`ete $\Ac$ est un espace bornologique complet muni d'une application bilin\'eaire born\'ee associative $\Ac\times\Ac\to\Ac$. Par exemple si $\Ac$ est un espace de Fr\'echet, alors la multiplication est automatiquement (jointement) continue et $\Ac$ est aussi une alg\`ebre de Fr\'echet. Si $\Ac$ et $\Bc$ sont deux alg\`ebres bornologiques compl\`etes, leur produit tensoriel bornologique $\Ac\hotimes\Bc$ est encore une alg\`ebre bornologique compl\`ete. \\

L'homologie cyclique d'une alg\`ebre bornologique compl\`ete $\Ac$ est d\'efinie au moyen des formes diff\'erentielles non-commutatives \cite{C0}. Soit $\Ac^+=\Ac\oplus \cc$ l'alg\`ebre obtenue en ajoutant une unit\'e (et ce, m\^eme si $\Ac$ est d\'ej\`a unitaire). L'espace des $n$-formes diff\'erentielles non-commutatives est le produit tensoriel compl\'et\'e
\be
\Om^n\Ac = \Ac^+\hotimes \Ac^{\hotimes n} \quad  n>0\ ,\qquad  \Om^0\Ac=\Ac\ ,
\ee
et la somme directe $\Om\Ac=\bigoplus_{n\geq 0}\Om^n\Ac$ est un espace bornologique complet. Nous adopterons la notation standard $a_0 \dd a_1 \ldots \dd a_n \in \Om^n\Ac$ pour le produit tensoriel $a_0\otimes a_1\ldots \otimes a_n$ et $\dd a_1 \ldots \dd a_n \in \Om^n\Ac$ pour $1\otimes a_1\ldots \otimes a_n$. La diff\'erentielle $\dd:\Om^n\Ac\to \Om^{n+1}\Ac$ est une application lin\'eaire born\'ee de carr\'e nul. On introduit de mani\`ere classique l'op\'erateur de Hochschild $b:\Om^n\Ac\to \Om^{n-1}\Ac$ et de Connes $B:\Om^n\Ac\to\Om^{n+1}\Ac$, tous deux born\'es et v\'erifiant $b^2=bB+Bb=B^2=0$.\\
L'homologie cyclique enti\`ere de $\Ac$ s'obtient en compl\'etant $\Om\Ac$ dans la \emph{bornologie enti\`ere} \cite{Me}. Une cha\^{\i}ne enti\`ere est alors une collection de formes diff\'erentielles $\om_n\in \Om^n\Ac$ donn\'ees pour tout $n\in\nn$ et v\'erifiant une condition de croissance lorsque $n\to \infty$. Nous noterons $\Ome\Ac$ l'espace bornologique des formes diff\'erentielles enti\`eres ainsi obtenu. On montre que $b$ et $B$ s'\'etendent en des applications lin\'eaires born\'ees sur $\Ome\Ac$, qui devient donc un complexe bornologique un fois muni de la diff\'erentielle totale $b+B$, naturellement $\zz_2$-gradu\'e par le degr\'e pair/impair des formes diff\'erentielles.

\begin{definition}[\cite{Me}]
Soit $\Ac$ une alg\`ebre bornologique compl\`ete. Son homologie cyclique enti\`ere est l'homologie du complexe $\zz_2$-gradu\'e $\Ome\Ac$ muni de la diff\'erentielle $b+B$:
\be
HE_i(\Ac) = H_i(\Ome\Ac)\ ,\qquad i\in\zz_2\ .
\ee
La cohomologie cyclique enti\`ere bivariante de deux alg\`ebres bornologiques compl\`etes $\Ac$ et $\Bc$ est l'homologie du complexe $\zz_2$-gradu\'e des applications lin\'eaires born\'ees de $\Ome\Ac$ vers $\Ome\Bc$:
\be
HE^n(\Ac,\Bc) = H_i(\hom(\Ome\Ac,\Ome\Bc))\ ,\qquad i\in \zz_2\ .
\ee
\end{definition}
En particulier $HE_i(\Ac)=HE^i(\cc,\Ac)$, et $HE^i(\Ac,\cc)$ s'identifie \`a la cohomologie cyclique enti\`ere de $\Ac$ d\'efinie par Connes \cite{C2}. Rappelons enfin une description \'equivalente de l'homologie cyclique enti\`ere reli\'ee au formalisme de Cuntz et Quillen \cite{CQ1}. Pour toute alg\`ebre bornologique compl\`ete $\Ac$, Meyer d\'efinit dans \cite{Me} la notion d'extension analytique universelle 
$$
0\to \Jc \to \Rc \to \Ac \to 0 \ ,
$$
o\`u $\Rc$ est une alg\`ebre bornologique compl\`ete analytiquement quasi-libre et l'id\'eal $\Jc$ est analytiquement nilpotent. Un exemple d'extension universelle est donn\'e par \emph{l'alg\`ebre tensorielle analytique} $\Rc=\Tc\Ac$, compl\'etion de l'alg\`ebre tensorielle $T\Ac$ dans une bornologie appropri\'ee \cite{Me}. Une \'equivalence de Goodwillie g\'en\'eralis\'ee implique alors l'isomorphisme $HE_i(\Ac)=HE_i(\Rc)$. De plus le $X$-complexe
\be
X(\Rc) \ :\ \Rc\ \mathop{\rightleftarrows}^{\nat\dd}_{\bb}\ \Om^1\Rc_{\nat}
\ee 
avec $\Om^1\Rc_{\nat}=\Om^1\Rc/b\Om^2\Rc$ calcule l'homologie cyclique enti\`ere de $\Rc$. Ces isomorphismes sont induits par des \'equivalences d'homotopies de complexes $\zz_2$-gradu\'es:
$$
X(\Rc) \stackrel{\sim}{\longleftarrow} \Ome\Rc \stackrel{\sim}{\longrightarrow} \Ome\Ac\ .
$$
La fl\^eche de gauche est la projection du complexe des formes diff\'erentielles enti\`eres sur le $X$-complexe, tandis que la fl\^eche de droite est induite par l'homomorphisme $\Rc\to\Ac$. Lorsque $\Rc=\Tc\Ac$, nous avons construit dans \cite{P4} un inverse explicite $\gamma: X(\Tc\Ac)\to \Ome\Tc\Ac$ r\'ealisant l'\'equivalence de Goodwillie, dont il sera fait constamment usage. Rappelons enfin que par hypoth\`ese, toute extension d'alg\`ebres bornologiques est scind\'ee par une application \emph{lin\'eaire} born\'ee.

\section{Bimodules non born\'es}\label{snobo}

Soient $\Ac$ et $\Bc$ deux alg\`ebres bornologiques compl\`etes. Nous appellerons \emph{$\Ac$-$\Bc$-bimodule non-born\'e} tout triplet $(H,\rho,D)$ v\'erifiant les propri\'et\'es suivantes:
\begin{itemize}
\item $H$ est un espace bornologique complet $\zz_2$-gradu\'e. Le produit tensoriel compl\'et\'e $H_{\Bc}=H\hotimes\Bc$ est donc naturellement muni d'une structure de $\Bc$-module \`a droite, et l'on note $\End(H_{\Bc})$ l'alg\`ebre des endomorphismes born\'es de $H_{\Bc}$ qui commutent avec l'action de $\Bc$.
\item $\rho:\Ac\to \End(H_{\Bc})$ est un homomorphisme born\'e qui repr\'esente l'alg\`ebre $\Ac$ dans les endomorphismes de $H_{\Bc}$ de degr\'e pair. $H_{\Bc}$ est donc un $\Ac$-module \`a gauche.
\item $D:\dom(D)\subset H_{\Bc}\to H_{\Bc}$ est un endomorphisme non born\'e de degr\'e impair, commutant avec l'action de $\Bc$. On appelle $D$ un \emph{op\'erateur de Dirac}.
\item Le commutateur $[D,\rho(a)]$ s'\'etend en un endomorphisme pair dans $\End(H_{\Bc})$ pour tout $a\in\Ac$.
\end{itemize}
Implicitement on supposera l'existence d'un sous-espace dense $\Hc\subset H$, complet dans sa propre bornologie et tel que $D$ soit un endomorphisme born\'e du $\Bc$-module \`a droite $\Hc_{\Bc}=\Hc\hotimes\Bc$. Dans les exemples concrets $H$ est un espace de Hilbert mais ce point importe peu au niveau de g\'en\'eralit\'e consid\'er\'e ici. Pour construire le caract\`ere de Chern bivariant nous aurons besoin d'imposer l'existence de l'op\'erateur de la chaleur associ\'e au laplacien de Dirac $D^2$:
\begin{itemize}
\item L'op\'erateur $\exp(-t D^2) \in \End(H_{\Bc})$ existe en tant qu'endomorphisme pour tout $t\geq 0$ et v\'erifie l'\'equation de la chaleur $\frac{d}{dt}\exp(-tD^2) = -D^2\exp(-tD^2)$. 
\end{itemize}

\begin{example}\textup{L'exemple classique d'un bimodule non-born\'e est une vari\'et\'e fibr\'ee $M \stackrel{X}{\longrightarrow} B$ au-dessus d'une base compacte $B$, \`a fibre compacte; $\Ac=\cinf(M)$ et $\Bc=\cinf(B)$ sont des alg\`ebres de Fr\'echet commutatives de fonctions lisses; $D$ est une famille lisse d'op\'erateurs de Dirac op\'erant sur les sections d'un espace fibr\'e vectoriel $E\to X$ et param\'etr\'ee par la base; $H$ est l'espace de Hilbert des sections de carr\'e sommable de $E$ et $\Hc\subset H$ est l'espace de Fr\'echet des sections lisses. Dans cet exemple la fibration $M$ est triviale mais on peut toujours se ramener \`a un bimodule du type $H_{\Bc}=H\hotimes\Bc$ par trivialisation locale et partition de l'unit\'e.  }
\end{example}

Contrairement \`a la situation des bimodules non-born\'es en $K$-th\'eorie bivariante de Kasparov \cite{Bl}, on ne peut pas imposer \`a la ``r\'esolvante'' de l'op\'erateur de Dirac d'\^etre compacte, car cette notion n'a pas de sens en bornologie. Par contre, la notion d'op\'erateurs \emph{tra\c{c}ables} y est bien d\'efinie en g\'en\'eral. Nous avons introduit dans \cite{P4} l'alg\`ebre des endomorphismes tra\c{c}ables $\ell^1(H_{\Bc})$, qui est naturellement un bimodule sur $\End(H_{\Bc})$. La supertrace d'op\'erateurs sur $H$ induit une application lin\'eaire born\'ee \cite{P4}
\be
\Tr: \ell^1(H_{\Bc}) \to \Bc\ ,\label{tr}
\ee
qui est une \emph{supertrace partielle} sur $\ell^1(H_{\Bc})$ vu comme $\End(H_{\Bc})$-bimodule. Nous allons donc remplacer la condition de compacit\'e sur la r\'esolvante de l'op\'erateur de Dirac par une condition de tra\c{c}abilit\'e sur l'op\'erateur de la chaleur. La formulation pr\'ecise est en fait un peu plus compliqu\'ee et m\`ene \`a la notion de $\te$-sommabilit\'e d\'efinie ci-dessous. \\

En rapport avec la p\'eriodicit\'e de Bott formelle on distingue deux types de bimodules, suivant leur parit\'e \cite{P4}. Un bimodule $(H,\rho,D)$ est \emph{pair} si l'espace $H$ se d\'ecompose en somme directe de deux sous-espaces distincts $H_+\oplus H_-$ relativement \`a sa $\zz_2$-graduation. Puisque l'image de $\rho$ est incluse dans la sous-alg\`ebre paire de $\End(H_{\Bc})$ et que $D$ est impair, on peut alors \'ecrire en notation matricielle
$$
H=\left(\begin{matrix} H_+ \\ H_- \end{matrix}\right) \ ,\quad \rho=\left(\begin{matrix} \rho_+ & 0 \\ 0 & \rho_- \end{matrix}\right) \ ,\quad D=\left(\begin{matrix} 0 & Q^* \\ Q & 0 \end{matrix}\right)\ .
$$
Ici $Q$ et $Q^*$ sont des op\'erateurs ind\'ependants, la notion d'adjonction n'ayant pas de sens pour le moment. Un bimodule $(H,\rho,D)$ est \emph{impair} si $H=L\otimes C_1$ est le produit tensoriel d'un espace trivialement gradu\'e $L$ avec l'alg\`ebre de Clifford $\zz_2$-gradu\'ee $C_1=\cc\oplus\cc\eps$, engendr\'ee par l'unit\'e $1$ en degr\'e zero et l'\'el\'ement $\eps$ en degr\'e un ($\eps^2=1$). Dans ce cas il existe un homomorphisme $\al:\Ac\to \End(L_{\Bc})$ et un endomorphisme non born\'e $Q$ sur $L_{\Bc}$ tels que 
$$
H=L\otimes C_1\ ,\qquad \rho = \al \otimes 1 \ ,\qquad D = Q\otimes \eps\ .
$$

La strat\'egie suivie dans \cite{P4} pour construire le caract\`ere de Chern $\ch(H,\rho,D)\in HE^*(\Ac,\Bc)$ consiste d'abord \`a relever $(H,\rho,D)$ en un bimodule sur des extensions universelles de $\Ac$ et $\Bc$. Nous avons choisi de travailler avec les alg\`ebres tensorielles analytiques dans \cite{P4}, mais il est en fait possible de g\'en\'eraliser la construction \`a des extensions quelconques sans trop d'effort. Choisissons donc une extension analytique universelle
$$
0 \to \Jc \to \Rc \to \Bc \to 0\ . 
$$
Elle induit une extension d'espaces vectoriels bornologiques $0\to H_{\Jc} \to H_{\Rc}\to H_{\Bc} \to 0$. Supposons d'abord que l'image de l'homomorphisme $\rho:\Ac\to \End(H_{\Bc})$ ainsi que l'op\'erateur de Dirac $D$ s'\'etendent en des endomorphismes du $\Bc^+$-module \`a droite $H_{\Bc^+}$, o\`u $\Bc^+$ est l'unitarisation de $\Bc$. Il existe alors un homomorphisme d'alg\`ebre canonique 
\be
\xymatrix{\End(H_{\Rc}) \ar[r] & \End(H_{\Bc}) \ar@/_1pc/@{.>}[l]_{\si} } \ , \label{morp}
\ee
scind\'e par une application \emph{lin\'eaire} born\'ee $\si$. Cette hypoth\`ese d'extension est assez restrictive. Dans certaines situations cependant, l'homomorphisme ci-dessus existe sans avoir besoin de passer par l'unitarisation, voir l'exemple \ref{ebott}. On supposera donc l'existence de (\ref{morp}) sans autre pr\'ecision. La propri\'et\'e universelle \cite{Me} de l'alg\`ebre tensorielle analytique $\Tc\Ac$ implique ensuite l'existence d'un homomorphisme born\'e $\rho_* : \Tc\Ac\to \End(H_{\Rc})$ en vertu du diagramme commutatif 
$$
\vcenter{\xymatrix{
0 \ar[r] & \Jc\Ac \ar[r] \ar[d]_{\rho_*} & \Tc\Ac \ar[r] \ar[d]_{\rho_*} & \Ac \ar[r] \ar[d]^{\rho} & 0 \\
0\ar[r] & \Nc^s \ar[r]  & \End(H_{\Rc}) \ar[r] & \End(H_{\Bc}) \ar[r] \ar@/_1pc/@{.>}[l]_{\si} & 0 }}
$$
o\`u $\Nc^s$ est le noyau de l'homomorphisme (\ref{morp}). La notation $^s$ rappelle que $\Nc^s$ est une alg\`ebre $\zz_2$-gradu\'ee ( = {\it s}upersym\'etrique). De fa\c{c}on analogue, l'op\'erateur de Dirac $D$ sur $H_{\Bc}$ se rel\`eve en un endomorphisme non-born\'e $\Dh$ sur $H_{\Rc}$.\\

On se place ensuite dans le formalisme des cocha\^{\i}nes d'alg\`ebre de Quillen \cite{Q2}. Soit $\Om^*\Rc=\Rc\oplus \Om^1\Rc$ l'alg\`ebre $\zz_2$-gradu\'ee des formes diff\'erentielles sur $\Rc$ tronqu\'ee en degr\'es $> 1$. Le $\Om^*\Rc$-module \`a droite $H_{\Om^*\Rc}$ est naturellement muni d'une diff\'erentielle. Dans \cite{P4} , nous avons d\'efini une compl\'etion bornologique de la cog\`ebre bar de $\Tc\Ac$ que l'on notera $\Cc$. Elle est munie de la codiff\'erentielle de Hochschild \cite{Q2}. D\'esignons par $\Om_1\Cc$ le bicomodule des 1-coformes universelles sur $\Cc$ et par $\Om_*\Cc= \Cc\oplus \Om_1\Cc$ la cog\`ebre $\zz_2$-gradu\'ee associ\'ee \cite{Q2}. L'espace des applications lin\'eaires born\'ees 
$$
\Fc = \hom(\Om_*\Cc, H_{\Om^*\Rc})
$$
est donc naturellement un module \`a droite sur l'alg\`ebre $\hom(\Om_*\Cc,\Om^*\Rc)$ munie du produit de convolution. $\Fc$ est aussi dot\'e d'une diff\'erentielle totale $\d$. L'observation fondamentale (\cite{P4}) est que l'homomorphisme $\rho_*:\Tc\Ac\to \End(H_{\Rc})$ et l'op\'erateur de Dirac $\Dh$ agissent par endomorphismes de degr\'e impair sur $\Fc$. R\'eunissons-les au sein d'une superconnexion de Quillen \cite{Q1}:
\be
\nabla =  \d + \rho_* + \Dh \ .
\ee
La courbure $\nabla^2=\d(\rho_*+\Dh) +[\Dh,\rho_*] +\Dh^2$ est donc un endomorphisme de degr\'e pair sur $\Fc$. Pr\'ecisons que $[\ ,\ ]$ est le commutateur gradu\'e. En gros, on cherche \`a obtenir le caract\`ere de Chern de $(H,\rho,D)$ en prenant l'exponentielle de cette courbure, comme dans le cas classique d'un espace fibr\'e vectoriel muni d'une connexion. Le lemme 7.1 de \cite{P4} montre comment d\'efinir l'op\'erateur de la chaleur $\exp(-t\Dh^2)$ au moyen d'un d\'eveloppement formel en s\'erie de Duhamel. La condition de $\te$-sommabilit\'e impose que l'exponentielle de la courbure $\nabla^2$ soit un endomorphisme \emph{tra\c{c}able} dans le sens suivant:

\begin{definition}[$\te$-sommabilit\'e \cite{P4}]
Un bimodule $(H,\rho,D)$ admissible relativement \`a une extension $\Rc$ est $\te$-sommable si la s\'erie de Duhamel 
$$
\exp(-\nabla^2):= \sum_{n\ge 0}(-)^n \int_{\Delta_n}dt_1\ldots dt_n\, e^{-t_0\Dh^2}\Te e^{-t_1\Dh^2}\ldots \Te e^{-t_n\Dh^2}\ ,
$$
avec $\Te=\d(\rho_*+\Dh)+[\Dh,\rho_*]$, d\'efinit un \'el\'ement pair de l'alg\`ebre de convolution $\hom(\Om_*\Cc, \ell^1(H_{\Om^*\Rc}))$ op\'erant par endomorphismes sur le $\hom(\Om_*\Cc,\Om^*\Rc)$-module \`a droite $\Fc$. 
\end{definition}
Ici $\Delta_n=\{(t_0,\ldots,t_n)\in [0,1]^{n+1}\ |\ \sum_i t_i=1\}$ d\'esigne le $n$-simplexe standard. En fait on ne s'int\'eresse qu'\`a la projection de $\exp(-\nabla^2)$ sur $\hom(\Om_1\Cc,\ell^1(H_{\Om^*\Rc}))$. On peut alors composer cette application lin\'eaire \`a l'entr\'ee par la cotrace $\nat:\Ome\Tc\Ac\to \Om_1\Cc$ (voir \cite{P4, Q2}), \`a la sortie par une supertrace partielle $\tau:\ell^1(H_{\Om^*\Rc})\to \Om^*\Rc$ afin d'obtenir une application lin\'eaire born\'ee
\be
\chi = \tau \exp(-\nabla^2) \nat\ :\ \Ome\Tc\Ac \to X(\Rc)\ . \label{map}
\ee
La normalisation de $\tau$ est fix\'ee de mani\`ere unique par p\'eriodicit\'e de Bott \cite{P4}, et sa parit\'e est la m\^eme que celle de $(H,\rho,D)$. L'identit\'e de Bianchi $[\nabla, \nabla^2]=0$ implique imm\'ediatement que (\ref{map}) est un \emph{morphisme} du $(b+B)$-complexe des formes diff\'erentielles enti\`eres sur $\Tc\Ac$, vers le $X$-complexe de $\Rc$. Comme on sait que $\Ome\Tc\Ac$ calcule l'homologie cyclique enti\`ere de $\Ac$ via l'\'equivalence de Goodwillie $\gamma: X(\Tc\Ac)\to \Ome\Tc\Ac$, on en d\'eduit que la compos\'ee
\be
\ch(H,\rho,D)\ :\ X(\Tc\Ac) \stackrel{\gamma}{\longrightarrow} \Ome\Tc\Ac \stackrel{\chi}{\longrightarrow} X(\Rc)
\ee
est un cocyle cyclique bivariant entier, dont la classe de cohomologie d\'efinit le caract\`ere de Chern. Soit maintenant $\cinf[0,1]$ l'alg\`ebre de Fr\'echet des fonctions lisses sur l'intervalle $[0,1]$ dont toutes les d\'eriv\'ees d'ordre $\geq 1$ s'annulent aux extr\'emit\'es, et notons $\Bc[0,1]$ le produit tensoriel compl\'et\'e $\Bc\hotimes\cinf[0,1]$. Un $\Ac$-$\Bc[0,1]$-bimodule d\'efinit une \emph{homotopie diff\'erentiable} entre les deux $\Ac$-$\Bc$-bimodules associ\'es aux extr\'emit\'es de l'intervalle. On d\'efinit de fa\c{c}on analogue une homotopie diff\'erentiable entre bimodules $\te$-sommables relativement \`a une extension $\Rc$. 

\begin{theorem}[\cite{P4}]
Soit $(H,\rho,D)$ un $\Ac$-$\Bc$-bimodule non-born\'e de parit\'e $i\in \zz_2$, et $\te$-sommable relativement \`a une extension analytique universelle $0\to \Jc\to\Rc\to\Bc\to 0$. Alors la classe de cohomologie cyclique enti\`ere bivariante du caract\`ere de Chern
\be
\ch(H,\rho,D)  \in HE^i(\Ac,\Bc)
\ee
est invariante par homotopie diff\'erentiable de bimodules $\te$-sommables. \hfill \cqfd
\end{theorem}

\begin{remark}\textup{La classe de cohomologie cyclique du caract\`ere de Chern d\'epend \emph{a priori} du choix de l'extension $\Rc$ et du rel\`evement de l'op\'erateur de Dirac $\Dh$, \`a moins que l'on puisse montrer que deux tels rel\`evements sont connect\'es par une homotopie de bimodules $\te$-sommables. Il convient donc de montrer, dans chaque situation concr\`ete, que l'on a effectivement construit le ``bon'' caract\`ere de Chern!}
\end{remark}

\begin{example} \textup{Lorsque $H$ est un espace de Hilbert, $\Ac$ une alg\`ebre quelconque, $\rho:\Ac\to \End(H)$ une repr\'esentation et $D$ un op\'erateur non born\'e autoadjoint \`a r\'esolvante compacte sur $H$, on obtient un triplet spectral \cite{C1}. Dans ce cas $\Bc=\cc$ est une alg\`ebre quasi-libre et il suffit de prendre l'extension triviale $\Rc=\cc$. La condition de $\te$-sommabilit\'e impose que l'op\'erateur de la chaleur $\exp(-D^2)$ soit tra\c{c}able, et le morphisme (\ref{map}) se r\'eduit \`a un cocycle cyclique entier $\chi: \Ome\Ac \to \cc$ qui correspond exactement au cocycle JLO bien connu \cite{JLO}.}
\end{example}

\begin{example}\label{ebott}\textup{Dans certaines situations il n'est pas n\'ecessaire de choisir une extension $\Rc$ analytique universelle, pourvu que le complexe $X(\Rc)$ porte suffisamment d'information sur l'homologie cyclique de $\Bc$. Illustrons cela dans le cas de la classe de Bott sur l'espace euclidien $\rr^n$. On prend pour $\Bc$ l'alg\`ebre de Fr\'echet commutative des fonctions lisses \`a d\'ecroissance rapide sur $\rr^n$, munie de sa bornologie de von Neumann. Soit $\Om^+(\rr^n)$ l'espace de Fr\'echet des formes diff\'erentielles lisses \`a d\'ecroissance rapide, de degr\'e pair, sur $\rr^n$. On d\'eforme le produit commutatif sur $\Om^+(\rr^n)$ en un produit non-commutatif introduit par Fedosov \cite{F}:
$$
\om_1\circ\om_2 = \om_1\om_2 - d\om_1 d\om_2\qquad \forall \om_i\in \Om^+(\rr^n)\ .
$$
Soit $\Rc$ cette alg\`ebre non-commutative, munie de sa bornologie de von Neumann. La projection de $\Rc$ sur l'espace des zero-formes $\Bc$ est un homomorphisme d'alg\`ebre, d'o\`u une extension
$$
0 \to \Jc \to \Rc \to \Bc \to 0\ .
$$
Le noyau $\Jc$ s'identifie \`a l'id\'eal nilpotent des formes diff\'erentielles paires de degr\'e $\geq 2$. $\Rc$ n'est pas quasi-libre. Cependant, le complexe $X(\Rc)$ contient toute l'information sur la cohomologie de de Rham \`a d\'ecroissance rapide sur $\rr^n$.\\
La classe de Bott est l'\'el\'ement de $K$-th\'eorie topologique de $\rr^n$ repr\'esent\'ee par le $\cc$-$\Bc$-bimodule $(H,\rho,D)$ suivant: $H$ est l'espace de dimension finie $S_n$ des spineurs complexes associ\'es \`a l'espace euclidien $\rr^n$; le $\Bc$-module $H_{\Bc}=S^n\otimes\Bc$ s'identifie aux sections de spineurs lisses et \`a d\'ecroissance rapide sur $\rr^n$; $\rho$ est la repr\'esentation triviale de $\cc$ par multiplication sur $H_{\Bc}$; et $D$ est la multiplication de Clifford du vecteur issu de l'origine de $\rr^n$ sur $H_{\Bc}$ (pour $n$ impair on doit prendre $H=S^n\otimes C_1$ et tensoriser $D$ par $\eps$). Puisque $\Ac=\cc$, le caract\`ere de Chern de $(H,\rho,D)$ se r\'eduit \`a une classe d'homologie dans $X(\Rc)$, et par cons\'equent \`a une classe de cohomologie de de Rham \`a d\'ecroissance rapide sur $\rr^n$. Le calcul explicite r\'ealis\'e dans \cite{P4} fait intervenir l'exponentielle $\exp(-\Dh^2)\in \Rc$ du laplacien de Dirac relev\'e au $\Rc$-module $H_{\Rc}$. Le caract\`ere de Chern est alors repr\'esent\'e par une forme diff\'erentielle d'allure gaussienne et de degr\'e maximal sur $\rr^n$. }
\end{example}

\section{Bimodules born\'es $p$-sommables}\label{sborn}

Soit $H$ un espace bornologique complet $\zz_2$-gradu\'e, $\End(H)$ l'alg\`ebre de ses endomorphismes born\'es et $\ell^1(H)$ le $\End(H)$-bimodule des endomorphismes tra\c{c}ables. Par souci de simplification on supposera que l'homomorphisme naturel $\ell^1(H)\to \End(H)$ est injectif, ce qui identifie $\ell^1(H)$ \`a un id\'eal bilat\`ere de $\End(H)$. On dit qu'une sous-alg\`ebre $\zz_2$-gradu\'ee $\Ic^s\subset \End(H)$ (compl\`ete dans sa propre bornologie) est \emph{$p$-sommable}, avec $p$ entier, si la puissance $p$-i\`eme de $\Ic^s$ d\'efinie comme l'image du produit de concat\'enation
$$
\underbrace{\Ic^s\hotimes \ldots \hotimes \Ic^s}_p \to \Ic^s\ ,
$$
est contenue dans $\ell^1(H)$. Un exemple bien connu est celui des classes de Schatten $\Ic^s=\ell^p(H)$ sur un espace de Hilbert $H$. \\

Soient $\Ac$ et $\Bc$ deux alg\`ebres bornologiques compl\`etes, $H$ un espace bornologique complet $\zz_2$-gradu\'e et $\Ic^s\subset\End(H)$ une sous-alg\`ebre $p$-sommable pour un entier $p\geq 1$ fix\'e. L'alg\`ebre $\Ic^s\hotimes\Bc$ agit de mani\`ere \'evidente par endomorphismes sur le $\Bc$-module \`a droite $H_{\Bc}=H\hotimes\Bc$. Tout triplet $(H,\rho,F)$ est appel\'e \emph{$\Ac$-$\Bc$-bimodule born\'e $p$-sommable} s'il v\'erifie les propri\'et\'es suivantes:
\begin{itemize}
\item Il existe une sous-alg\`ebre $\zz_2$-gradu\'ee $\Ec^s\subset \End(H_{\Bc})$, compl\`ete dans sa propre bornologie, contenant $\Ic^s\hotimes\Bc$ comme id\'eal bilat\`ere. On \'ecrira 
$$
\End(H_{\Bc}) \supset \Ec^s \triangleright \Ic^s\hotimes\Bc\ .
$$
\item  $\rho:\Ac\to \Ec^s$ est une repr\'esentation born\'ee de $\Ac$ dans la sous-alg\`ebre de degr\'e pair de $\Ec^s$.
\item $F\in\End(H_{\Bc})$ est un endomorphisme de degr\'e impair multiplicateur de $\Ec^s$ et tel que $F^2=1$.
\item $[F,\rho(a)]\in \Ic^s\hotimes\Bc$ pour tout $a\in\Ac$.
\end{itemize}
Les bimodules born\'es \emph{pairs} ou \emph{impairs} sont d\'efinis de mani\`ere analogue aux bimodules non-born\'es. La n\'ecessit\'e de consid\'erer une sous-alg\`ebre $\Ec^s\subset\End(H_{\Bc})$ est dict\'ee par le fait que $\Ic^s\hotimes\Bc$ n'est pas n\'ecessairement un id\'eal bilat\`ere de $\End(H_{\Bc})$. On peut penser \`a $\Ec^s$ comme \'etant la plus petite alg\`ebre d'endomorphismes stable sous multiplication par $F$ et qui contient l'image de $\rho$. \\

Pour construire le caract\`ere de Chern de $(H,\rho,F)$ en cohomologie cyclique bivariante, nous allons relever comme pr\'ec\'edemment le bimodule \`a des extensions universelles de $\Ac$ et $\Bc$. Choisissons donc une extension quasi-libre analytiquement nilpotente $0 \to \Jc \to \Rc \to \Bc \to 0$ et supposons dans un premier temps que l'op\'erateur $F\in \End(H_{\Bc})$ soit de la forme
\be
F=G\otimes 1 \quad \mbox{avec}\quad G^2=1\in \End(H)\ .
\ee
Alors $F$ se rel\`eve canoniquement en un endomorphisme $\Fh=F$ sur le $\Rc$-module \`a droite $H_{\Rc}$. La condition d'admissibilit\'e suivante est adapt\'ee de \cite{P8}:
 
\begin{definition}\label{dadm}
Le bimodule $(H,\rho,F)$ est \emph{admissible} relativement \`a l'extension $0 \to \Jc \to \Rc \to \Bc \to 0$ s'il existe deux sous-alg\`ebres $\zz_2$-gradu\'ees 
$$
\End(H_{\Rc})\supset \Mc^s\triangleright\Ic^s\hotimes\Rc\ ,\qquad \End(H_{\Rc})\supset \Nc^s\triangleright\Ic^s\hotimes\Jc
$$
avec $F$ multiplicateur de $\Mc^s$, ainsi qu'un diagramme commutatif d'extensions  
$$
\vcenter{\xymatrix{
0\ar[r] & \Nc^s \ar[r]  & \Mc^s \ar[r] & \Ec^s \ar[r] & 0 \\
0 \ar[r] & \Ic^s\hotimes \Jc \ar[r] \ar[u] & \Ic^s\hotimes \Rc \ar[r] \ar[u] & \Ic^s\hotimes\Bc \ar[r] \ar[u] & 0}}
$$
\end{definition}

\begin{example}\textup{Les bimodules consid\'er\'es dans \cite{P5} sont admissibles par rapport \`a n'importe quelle extension $\Rc$ et ont la forme suivante: $H$ est un espace de Hilbert et
$$
\Ic^s= \ell^p(H)\ ,\qquad \Ec^s= \End(H)\hotimes\Bc \ .
$$
En effet on peut alors prendre $\Mc^s=\End(H)\hotimes \Rc$ et $\Nc^s=\End(H)\hotimes \Jc$. Le fait d'imposer que l'image de $\rho$ appartienne au produit tensoriel $\End(H)\hotimes\Bc\subset \End(H_{\Bc})$ est assez restrictif. Il existe n\'eanmoins un certain nombre d'exemples int\'eressants v\'erifiant cette propri\'et\'e, comme celui du chapitre \ref{cgrou}. Cette situation ne couvre cependant pas tous les cas de figure importants, en particulier le repr\'esentant born\'e de la classe de Bott sur l'espace $\rr^n$ v\'erifie seulement la condition plus g\'en\'erale \ref{dadm}. }
\end{example}

La propri\'et\'e universelle de l'alg\`ebre tensorielle analytique permet ensuite de relever l'homomorphisme $\rho:\Ac\to \Ec^s$ en un homomorphisme $\rho_*:\Tc\Ac\to \Mc^s$ en vertu du diagramme commutatif
$$
\vcenter{\xymatrix{
0 \ar[r] & \Jc\Ac \ar[r] \ar[d]_{\rho_*} & \Tc\Ac \ar[r] \ar[d]_{\rho_*} & \Ac \ar[r] \ar[d]^{\rho} & 0 \\
0\ar[r] & \Nc^s \ar[r]  & \Mc^s \ar[r] & \Ec^s \ar[r] \ar@/_/@{.>}[l]_{\si} & 0 }}
$$
Notons que pour $x\in\Tc\Ac$, le commutateur $[F,x]\in \Mc^s$ est dans l'id\'eal $\Ic^s\hotimes\Rc$ et par cons\'equent le triplet $(H,\rho_*,F)$ d\'efinit un $\Tc\Ac$-$\Rc$-bimodule born\'e $p$-sommable, de m\^eme parit\'e $i\in \zz_2$ que le bimodule initial $(H,\rho,F)$. Son caract\`ere de Chern dans $HE^i(\Ac,\Bc)$ est repr\'esent\'e, pour tout choix d'entier $n\geq p$ de parit\'e $i\mod 2$, par un cocycle $\chih^n\in \hom(\Ome\Tc\Ac,X(\Rc))$ construit de la fa\c{c}on suivante \cite{P5}. $\chih^n$ s'annule sur les espaces $\Om^k\Tc\Ac$ si $k\neq n$ et $n+1$, et ses deux composantes non nulles $\chih^n_0:\Om^n\Tc\Ac\to \Rc$ et $\chih^n_1:\Om^{n+1}\Tc\Ac\to \Om^1\Rc_{\nat}$ sont d\'efinies par
\beq
&& \chih^n_0(x_0\dd x_1\ldots \dd x_n) = (-)^n\frac{\Gamma(1+\frac{n}{2})}{(n+1)!}\,\sum_{\la\in S_{n+1}}\eps(\la)\tau(x_{\la(0)}[F,x_{\la(1)}]\ldots[F,x_{\la(n)}])\ ,\non\\
&&\chih^n_1(x_0\dd x_1\ldots \dd x_{n+1}) = (-)^n\frac{\Gamma(1+\frac{n}{2})}{(n+1)!}\, \sum_{i=1}^{n+1}\tau\nat ( x_0[F,x_{1}]\ldots \dd x_i \ldots [F,x_{n+1}])\ ,\non\\
&& \label{chi0}
\eeq
avec $S_{n+1}$ le groupe des permutations cycliques sur $n+1$ \'el\'ements, $\eps(\la)$ la signature de la permutation $\la$, et $\tau$ la supertrace partielle provenant de la supertrace d'op\'erateurs sur la puissance $n$-i\`eme de $\Ic^s$. Ainsi la compos\'ee
\be
\ch^n(H,\rho,F)\ :\ X(\Tc\Ac) \stackrel{\gamma}{\longrightarrow} \Ome\Tc\Ac \stackrel{\chi^n}{\longrightarrow} X(\Rc)
\ee
est un cocycle cyclique bivariant entier pour tout $n\geq p$. Une homotopie diff\'erentiable entre deux bimodules born\'es est d\'efinie de mani\`ere analogue au cas des bimodules non-born\'es.

\begin{theorem}[\cite{P5}]
Soit $(H,\rho,F)$ un $\Ac$-$\Bc$-bimodule born\'e $p$-sommable et de parit\'e $i\in \zz_2$, avec $F=G\otimes 1$. On le suppose admissible relativement \`a une extension universelle de $\Bc$. Alors pour tout entier $n\geq p$, $n\equiv i\mod 2$, la classe de cohomologie cyclique bivariante enti\`ere du caract\`ere de Chern
\be
\ch^n(H,\rho,F) \in HE^i(\Ac,\Bc)
\ee
est invariante sous homotopie diff\'erentiable et ind\'ependante du choix de $n$. \cqfd
\end{theorem}

\begin{remark}\textup{Comme dans le cas non-born\'e, la classe de cohomologie cyclique du caract\`ere de Chern d\'epend en principe du choix des extensions $0 \to \Jc \to \Rc \to \Bc \to 0$ et $0 \to \Nc^s \to \Mc^s \to \Ec^s \to 0$.}
\end{remark}

Dans la situation o\`u $F$ n'est pas de la forme $G\otimes 1$, on peut toujours tenter de prendre un rel\`evement $\Fh\in \End(H_{\Rc})$ mais alors $\Fh^2\neq 1$ en g\'en\'eral. Nous avons montr\'e dans \cite{P5} comment modifier en cons\'equence les formules (\ref{chi0}). Les cocycles qui en r\'esultent sont alors plus difficiles \`a g\'erer et ne d\'efinissent pas de fa\c{c}on \'evidente une classe de cohomologie cyclique dans $HE^*(\Ac,\Bc)$. Il convient de remarquer que la condition $F=G\otimes 1$ n'est pas v\'eritablement restrictive puisqu'il s'agit de la repr\'esentation standard d'un \'el\'ement de $K$-th\'eorie bivariante sous la forme d'un \emph{quasihomomorphisme} \cite{Cu}
\be
\rho: \Ac\to \Ec^s\triangleright \Ic^s\hotimes\Bc\ . 
\ee
Nous reviendrons sur les formules (\ref{chi0}) dans le chapitre \ref{csec} d\'edi\'e aux images directes d'invariants primaires et secondaires pour les $m$-alg\`ebres de Fr\'echet. Notons que dans \cite{Ni1, Ni2} Nistor a construit un caract\`ere de Chern bivariant pour les quasihomomorphismes $p$-sommables. J'ignore s'il co\"{\i}ncide exactement avec celui consid\'er\'e ici. Cependant nous pr\'ef\'ererons utiliser les formules (\ref{chi0}) en raison de leur compatibilit\'e avec le caract\`ere de Chern des bimodules non-born\'es $\te$-sommables:\\

\noindent {\bf Lien entre bimodules born\'es et non-born\'es.} Consid\'erons maintenant un $\Ac$-$\Bc$-bimodule born\'e $p$-sommable $(H,\rho,F)$ v\'erifiant les hypoth\`eses du th\'eor\`eme ci-dessus. Soit $|D|$ un endomorphisme non born\'e de degr\'e pair sur $H_{\Bc}$ commutant avec $F$, et tel que $(H,\rho,D)$ avec $D=|D|F$ soit un bimodule non-born\'e $\te$-sommable. Dans \cite{P5} nous donnons des conditions formelles pour assurer l'\'egalit\'e des caract\`eres de Chern $\ch(H,\rho,D)\equiv \ch^n(H,\rho,F)$ dans $HE^i(\Ac,\Bc)$. Il s'agit d'une g\'en\'eralisation bivariante du proc\'ed\'e de r\'etraction introduit par Connes et Moscovici pour les triplets spectraux \cite{CM93}. Le principe est bas\'e sur des transgressions de Chern-Simons dans le complexe $\hom(\Ome\Tc\Ac,X(\Rc))$, suivies d'une homotopie entre $D$ et $F$ (dans un sens \`a pr\'eciser)
$$
D_t = D/|D|^t\ ,\quad t\in [0,1]\ .
$$
C'est d'ailleurs en suivant ce proc\'ed\'e que nous avons \'etabli dans \cite{P5} les formules (\ref{chi0}). Si l'on se place dans le cadre bornologique g\'en\'eral cette r\'etraction est purement formelle. Elle doit donc \^etre utilis\'ee au cas par cas, dans les situations concr\`etes o\`u tout est bien d\'efini. La comparaison des caract\`eres de Chern $\ch(H,\rho,D)$ et $\ch^n(H,\rho,F)$ fournit alors un outil efficace pour \'etablir des th\'eor\`emes de l'indice en g\'eom\'etrie non-commutative. Nous proposons dans le chapitre suivant une illustration de ces m\'ethodes par l'\'etude des actions propres et isom\'etriques de groupes localement compacts sur des vari\'et\'es de Riemann.

\chapter{Th\'eor\`eme de l'indice \'equivariant}\label{ceq}

Ce chapitre sert d'illustration aux formules de caract\`ere de Chern bivariant introduites pr\'ec\'edemment. On consid\`ere un groupe $G$ localement compact agissant proprement par isom\'etries sur une vari\'et\'e riemannienne $M$ compl\`ete et cocompacte. L'alg\`ebre de convolution des fonctions continues \`a support compact sur $G$ est compl\'et\'ee en une alg\`ebre de Banach ``admissible'' $\Bc$ (voir la d\'efinition \ref{dad}). A tout op\'erateur diff\'erentiel elliptique $G$-invariant $Q$ d'ordre 1 sur $M$ on peut associer un indice qui est une classe de $K$-th\'eorie topologique $\mu(Q)\in \Kt_*(\Bc)$. L'objectif est alors de calculer son caract\`ere de Chern en homologie cyclique enti\`ere 
\be
\ch(\mu(Q))\in HE_*(\Bc)\ .\label{indice}
\ee
On d\'esigne par $\Ac$ le produit crois\'e $\cinfc(M)\cp G$. C'est une alg\`ebre bornologique compl\`ete dot\'ee d'une classe canonique en $K$-th\'eorie $[e]\in \Kt_0(\Ac)$. Le th\'eor\`eme de l'indice d\'emontr\'e dans \cite{P6} exprime le caract\`ere de Chern (\ref{indice}) comme cup-produit de $\ch(e)\in HE_0(\Ac)$ avec le caract\`ere de Chern bivariant d'un $\Ac$-$\Bc$-bimodule non-born\'e $\te$-sommable naturellement attach\'e \`a l'op\'erateur elliptique $Q$. Le d\'eveloppement asymptotique du noyau de la chaleur \`a temps court permet ainsi d'obtenir une \emph{formule locale} pour (\ref{indice}), faisant intervenir les sous-vari\'et\'es des points fixes de l'action de $G$ sur $M$. Ces r\'esultats sont issus de l'article\\

\noindent \cite{P6} D. Perrot: The equivariant index theorem in entire cyclic cohomology, preprint arXiv:math/0410315, \`a para\^{\i}tre dans {\it J. K-Theory} (disponible en ligne).

\section{Actions propres}\label{sad}

Soit $G$ un groupe topologique localement compact s\'eparable. On note $C_c(G)$ l'alg\`ebre des fonctions continues \`a valeur complexe et \`a support compact sur $G$, munie du produit de convolution
\be
(b_1b_2)(g) = \int_G dh\, b_1(h)b_2(gh^{-1})\ ,\quad \forall\ b_i\in C_c(G)\ ,\ g\in G\ ,
\ee
o\`u $dh$ est une mesure de Haar invariante \`a droite. Nous aurons besoin de la compl\'eter en une alg\`ebre de Banach:

\begin{definition}[\cite{P6}]\label{dad}
Soit $G$ un groupe localement compact et $dg$ une mesure de Haar invariante \`a droite. Une mesure $d\nu$ sur $G$ est \emph{admissible} s'il existe une fonction $\si$ strictement positive et continue sur $G$ telle que 
\be
d\nu= \si\, dg \quad \mbox{et} \quad \si(gh)\leq \si(g)\si(h)\quad \forall g,h\in G\ .
\ee
La norme $\|b\|=\int_Gd\nu\, |b(g)|$ associ\'ee \`a cette mesure v\'erifie alors $\|b_1b_2\|\leq \|b_1\|\|b_2\|$ pour tous $b_1,b_2\in C_c(G)$, et l'alg\`ebre de Banach $\Bc=L^1(G,d\nu)$ ainsi obtenue est appel\'ee une \emph{compl\'etion admissible} de l'alg\`ebre de convolution.
\end{definition}

On peut construire un bon exemple de mesure admissible \`a partir d'une distance $d:G\times G\to \rr_+$ invariante \`a droite sur $G$ et d'un param\`etre $\al\in \rr_+$. En tout point $g\in G$ posons
$$
\si(g)= (1+d(g,1))^{\al}\ .
$$
La fonction $\si$ cro\^{\i}t donc comme une puissance de la distance qui s\'epare $g$ de l'identit\'e. Les \'el\'ements de $\Bc=L^1(G,d\nu)$ sont des fonctions localement int\'egrables sur $G$ qui v\'erifient une certaine condition de d\'ecroissance \`a l'infini, suivant la valeur de $\al$. En particulier lorsque $G$ est ab\'elien, on voit par transformation de Fourier que $\Bc$ est un espace de fonctions sur le dual de Pontrjagin $\Gh$, d'autant plus ``diff\'erentiables'' que $\al$ est grand. Cette n\'ecessit\'e de contr\^oler le degr\'e de r\'egularit\'e des fonctions est dict\'ee par l'utilisation de l'homologie cyclique.\\
 
Consid\'erons maintenant une vari\'et\'e diff\'erentielle $M$ compl\`ete, lisse et sans bord, sur laquelle $G$ agit proprement par diff\'eomorphismes. On suppose de plus que le quotient $X=G\backslash M$ est compact. Puisque l'action est propre, on peut sans perte de g\'en\'eralit\'e fixer une m\'etrique de Riemann sur $M$ telle que $G$ agisse par isom\'etries. De m\^eme, si $E\to M$ est un espace fibr\'e vectoriel $G$-\'equivariant, complexe et de rang fini, on peut toujours le supposer muni d'une structure hermitienne $G$-invariante.\\
Soient $E_+\to M$ et $E_-\to M$ deux fibr\'es vectoriels complexes $G$-\'equivariants. Pour tout $m\in\rr$ on d\'esigne par $\Psi^m_c(E_+,E_-)$ l'espace des op\'erateurs pseudodiff\'erentiels $G$-invariants d'ordre $m$ et \`a \emph{support propre} \cite{Ho, Ka}. Le groupe de $K$-homologie $G$-\'equivariante de degr\'e pair $K_0^G(M)$ est d\'efini comme l'ensemble des classes d'homotopie stable d'op\'erateurs pseudodiff\'erentiels elliptiques $G$-invariants $Q\in \Psi^0_c(E_+,E_-)$. De m\^eme le groupe de $K$-homologie \'equivariante de degr\'e impair $K_1^G(M)$ est l'ensemble des classes d'homotopie stable d'op\'erateurs pseudodiff\'erentiels elliptiques $G$-invariants autoadjoints $Q\in \Psi^0_c(E,E)$. L'addition sur est induite par somme directe de fibr\'es et d'op\'erateurs. Un op\'erateur pseudodiff\'erentiel elliptique $Q$ d'ordre $m$ quelconque d\'etermine aussi un \'el\'ement de $K$-homologie $[Q]$ en prenant la classe de l'op\'erateur d'ordre z\'ero $Q\cdot\delta_m$, o\`u $\delta_m$ est un op\'erateur elliptique $G$-invariant de symbole $s(x,p)=(1+\|p\|)^{-m}$ et d'ordre $-m$.\\

Choisissons maintenant une compl\'etion admissible $\Bc$ de l'alg\`ebre de convolution $C_c(G)$. A tout op\'erateur elliptique $Q\in \Psi_c^0(E_+,E_-)$ repr\'esentant une classe de $K$-homologie paire on peut associer son indice dans la $K$-th\'eorie topologique d'alg\`ebre de Banach $\Bc$ de la mani\`ere suivante. Choisissons une param\'etrix $G$-invariante $P\in \Psi_c^0(E_-,E_+)$ pour $Q$:
$$
PQ - 1 \in \Psi_c^{-\infty}(E_+,E_+)\ ,\qquad QP - 1 \in \Psi_c^{-\infty}(E_-,E_-)\ .
$$
D\'esignons par $E$ le fibr\'e vectoriel $E_+\oplus E_-$ et consid\'erons l'op\'erateur pseudodiff\'erentiel inversible $T\in \Psi_c^0(E,E)$ 
$$
T= \left(\begin{matrix} 1-PQ & P \\ 2Q-QPQ & QP-1 \end{matrix}\right)\ .
$$
Soit $\Psit_c^{-\infty}(E,E)$ l'alg\`ebre des op\'erateurs pseudodiff\'erentiels r\'egularisants, dont on a rajout\'e une unit\'e aux blocs diagonaux. Les idempotents
\be
e_1 = T^{-1}\left(\begin{matrix} 1 & 0 \\ 0 & 0 \end{matrix}\right)T \ ,\qquad e_0 = \left(\begin{matrix} 0 & 0 \\ 0 & 1 \end{matrix}\right)
\ee
sont des \'el\'ements de $\Psit_c^{-\infty}(E,E)$ dont la diff\'erence appartient \` a $\Psi_c^{-\infty}(E,E)$. Soit maintenant $\Kc$ l'alg\`ebre de Fr\'echet des ``op\'erateurs compacts lisses'' (matrices complexes infinies \`a d\'ecroissance rapide \cite{Cu2}), et $C_c(G;\Kc)$ l'alg\`ebre de convolution des fonctions continues \`a support compact sur $G$ et \`a valeurs dans $\Kc$. A l'aide d'une fonction ``cut-off'' sur $M$, on construit un homomorphisme (voir \cite{C1, P6})
\be 
\te: \Psi_c^{-\infty}(E,E) \to C_c(G,\Kc)\ .
\ee
Or $C_c(G;\Kc)$ est une sous-alg\`ebre du produit tensoriel projectif $\Kc\hotimes\Bc$. L'indice analytique $G$-\'equivariant de l'op\'erateur $Q$ est alors d\'efini comme la classe de $K$-th\'eorie topologique 
\be
\mu(Q) = [\te(e_1)] - [\te(e_0)] \in \Kt_0(\Kc\hotimes\Bc)\cong \Kt_0(\Bc)\ ,
\ee
qui est ind\'ependante des choix effectu\'es. Lorsque $Q$ repr\'esente une classe de $K$-homologie de degr\'e impair, on d\'efinit son indice $\mu(Q)\in K_1(\Bc)$ en se ramenant au cas pair par p\'eriodicit\'e de Bott. L'indice analytique descend en une application sur la $K$-homologie \'equivariante
\be
\mu:K_i^G(M)\to \Kt_i(\Bc)\qquad  i\in\zz_2\ .
\ee 
Notons qu'\`a ce niveau le choix d'une compl\'etion de $C_c(G)$ en alg\`ebre de Banach importe peu. La propri\'et\'e d'``admissibilit\'e'' ne sera vraiment utilis\'ee que dans la construction du caract\`ere de Chern bivariant.

\section{Le bimodule associ\'e}

Comme pr\'ec\'edemment consid\'erons un groupe localement compact $G$ et une vari\'et\'e riemannienne compl\`ete $M$ munie d'une $G$-action propre, cocompacte et isom\'etrique. Pour toute partie compacte $K\subset M$ d\'esignons par $\cinf_K(M)$ l'espace de Fr\'echet des fonctions lisses sur $M$ \`a support dans $K$, et par $C_c(G;\cinf_K(M))$ l'espace des fonctions continues \`a support compact sur $G$ et \`a valeurs dans $\cinf_K(M)$. La limite inductive
\be
\Ac = \varinjlim_{K\subset M} C_c(G;\cinf_K(M))
\ee
est un espace bornologique complet qui s'identifie \`a un espace de fonctions sur le groupo\"{\i}de $G\ltimes M$. On d\'efinit le produit de convolution de deux \'el\'ements $a_1,a_2\in\Ac$ par
$$
(a_1a_2)(g,x) = \int_G dh\, a_1(h,x)a_2(gh^{-1},hx)\ ,\quad \forall (g,x)\in G\times M\ ,
$$
o\`u $dh$ est la mesure de Haar invariante \`a droite. Alors $\Ac$ est une alg\`ebre bornologique compl\`ete \cite{P6}, que l'on notera sous la forme d'un produit crois\'e
\be
\Ac = \cinfc(M)\cp G\ ,
\ee
avec $\cinfc(M)$ l'alg\`ebre commutative des fonctions lisses \`a support compact sur $M$.\\

Choisissons maintenant une compl\'etion admissible $\Bc$ de l'alg\`ebre de convolution $C_c(G)$, d'apr\`es la d\'efinition \ref{dad}. Munie de sa bornologie de von Neumann, $\Bc$ est une alg\`ebre bornologique compl\`ete. Nous allons associer un $\Ac$-$\Bc$-bimodule non-born\'e \`a tout op\'erateur \emph{diff\'erentiel} elliptique $Q$ d'ordre 1 repr\'esentant une classe de $K$-homologie \'equivariante $[Q]\in K_*^G(M)$. Dans le cas pair, $Q\in \Psi_c^1(E_+,E_-)$ agit entre les sections de deux fibr\'es vectoriels hermitiens $G$-\'equivariants. Avec son adjoint formel $Q^*\in \Psi_c^1(E_-,E_+)$ on peut construire un op\'erateur diff\'erentiel elliptique de degr\'e impair agissant sur les sections du fibr\'e $\zz_2$-gradu\'e $E=E_+\oplus E_-$
\be
D= \left(\begin{matrix} 0 & Q^* \\ Q & 0 \end{matrix}\right) \ .
\ee
$D$ s'\'etend en un op\'erateur autoadjoint non born\'e sur l'espace de Hilbert $\zz_2$-gradu\'e $H=L^2(E)$ associ\'e \`a la m\'etrique de Riemann sur $M$ et la structure hermitienne sur $E$. Dans le cas impair, $Q\in \Psi_c^1(E,E)$ est un op\'erateur autoadjoint agissant sur les sections d'un fibr\'e $E$ trivialement gradu\'e; en tensorisant l'espace des sections de $E$ par l'alg\`ebre de Clifford $\zz_2$-gradu\'ee $C_1=\cc\oplus \cc\eps$, on construit l'op\'erateur diff\'erentiel elliptique de degr\'e impair
\be
D= Q\otimes \eps\ .
\ee
$D$ s'\'etend en un op\'erateur autoadjoint non born\'e sur l'espace de Hilbert $\zz_2$-gradu\'e $H=L^2(E)\otimes C_1$. On tra\^{\i}te maintenant les cas pair et impair simultan\'ement. Munissons $H$ de sa bornologie de von Neumann. L'alg\`ebre des fonctions lisses \`a support compact $\cinfc(M)$ agit par multiplication sur les sections de $E$, donc est repr\'esent\'ee par endomorphismes born\'es sur $H$. De plus l'action de $G$ \'etant isom\'etrique, elle induit une repr\'esentation unitaire $r:G\to U(H)$. On obtient un $\Ac$-$\Bc$-bimodule non-born\'e $(H,\rho,D)$ de m\^eme parit\'e que la classe $[Q]$ comme suit:
\begin{itemize}
\item $H_{\Bc}=H\hotimes\Bc$ est isomorphe \`a l'espace de Banach $L^1(G;H)$ des fonctions int\'egrables sur $G$ (relativement \`a la mesure admissible $d\nu$) et \`a valeurs dans $H$. Sa structure de $\Bc$-module \`a droite est donn\'ee par 
$$
(\xi b)(g)=\int_G dh\, \xi(h)b(gh^{-1})\ ,\quad \forall\ \xi\in H_{\Bc}\ ,\ b\in \Bc\ ,\ g\in G\ .
$$
\item $\rho:\Ac\to\End(H_{\Bc})$ est la repr\'esentation de $\Ac$ par endomorphismes pairs  
$$
(\rho(a)\xi)(g)=\int_G dh\, a(h)r(h)\cdot \xi(gh^{-1})\ ,\quad \forall\ a\in\Ac\ ,\ \xi\in H_{\Bc}\ ,\ g\in G\ ,
$$
o\`u l'\'evaluation de $a\in\Ac$ en un point $h\in G$ donne une fonction $a(h)\in\cinfc(M)$ vue comme op\'erateur born\'e sur $H$. 
\item $D: \dom(D)\subset H_{\Bc}\to H_{\Bc}$ est l'op\'erateur non born\'e impair de domaine dense 
$$
(D\xi)(g)=D(\xi(g))\ ,\quad \forall\ \xi\in \dom(D)\ ,\ g\in G\ ,
$$
qui commute avec l'action \`a droite de $\Bc$. Comme $D$ est un op\'erateur diff\'erentiel d'ordre 1, le commutateur $[D,\rho(a)]\in \End(H_{\Bc})$ s'\'etend en un endomorphisme born\'e pour tout $a\in\Ac$.
\end{itemize}
L'op\'erateur de la chaleur $\exp(-tD^2)\in \End(H)$, d\'efini par calcul fonctionnel sur l'extension autoadjointe de $D$, est aussi naturellement un endomorphisme du $\Bc$-module $H_{\Bc}$. Nous avons montr\'e dans \cite{P6} que le bimodule non born\'e $(H,\rho,D)$ ainsi obtenu a les propri\'et\'es requises pour construire son caract\`ere de Chern bivariant
\be
\ch(H,\rho,D)\in HE^*(\Ac,\Bc)\ .
\ee
On choisit ici l'extension universelle $\Rc=\Tc\Bc$ pour $\Bc$. L'analyse est grandement simplifi\'ee par le fait que l'op\'erateur de Dirac sur $H_{\Bc}$ provient d'un op\'erateur sur $H$. Son rel\`evement au $\Rc$-module $H_{\Rc}$ se r\'eduit donc \`a $\Dh=D$. Les propri\'et\'es de la compl\'etion admissible $\Bc$ sont essentielles dans la preuve de la proposition suivante:

\begin{proposition}[\cite{P6}]\label{pana}
Soit $Q$ un op\'erateur diff\'erentiel elliptique $G$-invariant d'ordre 1 repr\'esentant une classe de $K$-homologie $[Q]\in K_i^G(M)$, $i\in\zz_2$, et soit $D$ l'op\'erateur de Dirac associ\'e. Alors pour tout $t>0$, le bimodule $(\Ec,\rho,\sqrt{t}D)$ est $\te$-sommable. Son caract\`ere de Chern
\be 
\ch(H,\rho,\sqrt{t}D) \in HE^i(\Ac,\Bc)
\ee
est une classe de cohomologie cyclique bivariante enti\`ere ind\'ependante de $t$. \cqfd 
\end{proposition}
Les formules de d\'eveloppement asymptotique de l'op\'erateur de la chaleur $\exp(-tD^2)$ \`a la limite $t\downarrow 0$ permettront d'obtenir une formule \emph{locale} pour le caract\`ere de Chern bivariant.

\section{Th\'eor\`eme de l'indice}

Lorsque $M$ est munie d'une action propre et cocompacte de $G$, il existe une fonction cut-off $c\in \cinfc(M)$, c'est-\`a-dire telle que $\int_G c(gx)^2 dg=1$ pour tout $x\in M$. On peut alors construire un idempotent $e\in \Ac=\cinfc(M)\cp G$ en posant
\be
e(g,x) = c(x)c(gx)\qquad \forall\ (g,x)\in G\times M\ .
\ee
Cet idempotent d\'efinit une classe canonique $[e]\in \Kt_0(\Ac)$ dans la $K$-th\'eorie de l'alg\`ebre bornologique $\Ac$, ind\'ependante du choix de fonction cut-off. Son caract\`ere de Chern est une classe d'homologie cyclique enti\`ere de degr\'e pair
$$
\ch(e)\in HE_0(\Ac)\ .
$$
Choisissons une compl\'etion admissible $\Bc$ de l'alg\`ebre de convolution $C_c(G)$. Pour toute classe de cohomologie cyclique enti\`ere bivariante $\varphi\in HE^i(\Ac,\Bc)$, le cup-produit $\varphi\cdot\ch(e)$ d\'efinit une classe dans l'homologie cyclique enti\`ere $HE_i(\Bc)$. On a alors le th\'eor\`eme de l'indice suivant. 

\begin{theorem}[\cite{P6}]\label{tind}
Soit $[Q]\in K_i^G(M)$, $i\in\zz_2$, une classe de $K$-homologie \'equivariante repr\'esent\'ee par un op\'erateur diff\'erentiel elliptique $G$-invariant d'ordre 1, et soit $(H,\rho,D)$ le $\Ac$-$\Bc$-bimodule non-born\'e associ\'e. Alors l'image de $[Q]$ par l'application d'assemblage $\mu:K_i^G(M) \to K_i(\Bc)$ a un caract\`ere de Chern donn\'e par le cup-produit 
\be
\ch(\mu(Q))= \ch(H,\rho,D)\cdot \ch(e)\ \in HE_i(\Bc)\ ,
\ee
o\`u $\ch(H,\rho,D)\in HE^i(\Ac,\Bc)$ est le caract\`ere de Chern du bimodule et $\ch(e)\in HE_0(\Ac)$ le caract\`ere de Chern de la classe canonique $[e]\in \Kt_0(\Ac)$. 
\end{theorem}
\emph{Esquisse de d\'emonstration:} Premi\`erement on peut se ramener \`a un op\'erateur de Dirac inversible. Sa phase $F=D/|D|$, $F^2=1$, est un op\'erateur born\'e sur $H$. Comme expliqu\'e au chapitre \ref{ccar}, on r\'etracte ensuite le bimodule non born\'e $(H,\rho,D)$ sur le bimodule born\'e $(H,\rho,F)$ au moyen de l'homotopie op\'erateur 
$$
D_t = D/|D|^t\ ,\quad t\in [0,1]\ .
$$
$(H,\rho,F)$ est $p$-sommable pour tout $p> \dim M$. Les cocycles cycliques bivariants entiers $\ch^n(H,\rho,F)$ donn\'es par les formules (\ref{chi0}) d\'efinissent la m\^eme classe dans $HE^i(\Ac,\Bc)$ que le caract\`ere de Chern $\ch(H,\rho,D)$. Dans le cas $i=0$, le cup-produit de $\ch^n(H,\rho,F)$ avec $\ch(e)$ est une formule simple qui s'identifie au caract\`ere de Chern d'un idempotent $p$-sommable. On se ram\`ene ensuite \`a l'idempotent $\mu(Q)$ par d\'eformation. Le cas $i=1$ s'en d\'eduit par p\'eriodicit\'e de Bott. \cqfd \\

Soit $t>0$. Le cup-produit $\ch(H,\rho,\sqrt{t} D)\cdot\ch(e)$ repr\'esente $\ch(\mu(Q))$ en vertu de la proposition \ref{pana}. C'est une cha\^{\i}ne enti\`ere sur l'alg\`ebre $\Bc$. D\'esignons par $\ch_n(\sqrt{t}D)\in\Om^n\Bc$ sa composante de degr\'e $n$. Comme $\Bc=L^1(G,d\nu)$ est un espace de fonctions int\'egrables sur $G$, on voit que
$$
\Om^n\Bc = \Bc^{\hotimes n} \oplus \Bc^{\hotimes (n+1)} = L^1(G^n) \oplus L^1(G^{n+1})
$$
est un espace de fonctions int\'egrables sur l'ensemble $G^n\cup G^{n+1}$ muni de la mesure produit. On peut donc regarder $\ch_n(\sqrt{t}D)$ comme une fonction sur $G^n\cup G^{n+1}$. En fait, elle est continue et \`a support compact tant que $t>0$. Lorsque le fibr\'e $E\to M$ est un module de Clifford et $D$ un op\'erateur de Dirac g\'en\'eralis\'e, nous avons calcul\'e dans \cite{P6} que l'\'evaluation de $\ch_n(\sqrt{t}D)$ en un point $\gt\in G^n\cup G^{n+1}$ est donn\'ee par une formule locale explicite \`a la limite $t\downarrow 0$, faisant intervenir les \emph{points fixes} de l'action de $G$ sur $M$. Premi\`erement, \`a partir de l'idempotent $e\in \Ac$ on construit de mani\`ere purement alg\'ebrique une forme diff\'erentielle mixte
\be
\ch_n(e) \in \Omc^*(M)\otimes\Om^n\Bc\ .
\ee
Elle peut se voir comme une fonction sur $G^n\cup G^{n+1}$ \`a valeurs dans l'espace $\Omc^*(M)$ des formes diff\'erentielles (commutatives) \`a support compact sur $M$. Ensuite, pour tout $g\in G$ d\'esignons par $M_g\subset M$ l'ensemble des points fixes de $g$. C'est une r\'eunion de sous-vari\'et\'es de $M$, qui peuvent avoir diff\'erentes dimensions. Au moins localement, le module de Clifford $E$ restreint \`a $M_g$ se d\'ecompose en un produit $E=S\otimes E/S$ o\`u $S$ est un fibr\'e de spineurs. En supposant pour simplifier que $M_g$ a une structure de spin, on d\'efinit globalement sur $M$ un courant de de Rham 
\be
C_g = \widehat{A}(M_g)\frac{\ch(E/S,g)}{\ch(S_N,g)} \cap [M_g]\ ,
\ee
avec $\widehat{A}(M_g)$, $\ch(E/S,g)$ et $\ch(S_N,g)$ les classes caract\'eristiques entrant dans le th\'eor\`eme de Lefschetz g\'en\'eralis\'e (voir \cite{BGV} pour le cas non spin). $C_g$ est ferm\'e et invariant par le sous-groupe centralisateur de $g$. Les m\'ethodes de d\'eveloppement asymptotique du noyau de la chaleur \cite{BGV, Gi} conduisent \`a la formule de localisation suivante.

\begin{proposition}[\cite{P6}]\label{ploc}
Consid\'erons un module de Clifford $G$-\'equivariant $E$ et un op\'erateur de Dirac g\'en\'eralis\'e $G$-invariant $D:\cinfc(E)\to\cinfc(E)$. La composante de degr\'e $n$ du cup-produit $\ch(E,\rho,\sqrt{t}D)\cdot \ch(e)$ est une $n$-forme non commutative $\ch_n(\sqrt{t}D)\in \Om^n\Bc$ identifiable \`a une fonction sur $G^n\cup G^{n+1}$. Sa limite $t\downarrow 0$ en un point $\gt\in G^n\cup G^{n+1}$ est donn\'ee par la formule de localisation
\be
\lim_{t\downarrow 0}\ch_n(\sqrt{t}D)(\gt) = \sum_{M_g}\frac{(-)^{q/2}}{(2\pi i)^{d/2}}\int_{M_g}\widehat{A}(M_g)\frac{\ch(E/S,g)}{\ch(S_N,g)}\, \ch_n(e)(\gt)\ ,
\ee
o\`u $g\in G$ est le produit de concatenation $g_n\ldots g_1$ (resp. $g_n\ldots g_0$) si $\gt=(g_1,\ldots, g_n)$ (resp. $\gt=(g_0,\ldots, g_n)$). La somme est prise sur toutes les sous-vari\'et\'es fixes $M_g$ de dimension $d$ et codimension $q=\dim M-d$. \cqfd
\end{proposition}

\begin{remark}\textup{Puisque l'action de $G$ sur $M$ est propre, un \'el\'ement $g\in G$ ne peut avoir de points fixes que s'il appartient \`a un sous-groupe compact. Ainsi le support de la fonction limite $\lim_{t\downarrow 0}\ch_n(\sqrt{t}D)$ est restreint aux points $\gt\in G^n\cup G^{n+1}$ dont le produit de concat\'enation appartient \`a un sous-groupe compact de $G$. }
\end{remark}

\begin{example}\textup{Lorsque $G$ est un groupe compact agissant sur une vari\'et\'e compacte $M$, toute compl\'etion admissible de $C_c(G)$ est topologiquement \'equivalente \`a l'alg\`ebre de Banach des fonctions int\'egrables $\Bc=L^1(G)$ relativement \`a la mesure de Haar. En dimension paire l'indice analytique $\mu(Q)\in \Kt_0(\Bc)$ est une repr\'esentation virtuelle de $G$, et toute l'information int\'eressante sur le caract\`ere de Chern $\ch(\mu(Q))$ est concentr\'ee dans sa composante de degr\'e z\'ero $\ch_0(\sqrt{t}D)\in \Bc$. La classe de $K$-th\'eorie canonique $[e]\in \Kt_0(\Ac)$ est repr\'esent\'ee par l'idempotent 
$$
e(g,x)=1\qquad \forall\ (g,x)\in G\times M\ ,
$$
en supposant la mesure de Haar normalis\'ee. Alors la forme $\ch_0(e)\in \Omc^*(M)\otimes\Bc$ correspond simplement \`a la fonction $G\to \Omc^*(M)$ identiquement \'egale \`a 1. Par cons\'equent, la formule de localisation donne, en tout point $g\in G$,
\be
\lim_{t\downarrow 0} \ch_0(tD)(g) = \sum_{M_g}\frac{(-)^{q/2}}{(2\pi i)^{d/2}}\int_{M_g}\widehat{A}(M_g)\frac{\ch(E/S,g)}{\ch(S_N,g)}\ .
\ee
On retrouve ainsi le th\'eor\`eme de Lefschetz g\'en\'eralis\'e par Atiyah-Segal-Singer \cite{AS}.}
\end{example}

\begin{example}\textup{Lorsque $G$ est un groupe discret d\'enombrable agissant proprement et librement sur $M$, le quotient $X=G\backslash M$ est une vari\'et\'e compacte. On est donc en pr\'esence d'une fibration principale $M\stackrel{G}{\to} X$. Notons que pour les actions libres le support de la fonction $\lim_{t\downarrow 0}\ch_n(\sqrt{t}D)$ est restreint aux points $\gt\in G^n\cup G^{n+1}$ dont le produit de concat\'enation est $g=1$. Pour extraire l'information du caract\`ere de Chern $\ch(\mu(Q))$, il suffit donc de l'\'evaluer sur la cohomologie cyclique localis\'ee en l'unit\'e, autrement dit la cohomologie de groupe. Soit $v\in Z^n(G;\cc)$ un cocycle de groupe \`a croissance polyn\^omiale relativement \`a une distance sur $G$, et soit $\varphi_v$ le $n$-cocycle cyclique sur l'alg\`ebre du groupe $C_c(G)$ qui lui est associ\'e (\cite{CM90}). En choisissant une compl\'etion admissible $\Bc$ au moyen de la distance et d'un param\`etre $\al$ suffisamment \'elev\'e, $\varphi_v$ s'\'etend en un $n$-cocycle cyclique continu sur $\Bc$. Son couplage avec le caract\`ere de Chern $\ch_n(e)\in \Omc^*(M)\otimes\Om^n\Bc$ donne une forme diff\'erentielle \`a support compact
$$
\langle \varphi_v\, , \,\ch_n(e) \rangle \in \Omc^*(M)\ ,
$$
dont l'image directe par la projection \'etale $M\to X$ est une forme diff\'erentielle ferm\'ee. On calcule que sa classe de cohomologie dans $H^*(X)$ co\"incide avec l'image r\'eciproque de la classe de cohomologie de groupe $v\in H^n(G)\cong H^n(BG)$ sous l'application classifiante $f: X\to BG$. Tout op\'erateur de Dirac g\'en\'eralis\'e op\'erant sur les sections d'un module de Clifford $E\to X$ peut se relever en une classe de $K$-homologie $[Q]\in K^G_*(M)$, et la formule de localisation donne 
\be
\langle \varphi_v\, , \,\ch(\mu(Q)) \rangle = \frac{1}{(2\pi i)^{d/2}} \int_X \widehat{A}(X) \ch(E/S) f^*(v)
\ee
avec $d$ la dimension de $X$. Ici le genre $\widehat{A}(X)$ et le caract\`ere de Chern $\ch(E/S)$ sont des classes de cohomologie sur $X$. On retrouve ainsi le th\'eor\`eme de Connes et Moscovici pour les $G$-recouvrements \cite{CM90}.  }
\end{example}

\chapter{Images directes}\label{csec}

On \'etablit dans ce chapitre la compatibilit\'e entre le caract\`ere de Chern bivariant construit dans le chapitre \ref{ccar} et les morphismes d'image directe en $K$-th\'eorie. A cette fin il est n\'ecessaire de se restreindre \`a la cat\'egorie des $m$-alg\`ebres de Fr\'echet. Pour de telles alg\`ebres on distingue deux types d'invariants. D'un c\^ot\'e la $K$-th\'eorie topologique \cite{Ph} et l'homologie cyclique p\'eriodique sont des invariants primaires, c'est-\`a-dire stables sous homotopie et v\'erifiant la p\'eriodicit\'e de Bott. Elles portent une information de nature essentiellement topologique. D'un autre c\^ot\'e, les filtrations naturelles du complexe cyclique permettent de d\'efinir la $K$-th\'eorie multiplicative \cite{K1, K2} et les versions instables de l'homologie cyclique. Ce sont des invariants secondaires qui portent une information plus fine de nature g\'eom\'etrique. La relation entre invariants primaires et secondaires se fait au moyen de suites exactes longues. \\
Le th\'eor\`eme de Grothendieck-Riemann-Roch formul\'e dans \cite{P8} donne une version pr\'ecis\'ee du th\'eor\`eme de l'indice. En effet on montre qu'un quasihomomorphisme $p$-sommable entre deux $m$-alg\`ebres de Fr\'echet $\Ac$ et $\Bc$ induit des morphismes d'image directe simultan\'ement pour les invariants primaires et secondaires, et leur compatibilit\'e s'exprime au moyen d'un diagramme commutatif reliant les suites exactes longues. Comme on travaille ici avec les versions filtr\'ees de l'homologie cyclique, le caract\`ere de Chern du quasihomomorphisme est interpr\'et\'e comme une classe dans la cohomologie cyclique bivariante non-p\'eriodique $HC^*(\Ac,\Bc)$ au lieu de la cohomologie cyclique bivariante enti\`ere. Cela permet d'exploiter au mieux l'information port\'ee par les invariants secondaires. Mentionnons que les morphismes d'images directe en $K$-th\'eorie multiplicative fournissent un analogue non-commutatif des morphismes d'image directe en cohomologie de Deligne. Ces \'el\'ements sont expos\'es en d\'etail dans l'article\\

\noindent \cite{P8} D. Perrot: Secondary invariants for Fr\'echet algebras and quasihomomorphisms, {\it Documenta Math.} {\bf 13} (2008) 275-363.

\section{Invariants primaires et secondaires}

Soit $\Ac$ une $m$-alg\`ebre de Fr\'echet. Sa topologie est engendr\'ee par une famille d\'enombrable de semi-normes $q$ v\'erifiant la propri\'et\'e de sous-multiplicativit\'e 
$$
q(ab)\leq q(a)q(b)\qquad \forall\ a,b\in \Ac\ . 
$$
On peut aussi repr\'esenter $\Ac$ comme limite projective d'une suite d'alg\`ebres de Banach. Phillips d\'efinit dans \cite{Ph} la $K$-th\'eorie topologique des $m$-alg\`ebres de Fr\'echet par analogie avec la $K$-th\'eorie topologique des alg\`ebres de Banach. L'alg\`ebre $\Kc$ des op\'erateurs compacts lisses est une $m$-alg\`ebre de Fr\'echet, ainsi que le produit tensoriel projectif compl\'et\'e $\Kc\hotimes\Ac$. Soit $(\Kc\hotimes\Ac)^+$ son unitarisation, et $p_0\in M_2(\Kc\hotimes\Bc)^+$ la matrice idempotente $p_0= \bigl(\begin{smallmatrix} 1 & 0 \\ 0 & 0 \end{smallmatrix}\bigr) $. Les groupes de $K$-th\'eorie topologique en degr\'e pair et impair sont d\'efinis par
\beq
\Kt_0(\Ac) &=& \{\mbox{classes d'homotopie diff\'erentiable d'idempotents $e\in M_2(\Kc\hotimes\Ac)^+$}\non\\
 && \qquad \mbox{tels que}\ e- p_0 \in M_2(\Kc\hotimes\Ac)\ \}\non\\
\Kt_1(\Ac) &=& \{\mbox{classes d'homotopie diff\'erentiable d'inversibles $u\in (\Kc\hotimes\Ac)^+$}\non\\
 && \qquad  \mbox{tels que}\ u-1 \in \Kc\hotimes\Ac\ \}\non
\eeq
Soit $\cinf(0,1)$ l'alg\`ebre de Fr\'echet des fonctions lisses sur l'intervalle, qui s'annulent aux extr\'emit\'es ainsi que toutes leurs d\'eriv\'ees. La suspension lisse de $\Ac$ est le produit tensoriel projectif $S\Ac=\Ac\hotimes\cinf(0,1)$. Alors la $K$-th\'eorie topologique v\'erifie la p\'eriodicit\'e de Bott $\Kt_0(S\Ac)\cong \Kt_1(\Ac)$ et $\Kt_1(S\Ac)\cong\Kt_0(\Ac)$. On peut donc par commodit\'e introduire les groupes $\Kt_n(\Ac)$ en tout degr\'e $n\in\zz$ par
\be
\Kt_n(\Ac) = \left\{ \begin{array}{cl} 
\Kt_0(\Ac) & n\ \mbox{pair} \\
\Kt_1(\Ac) & n\ \mbox{impair} \end{array} \right.
\ee
En relation avec la th\'eorie de l'indice non-commutative on est souvent amen\'e \`a consid\'erer la $K$-th\'eorie du produit tensoriel projectif $\Ic\hotimes\Ac$, o\`u $\Ic$ est une $m$-alg\`ebre de Fr\'echet $p$-sommable ($p$ entier), typiquement une classe de Schatten $\ell^p$. En utilisant un cocycle cyclique bivariant associ\'e \`a la trace sur la puissance $p$-i\`eme de $\Ic$ nous avons construit dans \cite{P8} un caract\`ere de Chern
\be
\Kt_n(\Ic\hotimes\Ac) \to HP_n(\Ac)\label{caract}
\ee
\`a valeurs dans l'homologie cyclique p\'eriodique de $\Ac$. Rappelons que $HP_n(\Ac)$ est l'homologie en degr\'e $n\mod 2$ du $(b+B)$-complexe des formes diff\'erentielles compl\'et\'e en prenant le produit direct $\Omh\Ac= \prod_{k\geq 0}\Om^k\Ac$. De mani\`ere \'equivalente \cite{CQ1}, pour toute extension quasi-libre de $m$-alg\`ebres de Fr\'echet $0 \to \Jc \to \Rc \to \Ac \to 0$ l'homologie cyclique p\'eriodique de $\Ac$ est calcul\'ee par le $X$-complexe de pro-alg\`ebre
$$
X(\Rch) \ :\ \Rch\ \rightleftarrows\ \Om^1\Rch_{\nat}\ ,
$$
o\`u $\Rch$ est la compl\'etion $\Jc$-adique de $\Rc$. Tout comme la $K$-th\'eorie topologique, l'homologie cyclique p\'eriodique v\'erifie l'isomorphisme $HP_{n+2}(\Ac)\cong HP_n(\Ac)$ par construction. De plus, les groupes $\Kt_n(\Ic\hotimes\Ac)$ et $HP_n(\Ac)$ sont stables sous homotopie diff\'erentiable et d\'efinissent les invariants primaires de $\Ac$.\\ 

Les invariants secondaires de $\Ac$ sont un m\'elange de $K$-th\'eorie topologique et des versions instables de l'homologie cyclique. Ces derni\`eres sont d\'efinies \`a partir des filtrations naturelles du complexe cyclique $\Omh\Ac$ par le degr\'e des formes diff\'erentielles \cite{C0}, ou de mani\`ere \'equivalente \`a partir de la filtration $\Jc$-adique sur $X(\Rch)$. En particulier l'homologie cyclique non-p\'eriodique de Connes $HC_n(\Ac)$ est calcul\'ee par un complexe quotient, et sa relation avec l'homologie cyclique p\'eriodique s'inscrit dans une suite exacte longue de type $SBI$ (voir \cite{P8})
$$
\ldots \longrightarrow  HP_{n+1}(\Ac) \stackrel{S}{\longrightarrow}   HC_{n-1}(\Ac) \stackrel{B}{\longrightarrow}   HN_n(\Ac) \stackrel{I}{\longrightarrow}  HP_n(\Ac) \longrightarrow  \ldots
$$
o\`u $HN_n(\Ac)$ est l'homologie cyclique n\'egative de $\Ac$. Rappelons que $HC_n(\Ac)=0$ pour $n<0$ et par cons\'equent $HN_n(\Ac)\cong HP_n(\Ac)$ pour $n\leq 0$. A l'aide de la filtration du complexe cyclique de $\Ac$ on peut alors fabriquer une version secondaire de la $K$-th\'eorie. En s'inspirant du travail de Karoubi \cite{K1, K2} nous avons d\'efini dans \cite{P8}, pour toute alg\`ebre $p$-sommable $\Ic$, des groupes de $K$-th\'eorie multiplicative $MK^{\Ic}_n(\Ac)$, $n\in\zz$, qui s'ins\`erent dans une suite exacte longue 
$$
\ldots \Kt_{n+1}(\Ic\hotimes\Ac) \to HC_{n-1}(\Ac) \stackrel{\delta}{\to} MK^{\Ic}_n(\Ac)  \to \Kt_n(\Ic\hotimes\Ac)  \to HC_{n-2}(\Ac)  \ldots
$$
Pour $n$ pair, toute classe dans $MK^{\Ic}_n(\Ac)$ est repr\'esent\'ee par un couple $(e,\te)$ avec $e$ un idempotent qui d\'efinit une classe $[e]\in \Kt_0(\Ic\hotimes\Ac)$, et $\te$ une cha\^{\i}ne cyclique qui transgresse le caract\`ere de Chern de $e$ dans l'homologie de de Rham non-commutative $HD_{n-2}(\Ac)$ (voir \cite{P8}). Pour $n$ impair, toute classe dans $MK^{\Ic}_n(\Ac)$ est repr\'esent\'ee par un couple $(u,\te)$ avec $u$ un inversible et $\te$ une cha\^{\i}ne cyclique munis de propri\'et\'es analogues. Lorsque $n\leq 0$ la suite exacte implique l'isomorphisme $MK^{\Ic}_n(\Ac)  \cong \Kt_n(\Ic\hotimes\Ac)$, tandis que les invariants secondaires proprement dits n'appara\^issent que pour $n>0$. Mentionnons que la $K$-th\'eorie multiplicative est intimement reli\'ee \`a la $K$-th\'eorie alg\'ebrique sup\'erieure \cite{R}. Il existe aussi un caract\`ere de Chern n\'egatif
\be
MK^{\Ic}_n(\Ac)\to HN_n(\Ac)
\ee 
qui rel\`eve le caract\`ere de Chern en $K$-th\'eorie topologique (\ref{caract}). Plus pr\'ecis\'ement on a le r\'esultat suivant.

\begin{proposition}[\cite{P8}]\label{pmul}
Soit $\Ac$ une $m$-alg\`ebre de Fr\'echet. Pour toute $m$-alg\`ebre de Fr\'echet $p$-sommable $\Ic$, les caract\`eres de Chern en $K$-th\'eorie topologique et multiplicative induisent une transformation entre les suites exactes longues 
\be
\vcenter{\xymatrix{
\Kt_{n+1}(\Ic\hotimes\Ac) \ar[r] \ar[d] & HC_{n-1}(\Ac) \ar[r]^{\delta} \ar@{=}[d] & MK^{\Ic}_n(\Ac)  \ar[r] \ar[d] & \Kt_n(\Ic\hotimes\Ac)  \ar[d] \\
HP_{n+1}(\Ac) \ar[r]^S & HC_{n-1}(\Ac) \ar[r]^{\widetilde{B}} & HN_n(\Ac) \ar[r]^I & HP_n(\Ac) }} \label{mul}
\ee
o\`u $\widetilde{B}$ est l'application de bord $-\sqrt{2\pi i}\, B$. \cqfd
\end{proposition}
Dans le cas o\`u $\Ic=\cc$ est l'alg\`ebre $1$-sommable des nombres complexes, on \'ecrira simplement $MK^{\cc}_n(\Ac)=MK_n(\Ac)$. Pour une alg\`ebre de Banach $\Ac$ on retrouve ainsi la $K$-th\'eorie multiplicative de Karoubi \cite{K1, K2}.

\begin{example}\textup{Prenons $\Ac=\cc$ et $\Ic=\ell^p$ une classe de Schatten sur un espace de Hilbert s\'eparable. La $K$-th\'eorie topologique de $\Ic$ est connue: $\Kt_0(\Ic)=\zz$ et $\Kt_1(\Ic)=0$. La suite exacte longue implique donc
$$
MK^{\Ic}_n(\cc)= \left\{ \begin{array}{ll}
\zz & \mbox{$n\leq 0$ pair} \\
\cc^{\times} & \mbox{$n>0$ impair} \\
0 & \mbox{autrement} \end{array} \right.
$$
La $K$-th\'eorie multiplicative de $\cc$ en degr\'e $n>0$ impair est le r\'eceptacle naturel des morphismes r\'egulateurs, voir l'exemple \ref{ereg} et la fin du chapitre \ref{cano}. }
\end{example}

\begin{example}\label{edel}\textup{Lorsque $\Ac=\cinf(M)$ est l'alg\`ebre commutative des fonctions lisses sur une vari\'et\'e compacte $M$, la $K$-th\'eorie multiplicative $MK_n(\Ac)$ a des propri\'et\'es analogues \`a la cohomologie de Deligne de $M$. Pour tout demi-entier $q$ notons $\zz(q)$ le sous-groupe additif $(2\pi i)^q\zz\subset \cc$. Par d\'efinition, le groupe de cohomologie de Deligne $H_{\Dc}^n(M;\zz(n/2))$ est l'hyperhomologie en degr\'e $n$ du complexe de faisceaux
$$
0\longrightarrow \zzb(n/2) \longrightarrow \Omb^0 \stackrel{d}{\longrightarrow} \Omb^1 \stackrel{d}{\longrightarrow} \ldots \stackrel{d}{\longrightarrow} \Omb^{n-1} \longrightarrow 0 
$$
avec le faisceau constant $\zzb(n/2)$ situ\'e en degr\'e 0 et le faisceau $\Omb^k$ des $k$-formes diff\'erentielles sur $M$ situ\'e en degr\'e $k+1$. De ce complexe on extrait imm\'ediatement un morphisme horizontal de la cohomologie de Deligne $H_{\Dc}^n(M;\zz(n/2))$ vers la cohomologie de \v{C}ech $\check{H}^n(Q;\zz(n/2))$, et un morphisme vertical vers les $n$-formes diff\'erentielles ferm\'ees $\Zdr^n(M)$, qui co\"{\i}ncident en cohomologie de de Rham:
$$
\vcenter{\xymatrix{
H^n_{\Dc}(M;\zz(n/2)) \ar[r] \ar[d]_d  & \check{H}^n(M;\zz(n/2)) \ar[d]^{\otimes\cc}  \\
\Zdr^n(M) \ar[r]  & \Hdr^n(M) }} 
$$
$\Zdr^n(M)$ et $\Hdr^n(M)$ sont inclus comme facteurs directs respectivement dans $HN_n(\Ac)$ et $HP_n(\Ac)$ pour l'alg\`ebre $\Ac=\cinf(M)$. De plus, nous avons construit explicitement dans \cite{P8} un morphisme $MK_n(\Ac)\to H^n_{\Dc}(M;\zz(n/2))$ en degr\'e $n\leq 2$ qui envoie le diagramme ci-dessus dans le carr\'e commutatif extrait de (\ref{mul})
$$
\vcenter{\xymatrix{
 MK_n(\Ac)  \ar[r] \ar[d] & \Kt_n(\Ac)  \ar[d] \\
 HN_n(\Ac) \ar[r] & HP_n(\Ac) }} 
$$
Pousser la comparaison en degr\'e $n\geq 3$ est plus d\'elicat, car il n'y a pas d'application \'evidente de la cohomologie de Deligne vers la $K$-th\'eorie multiplicative pour des raisons d'int\'egralit\'e (en d'autres termes les r\'eseaux de cohomologie de \v{C}ech et de $K$-th\'eorie topologique deviennent incompatibles). Cela montre qu'il serait plus int\'eressant de comparer la $K$-th\'eorie multiplicative \`a une version mieux adapt\'ee \`a la $K$-th\'eorie que la cohomologie de Deligne, comme par exemple les $K$-caract\`eres diff\'erentiels de \cite{BM}. }
\end{example}

\section{Quasihomomorphismes}

Soient $\Ac$ et $\Bc$ deux $m$-alg\`ebres de Fr\'echet. Nous dirons qu'un $\Ac$-$\Bc$-bimodule born\'e $(H,\rho,F)$ de degr\'e pair est mis sous la forme d'un \emph{quasihomomorphisme $p$-sommable} s'il existe un espace de Fr\'echet trivialement gradu\'e $L$, une alg\`ebre d'op\'erateurs $p$-sommables $\Ic\subset \End(L)$ et une alg\`ebre d'endomorphismes $\Ec\subset \End(L_{\Bc})$ contenant $\Ic\hotimes\Bc$ comme id\'eal bilat\`ere tels que 
$$
H=\left(\begin{matrix} L \\ L \end{matrix}\right) \ ,\qquad \rho=\left(\begin{matrix} \rho_+ & 0 \\ 0 & \rho_- \end{matrix}\right) \ ,\qquad F=\left(\begin{matrix} 0 & 1 \\ 1 & 0 \end{matrix}\right)
$$
avec $\rho_{\pm}: \Ac\to\Ec$ et $\rho_+-\rho_-: \Ac\to \Ic\hotimes\Bc$. De la m\^eme mani\`ere, un bimodule $(H,\rho,D)$ de degr\'e impair est un quasihomomorphisme $p$-sommable s'il est de la forme
$$
H=\left(\begin{matrix} L \\ L \end{matrix}\right)\otimes C_1 \ ,\quad \rho=\left(\begin{matrix} \rho_{++} & \rho_{+-} \\ \rho_{-+} & \rho_{--} \end{matrix}\right)\otimes 1 \ ,\quad F=\left(\begin{matrix} 1 & 0 \\ 0 & -1 \end{matrix}\right)\otimes\eps
$$
avec $\rho_{++},\rho_{--}:\Ac\to \Ec$ et $\rho_{+-},\rho_{-+}:\Ac\to \Ic\hotimes\Bc$. Dans tous les cas on peut canoniquement construire \`a partir de $\Ec$ et son id\'eal $\Ic\hotimes\Bc$ une sous-alg\`ebre $\zz_2$-gradu\'ee $\Ec^s \triangleright\Ic^s\hotimes\Bc$ des endomorphismes de $H_{\Bc}$ qui prend automatiquement en compte les conditions impos\'ees sur l'expression matricielle de $\rho$, voir \cite{P8}. Comme l'op\'erateur $F$ est mis sous une forme canonique, toute r\'ef\'erence \`a l'espace $H$ devient inutile et il suffit de ne retenir que l'homomorphisme 
\be
\rho:\Ac\to \Ec^s\triangleright\Ic^s\hotimes\Bc\ ,\label{quasi}
\ee
dont l'image est contenue dans la sous-alg\`ebre paire de $\Ec^s$. Plus g\'en\'eralement nous d\'efinissons dans \cite{P8} un quasihomomorphisme $p$-sommable de $\Ac$ vers $\Bc$ comme la donn\'ee d'une $m$-alg\`ebre de Fr\'echet trivialement gradu\'ee abstraite $\Ec\triangleright \Ic\hotimes\Bc$ avec $\Ic$ une $m$-alg\`ebre de Fr\'echet $p$-sommable, et d'un homomorphisme continu de la forme (\ref{quasi}). L'op\'erateur $F$ mis sous la forme matricielle ci-dessus agit alors comme multiplicateur de degr\'e impair sur $\Ec^s$. \\

Pour construire le caract\`ere de Chern d'un quasihomomorphisme $p$-sommable $\rho:\Ac\to \Ec^s\triangleright\Ic^s\hotimes\Bc$ dans la cohomologie cyclique bivariante $HC^*(\Ac,\Bc)$, il suffit de reprendre les formules \'etablies pour les bimodules born\'es de la section \ref{sborn} et de contr\^oler leurs propri\'et\'es adiques \cite{P8}. Comme auparavant il s'agit de relever $\rho$ en un quasihomomorphisme entre des extensions universelles de $\Ac$ et $\Bc$; il est donc n\'ecessaire d'imposer une condition d'admissibilit\'e analogue \`a \ref{dadm}.  

\begin{definition}[\cite{P8}]
L'alg\`ebre $\Ec\triangleright \Ic\hotimes\Bc$ est \emph{admissible} relativement \`a une extension $0\to \Jc \to \Rc \to \Bc \to 0$ s'il existe deux $m$-alg\`ebres de Fr\'echet $\Mc\triangleright \Ic\hotimes\Rc$ et $\Nc\triangleright\Ic\hotimes\Jc$ telles que $\Nc^n\cap \Ic\hotimes\Rc = \Ic\hotimes\Jc^n$ pour toute puissance $n\geq 1$, ainsi qu'un diagramme d'extensions
$$
\xymatrix{
0\ar[r] & \Nc \ar[r]  & \Mc \ar[r] & \Ec \ar[r] & 0 \\
0 \ar[r] & \Ic\hotimes \Jc \ar[r] \ar[u] & \Ic\hotimes \Rc \ar[r] \ar[u] & \Ic\hotimes \Bc \ar[r] \ar[u] & 0}
$$
\end{definition}

\begin{example}\textup{Ici encore l'arch\'etype d'alg\`ebre admissible est donn\'e par un produit tensoriel $\Ec=\End(L)\hotimes\Bc$ avec $L$ espace de Hilbert trivialement gradu\'e et $\Ic=\ell^p(L)$. Il suffit alors de prendre $\Mc=\End(L)\hotimes\Rc$ et $\Nc=\End(L)\hotimes\Jc$. Bien entendu il existe d'autres exemples importants d'alg\`ebres admissibles qui ne peuvent pas s'\'ecrire comme produit tensoriel.}
\end{example}

\noindent On construit ensuite les alg\`ebres $\zz_2$-gradu\'ees associ\'ees $\Mc^s\triangleright \Ic^s\hotimes\Rc$ et $\Nc^s\triangleright\Ic^s\hotimes\Jc$ en tenant compte de la parit\'e $i\in\zz_2$ du quasihomomorphisme. La propri\'et\'e universelle de l'alg\`ebre tensorielle $T\Ac$ 
$$
\vcenter{\xymatrix{
0 \ar[r] & J\Ac \ar[r] \ar[d]^{\rho_*} & T\Ac \ar[r] \ar[d]^{\rho_*} & \Ac \ar[r] \ar[d]^{\rho} & 0 \\
0\ar[r] & \Nc^s \ar[r]  & \Mc^s \ar[r] & \Ec^s \ar[r] & 0 }}
$$
donne un rel\`evement $\rho_*: T\Ac\to \Mc^s\triangleright\Ic^s\hotimes\Rc$ qui envoie l'id\'eal $J\Ac$ sur $\Nc^s$. On obtient ainsi un quasihomomorphisme $p$-sommable pour les compl\'etions adiques de $T\Ac$, $\Mc^s$ et $\Rc$ par rapport \`a leurs id\'eaux $J\Ac$, $\Nc^s$ et $\Jc$:
\be
\rho_*: \Th\Ac\to \Mch^s\triangleright\Ic^s\hotimes\Rch\ .
\ee
Fixons maintenant un entier $n\geq p$ de parit\'e $i\mod 2$. Soit $\chih^n\in \hom(\Omh \Th\Ac,X(\Rch))$ l'application lin\'eaire dont les deux composantes non nulles $\chih^n_0:\Om^n \Th\Ac\to \Rch$ et $\chih^n_1:\Om^{n+1}\Th\Ac\to \Om^1\Rch_{\nat}$ d\'efinies par (\ref{chi0})
\beq
&& \chih^n_0(x_0\dd x_1\ldots\dd x_n) = (-)^n\frac{\Gamma(1+\frac{n}{2})}{(n+1)!} \sum_{\la\in S_{n+1}} \eps(\la)\, \tau(x_{\la(0)}[F,x_{\la(1)}]\ldots [F,x_{\la(n)}])\ ,\non\\
&& \chih^n_1(x_0\dd x_1\ldots\dd x_{n+1}) = (-)^n\frac{\Gamma(1+\frac{n}{2})}{(n+1)!} \sum_{i=1}^{n+1}  \tau\nat(x_0[F,x_1]\ldots\dd x_i \ldots [F,x_{n+1}])\ ,\non\\
&& \label{chi}
\eeq
avec $\tau$ la supertrace sur la $n$-i\`eme puissance de $\Ic^s$. On sait que $\chih^n$ est un morphisme du $(b+B)$-complexe des formes diff\'erentielles sur $\Th\Ac$ vers le $X$-complexe de $\Rch$. Sa compos\'ee par l'\'equivalence de Goodwillie $\gamma: X(\Th\Ac)\to \Omh\Th\Ac$ d\'efinit alors un cocycle cyclique bivariant de degr\'e $n$ lorsque $\Rc$ est quasi-libre:

\begin{proposition}[\cite{P8}]\label{padic}
Soit $\rho:\Ac\to\Ec^s\triangleright\Ic^s\hotimes\Bc$ un quasihomomorphisme $p$-sommable de parit\'e $i\in\zz_2$, admissible relativement \`a une extension quasi-libre $0\to \Jc \to \Rc \to \Bc \to 0$, et soit $n\geq p$ un entier, $n\equiv i\mod 2$. Alors le caract\`ere de Chern du quasihomomorphisme d\'efini par la composition 
\be
\ch^n(\rho)\ :\ X(\Th\Ac) \stackrel{\gamma}{\longrightarrow} \Omh \Th\Ac \stackrel{\chih^n}{\longrightarrow} X(\Rch)
\ee
est une classe de cohomologie cyclique bivariante $\ch^n(\rho)\in HC^n(\Ac,\Bc)$ de degr\'e $n$. Les caract\`eres de Chern de degr\'es successifs sont reli\'es par l'op\'eration de suspension $S$ en cohomologie cyclique bivriante
\be
\ch^{n+2}(\rho) \equiv S \ch^n(\rho)\ \in\ HC^{n+2}(\Ac,\Bc)\ ,
\ee
et ainsi d\'efinissent tous la m\^eme classe de cohomologie cyclique bivariante p\'eriodique $\ch(\rho)\in HP^i(\Ac,\Bc)$. \cqfd
\end{proposition}

C'est donc le repr\'esentant $\ch^n(\rho)$ de degr\'e $n$ minimal qui v\'ehicule le maximum d'informations sur le quasihomomorphisme. Les propri\'et\'es int\'eressantes du caract\`ere de Chern bivariant concernent sa stabilit\'e par rapport aux deux relations d'\'equivalences possibles \cite{P8}. Rappelons que $\Bc[0,1]$ d\'esigne le produit tensoriel projectif $ \Bc\hotimes\cinf[0,1]$. C'est aussi l'alg\`ebre des fonctions lisses sur $[0,1]$ \`a valeurs dans $\Bc$ et dont toutes les d\'eriv\'ees d'ordre $\geq 1$ s'annulent aux extr\'emit\'es. On dit que deux quasihomomorphismes $\rho_0:\Ac\to\Ec^s\triangleright\Ic^s\hotimes\Bc$ et $\rho_1:\Ac\to\Ec^s\triangleright\Ic^s\hotimes\Bc$ sont \\

\noindent {\bf homotopes} s'il existe un quasihomomorphisme $\rho:\Ac\to\Ec[0,1]^s\triangleright\Ic^s\hotimes\Bc[0,1]$ dont l'\'evaluation aux extr\'emit\'es du segment donne respectivement $\rho_0$ et $\rho_1$; \\

\noindent {\bf conjugu\'es} s'il existe un \'el\'ement inversible de degr\'e pair $U\in (\Ec^s)^+$ tel que $\rho_1=U^{-1}\rho_0 U$ en tant qu'homomorphisme $\Ac\to \Ec^s$.\\

\noindent Apr\`es stabilisation par des matrices la relation de conjugaison est strictement plus forte que la relation d'homotopie. On montre dans \cite{P8} que deux quasihomomorphismes conjugu\'es d\'efinissent le m\^eme caract\`ere de Chern dans $HC^n(\Ac,\Bc)$ pour toute valeur de $n\geq p$, alors que les caract\`eres de Chern de deux quasihomomorphismes homotopes ne co\"{\i}ncident qu'apr\`es stabilisation par l'op\'eration $S$. On peut r\'esumer ces propri\'et\'es dans le corollaire suivant.

\begin{corollary}[\cite{P8}]\label{cbiv}
Soit $\rho:\Ac\to \Ec^s\triangleright \Ic^s\hotimes\Bc$ un quasihomomorphisme $p$-sommable de parit\'e $i\in\zz_2$, admissible relativement \`a une extension quasi-libre de $\Bc$. Supposons $p\equiv i\mod 2$. Alors le caract\`ere de Chern de degr\'e minimal $\ch^p(\rho)\in HC^p(\Ac,\Bc)$ induit une transformation entre les suites exactes SBI 
$$
\vcenter{\xymatrix{
 HP_{n+1}(\Ac) \ar[r]^S \ar[d] & HC_{n-1}(\Ac) \ar[r]^B \ar[d] & HN_n(\Ac)  \ar[r]^I \ar[d] & HP_n(\Ac)  \ar[d]  \\
 HP_{n-p+1}(\Bc) \ar[r]^S & HC_{n-p-1}(\Bc) \ar[r]^B & HN_{n-p}(\Bc) \ar[r]^I & HP_{n-p}(\Bc) &  }}
$$
invariante sous conjugaison du quasihomomorphisme. De plus la fl\^eche en homologie cyclique p\'eriodique $HP_n(\Ac)\to HP_{n-d}(\Bc)$ est invariante sous homotopie. \cqfd
\end{corollary}
\begin{remark}\textup{Le r\'esultat ci-dessus reste inchang\'e si l'on suppose le quasihomomorphisme seulement $(p+1)$-sommable, avec toujours $p\equiv i\mod 2$. }
\end{remark}

\section{Grothendieck-Riemann-Roch}

Nous pouvons maintenant \'enoncer le r\'esultat principal de ce chapitre, \`a savoir qu'un quasihomomorphisme $p$-sommable admissible induit des morphismes d'image directe pour les invariants primaires et secondaires. Une $m$-alg\`ebre de Fr\'echet $p$-sommable $\Ic$ est \emph{multiplicative} s'il existe un homomorphisme continu (produit externe)
$$
\boxtimes:\Ic\hotimes\Ic\to \Ic
$$
compatible avec la trace (\cite{P8}). Deux produits externes $\boxtimes$ et $\boxtimes'$ sont \emph{\'equivalents} s'il existe un multiplicateur inversible $U$ sur $\Ic$ qui conjugue ces deux produits. L'exemple type d'alg\`ebre $p$-sommable multiplicative est l'id\'eal de Schatten $\Ic=\ell^p(L)$ sur un espace de Hilbert s\'eparable de dimension infinie trivialement gradu\'e, le produit externe \'etant induit par l'identification du produit tensoriel d'espaces de Hilbert $L\otimes_2 L$ avec $L$ \`a un isomorphisme unitaire pr\`es \cite{CT}. \\
Soient maintenant $\Ac$, $\Bc$ deux $m$-alg\`ebres de Fr\'echet et $\rho:\Ac\to \Ec^s\triangleright \Ic^s\hotimes\Bc$ un quasihomomorphisme $p$-sommable de parit\'e $i\in\zz_2$, avec $p\equiv i\mod 2$. Si $\Ic$ est multiplicative alors $\rho$ induit un morphisme d'image directe sur la $K$-th\'eorie topologique
\be
\rho_!: \Kt_n(\Ic\hotimes\Ac) \to \Kt_{n-p}(\Ic\hotimes\Bc)\qquad \forall n\in\zz\ , \label{morK}
\ee
comme en th\'eorie de Kasparov. Par p\'eriodicit\'e de Bott il suffit en effet de d\'efinir cette application pour $n$ impair seulement. Supposons d'abord que $i=0$. Le quasihomomorphisme est alors d\'ecrit par un couple d'homomorphismes $(\rho_+,\rho_-):\Ac\rightrightarrows \Ec$ qui co\"{\i}ncident modulo l'id\'eal $\Ic\hotimes\Bc$. Pour tout inversible $u \in (\Kc\hotimes\Ic\hotimes\Ac)^+$ repr\'esentant une classe dans $\Kt_1(\Ic\hotimes\Ac)$, l'\'el\'ement
$$
\rho_!(u)= \rho_+(u)\rho_-(u)^{-1}\ \in (\Kc\hotimes\Ic\hotimes\Ic\hotimes\Bc)^+
$$
est un inversible repr\'esentant une classe dans $\Kt_1(\Ic\hotimes\Bc)$ apr\`es usage du produit externe $\boxtimes:\Ic\hotimes\Ic\to\Ic$. Supposons maintenant que $i=1$. Le quasihomomorphisme est alors d\'ecrit par un homomorphisme $\rho:\Ac\to \bigl( \begin{smallmatrix} \Ec & \Ic\hotimes\Bc \\ \Ic\hotimes\Bc & \Ec \end{smallmatrix} \bigr)$. Soit $p_0$ l'idempotent $\bigl( \begin{smallmatrix} 1 & 0 \\ 0 & 0 \end{smallmatrix} \bigr)$. Pour tout inversible $u \in (\Kc\hotimes\Ic\hotimes\Ac)^+$, l'\'el\'ement
$$
\rho_!(u)= \rho(u)^{-1}p_0\rho(u)\ \in M_2(\Kc\hotimes\Ic\hotimes\Ic\hotimes\Bc)^+
$$
est un idempotent repr\'esentant une classe dans $\Kt_0(\Ic\hotimes\Bc)$ apr\`es usage du produit externe, d'o\`u l'application (\ref{morK}). \\
Si de plus le quasihomomorphisme est admissible relativement \`a une extension quasi-libre de $\Bc$, on sait d'apr\`es le corollaire \ref{cbiv} que le caract\`ere de Chern de degr\'e minimal $\ch^p(\rho)\in HC^p(\Ac,\Bc)$ induit un morphisme
\be
\ch^p(\rho): HC_n(\Ac) \to HC_{n-p}(\Bc)\qquad \forall n\in\zz\ .\label{morC}
\ee
En combinant (\ref{morK}) et (\ref{morC}) on peut aussi construire des images directes en $K$-th\'eorie multiplicative. Le th\'eor\`eme de Grothendieck-Riemann-Roch d\'emontr\'e dans \cite{P8} exprime la compatibilit\'e entre tous ces morphismes et les suites exactes longues reliant $K$-th\'eorie et homologie cyclique:

\begin{theorem}[\cite{P8}]\label{trr}
Soit $\Ic$ une alg\`ebre $p$-sommable et multiplicative, et $\rho:\Ac\to \Ec^s\triangleright \Ic^s\hotimes\Bc$ un quasihomomorphisme de parit\'e $i\in \zz_2$ admissible relativement \`a une extension quasi-libre de $\Bc$. Supposons $p\equiv i\mod 2$. Alors $\rho$ induit un morphisme d'image directe $\rho_!:MK^{\Ic}_n(\Ac)\to MK^{\Ic}_{n-p}(\Bc)$ qui s'ins\`ere dans un diagramme gradu\'e-commutatif 
\be
\vcenter{\xymatrix{
\Kt_{n+1}(\Ic\hotimes\Ac) \ar[r] \ar[d]^{\rho_!} & HC_{n-1}(\Ac) \ar[r] \ar[d]^{\ch^p(\rho)} & MK^{\Ic}_n(\Ac)  \ar[r] \ar[d]^{\rho_!}  & \Kt_n(\Ic\hotimes\Ac)  \ar[d]^{\rho_!}  \\
\Kt_{n+1-p}(\Ic\hotimes\Bc) \ar[r]  & HC_{n-1-p}(\Bc) \ar[r]  & MK^{\Ic}_{n-p}(\Bc)  \ar[r]  & \Kt_{n-p}(\Ic\hotimes\Bc) }}  \label{rr}
\ee
Les fl\^eches verticales sont invariantes sous conjugaison de $\rho$; la fl\^eche en $K$-th\'eorie topologique est aussi invariante sous homotopie. De plus (\ref{rr}) est compatible avec le diagramme du corollaire \ref{cbiv} (avec $B$ multipli\'e par un facteur $-\sqrt{2\pi i}$) en prenant les caract\`eres de Chern $MK^{\Ic}_n\to HN_n$ et $\Kt_n(\Ic\hotimes\ .\ )\to HP_n$. \cqfd
\end{theorem}
\begin{remark}\textup{Dans le cas pair $i=0$, le r\'esultat reste inchang\'e si le quasihomomorphisme est seulement $(p+1)$-sommable avec toujours $p\equiv i\mod 2$. Par contre dans le cas impair $i=1$ on doit garder la condition de $p$-sommabilit\'e.}
\end{remark}
La d\'emonstration est une v\'erification purement calculatoire de la commutativit\'e du diagramme (\ref{rr}), utilisant les formules explicites pour $\ch^p(\rho)$. A la lumi\`ere de l'exemple \ref{edel}, le morphisme en $K$-th\'eorie multiplicative peut se voir comme une version non-commutative de l'op\'eration ``int\'egration des classes de cohomologie de Deligne le long des fibres d'une submersion'', l'entier $p$ correspondant \`a la dimension des fibres. Les caract\`eres de Chern p\'eriodique et n\'egatif ne sont pas repr\'esent\'es dans le diagramme (\ref{rr}). En fait la compatibilit\'e entre le corollaire \ref{cbiv} et le th\'eor\`eme ci-dessus s'exprime au moyen de diagrammes cubiques, par exemple
$$
\vcenter{\xymatrix@!0@=3.5pc{ & \Kt_n(\Ic\hotimes\Ac) \ar[rr] \ar'[d][dd] & & \Kt_{n-p}(\Ic\hotimes\Bc) \ar[dd] \\
MK^{\Ic}_n(\Ac) \ar[ur] \ar[rr] \ar[dd] & & MK^{\Ic}_{n-p}(\Bc) \ar[ur] \ar[dd] \\
 & HP_n(\Ac) \ar'[r][rr] & & HP_{n-p}(\Bc) \\
HN_n(\Ac) \ar[ur] \ar[rr] & & HN_{n-p}(\Bc) \ar[ur] & }}
$$
Le carr\'e en arri\`ere-plan d\'ecrit la partie purement topologique, c'est-\`a-dire invariante d'homotopie, du th\'eor\`eme de Grothendieck-Riemann-Roch. Le carr\'e en avant-plan rel\`eve donc cette situation au niveau des invariants secondaires.

\begin{example}\label{ereg}\textup{Soit $\Ec=\End(L)$ l'alg\`ebre de Banach des endomorphismes born\'es sur un espace de Hilbert trivialement gradu\'e $L$, et soit $\Ic=\ell^p(L)$ l'id\'eal de Schatten des op\'erateurs $p$-sommables. Lorsque $\Ac$ est quelconque et $\Bc=\cc$, un quasihomomorphisme $\rho:\Ac\to \Ec^s\triangleright\Ic^s$ est un module de Fredholm $p$-sommable sur $\Ac$ \cite{C0}.  En choisissant $n=p+1$ le diagramme (\ref{rr}) se r\'eduit \`a 
$$
\xymatrix{
 & \Kt_{p+2}(\Ic\hotimes\Ac) \ar[r] \ar[d]^{\rho_!} & HC_{p}(\Ac) \ar[r] \ar[d]^{\ch^p(\rho)} & MK^{\Ic}_{p+1}(\Ac)  \ar[r] \ar[d]^{\rho_!}  & \Kt_{p+1}(\Ic\hotimes\Ac)  \ar[d]^{\rho_!}  \\
0 \ar[r] & \zz \ar[r]  & \cc \ar[r]  & \cc^{\times}  \ar[r]  & 0  }
$$
Le morphisme en $K$-th\'eorie topologique $\Kt_p(\Ic\hotimes\Ac)\to\zz$ correspond \`a la fl\^eche d'indice, tandis qu'en $K$-th\'eorie multiplicative $MK^{\Ic}_{p+1}(\Ac)\to \cc^{\times}$ est un \emph{r\'egulateur}. Un morphisme analogue a \'et\'e introduit par Connes et Karoubi dans le contexte de la $K$-th\'eorie alg\'ebrique \cite{CK}. Le th\'eor\`eme \ref{trr} permet donc de g\'en\'eraliser la notion de r\'egulateur \`a des alg\`ebres cibles $\Bc$ quelconques. Dans cette optique il serait int\'eressant de chercher \`a fabriquer des analogues non-commutatifs de la torsion analytique sup\'erieure associ\'ee \`a une submersion \cite{BL}.}
\end{example}

\begin{example}\textup{Soit $\Ac$ une $m$-alg\`ebre de Fr\'echet munie d'une trace continue $\tau:\Ac\to \cc$. Regardons $\Ec=\Ic=\Ac$ comme une alg\`ebre $1$-sommable et posons $\Bc=\cc$. Le quasihomomorphisme $\rho:\Ac\to \Ac^s\triangleright \Ac^s\hotimes\cc$ de degr\'e pair d\'efini par $\rho_+=\id$ et $\rho_-=0$ est donc $1$-sommable et donne un diagramme commutatif dont les lignes sont exactes:
$$
\xymatrix{
 \Kt_0(\Ac) \ar[r] \ar@{=}[d] & \Ac/[\Ac,\Ac] \ar[r] \ar[d]^{\tau} & MK_1(\Ac)  \ar[r] \ar[d]^{\rho_!}  & \Kt_1(\Ac)  \ar@{=}[d]  \\
\Kt_0(\Ac)  \ar[r]^{\tau_*}  & \cc \ar[r]  & MK^{\Ac}_1(\cc)  \ar[r]  & \Kt_1(\Ac)   }
$$
Par ailleurs il existe un morphisme naturel du groupe $\GL_{\infty}(\Ac)$ des matrices inversibles vers la $K$-th\'eorie multiplicative $MK_1(\Ac)$, qui envoie une matrice $u$ sur le couple $(u,0)$. En composant par $\rho_!$ on obtient donc un morphisme de groupes
$$
\Det_{\tau}: \GL_{\infty}(\Ac) \to MK^{\Ac}_1(\cc)
$$
qui se factorise \`a travers la $K$-th\'eorie alg\'ebrique $\Ka_1(\Ac)$. D\'esignons maintenant par $\GL^0_{\infty}(\Ac)$ le noyau du morphisme $\GL_{\infty}(\Ac)\to \Kt_1(\Ac)$ (ce sont moralement les matrices inversibles homotopes \`a l'unit\'e). Alors l'image de l'application $\Det_{\tau}: \GL^0_{\infty}(\Ac)\to MK^{\Ac}_1(\cc)$ est incluse dans le sous-groupe $\cc/\tau_*\Kt_0(\Ac)$, et l'on retrouve ainsi le d\'eterminant introduit par de la Harpe et Skandalis pour les alg\`ebres de Banach \cite{HS}.}
\end{example}

\chapter{Formule locale d'anomalie}\label{cano}

On explique ici la mani\`ere obtenir des formules \emph{locales} pour le caract\`ere de Chern d'un quasihomomorphisme, par un proc\'ed\'e de renormalisation. Le r\^ole central est jou\'e par la \emph{cocha\^{\i}ne \^eta renormalis\'ee}, introduite comme s\'erie formelle (non convergente) dans le complexe cyclique bivariant. Son bord est toujours une somme finie de termes locaux qui repr\'esente le caract\`ere de Chern, quelle que soit la renormalisation choisie. Par exemple, en pr\'esence d'un op\'erateur de Dirac la renormalisation par fonction z\^eta donne une formule de r\'esidus qui g\'en\'eralise celle de Connes et Moscovici \cite{CM95} au cas bivariant. D'autres renormalisations sont bien entendu possibles, y compris sans op\'erateur de Dirac (voir le chapitre suivant), le choix optimal \'etant dict\'e par la situation g\'eom\'etrique \`a laquelle on est confront\'e. L'id\'ee de repr\'esenter une classe de cohomologie en prenant le bord d'une s\'erie formelle renormalis\'ee s'inspire des anomalies chirales en th\'eorie quantique des champs. En fait on montre que toute formule locale de l'indice peut se r\'eduire \`a un calcul d'anomalie. Le mat\'eriel expos\'e dans ce chapitre est tir\'e des articles \\

\noindent \cite{P7} D. Perrot: Anomalies and noncommutative index theory, cours donn\'e \`a Villa de Leyva, Colombie (2005), S. Paycha and B. Uribe ed., {\it Contemp. Math.} {\bf 434} (2007) 125-160.\\

\noindent \cite{P9} D. Perrot: Quasihomomorphisms and the residue Chern character, preprint arXiv:0804.1048.\\

Dans le cas de triplets spectraux, nous avions initialement obtenu la relation entre th\'eor\`eme de l'indice et anomalie chirale (corollaire \ref{cta}) en th\`ese de doctorat suivant une approche diff\'erente. La d\'emonstration n'\'etait pas directement reli\'ee \`a la renormalisation mais reposait sur des consid\'erations g\'eom\'etriques dans l'espace des potentiels de jauge, analogues \`a l'approche d'Atiyah et Singer \cite{AS2, S}. On ne d\'ecrira pas ici ce travail publi\'e dans\\

\noindent \cite{P1} D. Perrot: BRS cohomology and the Chern character in non-commutative geometry, {\it Lett. Math. Phys.} {\bf 50} (1999) 135-144.\\

\noindent \cite{P3} D. Perrot: On the topological interpretation of gravitational anomalies, {\it J. Geom. Phys.} {\bf 39} (2001) 81-95.

\section{Le principe}\label{sprin}

La strat\'egie utilis\'ee pour construire des repr\'esentants \emph{locaux} du caract\`ere de Chern d'un quasihomomorphisme $\rho:\Ac\to \Ec^s\triangleright \Ic^s\hotimes\Bc$ est bas\'ee sur une formule de transgression \cite{P9}. Dans tout ce qui suit on consid\`ere que le quasihomomorphisme est $(p+1)$-sommable et de parit\'e $i\equiv p\mod 2$. D'apr\`es la proposition \ref{padic}, on sait que la caract\`ere de Chern est repr\'esent\'e, en tout degr\'e $n\geq p$, $n\equiv i\mod 2$, par les classes $\ch^n(\rho)\in HC^n(\Ac,\Bc)$ d\'efinies au moyen d'une composition de morphismes 
$$
\ch^n(\rho)\ :\ X(\Th\Ac) \stackrel{\gamma}{\longrightarrow} \Omh \Th\Ac \stackrel{\chih^n}{\longrightarrow} X(\Rch)
$$
avec $0\to \Jc\to \Rc\to \Bc\to 0$ une extension quasi-libre et $\chih^n$ la formule non-locale (\ref{chi}). L'\'egalit\'e $\ch^{n+2}(\rho)\equiv S\ch^n(\rho)$ dans $HC^{n+2}(\Ac,\Bc)$ provient d'une relation entre les cocycles de degr\'es successifs $\chih^{n}$ et $\chih^{n+2}$. Plus pr\'ecis\'ement nous avons introduit dans \cite{P8} une \emph{cocha\^{\i}ne \^eta} $\etah^{n+1}\in \hom(\Omh \Th\Ac,X(\Rch))$ pour tout $n\geq p$ de parit\'e $i\mod 2$. Elle s'annule sur les espaces $\Om^k\Th\Ac$ pour $k\neq n+1$ et $n+2$, tandis que ses deux composantes non nulles $\etah^{n+1}_0:\Om^{n+1}\Th\Ac\to \Rch$ et $\etah^{n+1}_1:\Om^{n+2}\Th\Ac\to \Om^1\Rch_{\nat}$ sont d\'efinies par les formules
\beq
\etah^{n+1}_0(x_0\dd x_1\ldots\dd x_{n+1}) &=& \frac{\Gamma(\frac{n}{2}+1)}{(n+2)!} \, \frac{1}{2}\tau\Big(F x_0[F,x_1]\ldots [F,x_{n+1}] + \non\\
&& \sum_{i=1}^{n+1} (-)^{(n+1)i} [F,x_i]\ldots [F,x_{n+1}]Fx_0 [F,x_1]\ldots [F,x_{i-1}] \Big) \non
\eeq
\beq
\lefteqn{\etah^{n+1}_1(x_0\dd x_1\ldots\dd x_{n+2}) =}\label{eta}\\
&& \frac{\Gamma(\frac{n}{2}+1)}{(n+3)!} \sum_{i=1}^{n+2}  \frac{1}{2}\tau\nat\Big( (ix_0 F + (n+3-i)Fx_0)[F,x_1]\ldots\dd x_i \ldots [F,x_{n+2}]\Big)\non
\eeq
Un calcul direct montre la relation de transgression $\chih^{n}-\chih^{n+2}= [\d, \etah^{n+1}]$ dans le complexe $\hom(\Omh \Th\Ac,X(\Rch))$. On en d\'eduit le r\'esultat suivant.

\begin{proposition}[\cite{P8}]
Soit $\rho:\Ac\to\Ec^s\triangleright\Ic^s\hotimes\Bc$ un quasihomomorphisme $(p+1)$-sommable de parit\'e $i\equiv p\mod 2$, admissible relativement \`a une extension quasi-libre de $\Bc$, et soit un entier $n\geq p$, $n\equiv i\mod 2$. Alors la cochaine d\'efinie par la compos\'ee
\be
\tch^{n+1}(\rho)\ :\ X(\Th\Ac) \stackrel{\gamma}{\longrightarrow} \Omh \Th\Ac \stackrel{\etah^{n+1}}{\longrightarrow} X(\Rch)
\ee
v\'erifie la relation de transgression $\ch^n(\rho)-\ch^{n+2}(\rho)=[\d, \tch^{n+1}(\rho)]$ dans le complexe $\hom(X(\Th\Ac),X(\Rch))$. \cqfd
\end{proposition}

Pour $n<p$ la cochaine $\etah^{n+1}$ n'existe pas car le produit des commutateurs $[F,x]$ n'appartient pas au domaine de la supertrace $\tau$. Supposons cependant que l'on parvienne \`a \'etendre le domaine de $\tau$ par une m\'ethode quelconque, quitte \`a perdre ses propri\'et\'es de supertrace. Par exemple, on peut introduire artificiellement un op\'erateur r\'egularisant dans les formules (\ref{eta}). On d\'efinit ainsi les composantes de la \emph{cocha\^{\i}ne \^eta renormalis\'ee} $\etah_R^{n+1} \in \hom(\Omh \Th\Ac,X(\Rch))$ en tout degr\'e $n<p$ de parit\'e $i\mod 2$. Nous avons introduit dans \cite{P9} la s\'erie 
\be
\eta_R = \sum_{n<p} \etah_R^{n+1} + \sum_{n\geq p} \etah^{n+1}\ , \label{etar}
\ee
vue comme application lin\'eaire de la \emph{somme} directe $\Om \Th\Ac=\bigoplus_n\Om^n\Th\Ac$ vers $X(\Rch)$. Elle ne peut pas s'\'etendre en une cocha\^{\i}ne sur le \emph{produit} direct $\Omh \Th\Ac=\prod_n\Om^n\Th\Ac$ puisque ses composantes $\etah^{n+1}$ ne s'annulent pas pour $n$ suffisamment grand. Selon la terminologie de Higson \cite{Hi} $\eta_R$ est donc une cocha\^{\i}ne impropre. La relation de transgression $\chih^{n}-\chih^{n+2}= [\d, \etah^{n+1}]$ valide pour $n\geq p$ montre que son bord $\chi_R=[\d,\eta_R]$, bien d\'efini dans le $\hom$-complexe $\hom(\Om \Th\Ac,X(\Rch))$, a seulement un nombre fini de composantes non nulles et s'\'etend en un \emph{cocycle} g\'en\'eralement non-trivial dans $\hom(\Omh \Th\Ac,X(\Rch))$. Par exemple dans le cas pair $i=0$ on obtient graphiquement \be
\vcenter{\xymatrix@!0@=2.5pc{
\eta_R = \ar[d]_{[\d,\ ]} &  & \etah_R^1 \ar[dl] \ar[dr] & + & \etah_R^3 \ar[dl] \ar[dr] & \ldots  & + & \etah_R^{p-1} \ar[dl] \ar[dr] & + & \etah^{p+1} \ar[dl] \ar[dr] & + & \etah^{p+3} \ar[dl] \ar[dr] & +\  \ldots  \\
\chi_R  = & \chi_R^0 & + & \chi_R^2 & +& \stackrel{{}}{\ldots}  & \chi_R^{p-2} & + & \chi_R^{p} & + & 0 & + & \ 0\  \ldots }} \label{ser}
\ee
o\`u en chaque degr\'e $n$ pair la cocha\^{\i}ne $\chi_R^n$ poss\`ede deux composantes non nulles $\chi_{R0}^{n}:\Om^{n}T\Ac\to \Rc$ et $\chi_{R1}^{n}:\Om^{n+1}T\Ac\to \Om^1\Rc_{\nat}$. Le point crucial est le lemme suivant, dont la d\'emonstration est imm\'ediate.

\begin{lemma}[\cite{P9}]\label{lan}
Pour tout choix de cocha\^{\i}ne \^eta renormalis\'ee, le bord $\chi_R=[\d,\eta_R]$ est un cocycle cohomologue \`a $\chih^n$ dans le complexe $\hom(\Omh \Th\Ac,X(\Rch))$ pour tout $n\geq p$, $n \equiv i\mod 2$. Par cons\'equent la compos\'ee
\be
\ch_R(\rho)\ :\ X(\Th\Ac) \stackrel{\gamma}{\longrightarrow} \Omh \Th\Ac \stackrel{\chi_R}{\longrightarrow} X(\Rch)
\ee
repr\'esente le caract\`ere de Chern du quasihomomorphisme en cohomologie cyclique bivariante p\'eriodique $\ch(\rho)\in HP^i(\Ac,\Bc)$. \cqfd
\end{lemma}
Dans certaines circonstances (voir le th\'eor\`eme \ref{tres}) on peut aussi contr\^oler finement les propri\'et\'es adiques de $\ch_R(\rho)$ et d\'eterminer sa classe de cohomologie cyclique dans $HC^n(\Ac,\Bc)$, $n\equiv i\mod 2$, qui d\'epend a priori de la renormalisation choisie contrairement \`a son image dans $HP^i(\Ac,\Bc)$. On appellera le cocycle $\chi_R$ une \emph{anomalie} car il est obtenu comme bord de la s\'erie renormalis\'ee (\ref{etar}) et g\'en\'eralise un ph\'enom\`ene bien connu en th\'eorie quantique des champs \cite{P7}. Pour cette raison $\ch_R(\rho)$ est automatiquement donn\'e par une formule locale.
\begin{remark}\textup{En pratique il arrive souvent que le proc\'ed\'e de renormalisation ne soit d\'efini que sur une sous-alg\`ebre dense $\Ac_0\subset \Ac$. Le caract\`ere de Chern renormalis\'e $\ch_R(\rho)$ est alors dans la cohomologie cyclique bivariante $HP^i(\Ac_0,\Bc)$. On doit en tenir compte dans la formulation de th\'eor\`emes de l'indice locaux.}
\end{remark}

Le lien pr\'ecis avec la th\'eorie quantique des champs appara\^it dans la situation suivante. Soit $\rho:\Ac\to\Ec^s\triangleright\Ic^s\hotimes\Bc$ un quasihomomorphisme $(p+1)$-sommable de parit\'e $i=0$. Le th\'eor\`eme de Grothendieck-Riemann-Roch \ref{trr} restreint aux invariants primaires implique la commutativit\'e du diagramme 
\be
\vcenter{\xymatrix{
\Kt_j(\Ic\hotimes\Ac) \ar[r]^{\rho_!} \ar[d] \ar@{.>}[rd]^{\Delta} & \Kt_j(\Ic\hotimes\Bc) \ar[d] \\
HP_j(\Ac) \ar[r]_{\ch(\rho)} & HP_j(\Bc) }} \qquad j\in\zz_2 \label{diag} 
\ee
et l'on se donne comme objectif de fournir une formule explicite pour l'application diagonale $\Delta$. Ainsi la donn\'ee d'une classe de cohomologie cyclique sur $\Bc$ permet de construire un invariant de la $K$-th\'eorie topologique de $\Ac$. C'est g\'en\'eralement sous cette forme que sont \'enonc\'es les th\'eor\`emes de l'indice sup\'erieurs (voir par exemple \cite{CM90}, ou le chapitre \ref{ceq}). Par p\'eriodicit\'e de Bott, on peut toujours se ramener aux groupes de $K$-th\'eorie impairs ($j=1$). Donc si $u\in (\Ic\hotimes\Ac)^+$ est un \'el\'ement inversible repr\'esentant une classe $[u]\in \Kt_1(\Ic\hotimes\Ac)$, son image par la diagonale est repr\'esent\'ee par le cycle impair du $X$-complexe associ\'e \`a une extension quasi-libre $0\to \Jc\to \Rc\to \Bc\to 0$ (le facteur $\sqrt{2\pi i}$ est n\'ecessaire pour assurer la compatibilit\'e avec la p\'eriodicit\'e de Bott)
\be
\Delta(u) =\sqrt{2\pi i}\, \ch_R(\rho)\cdot\ch(u)\ \in\Om^1\Rch_{\nat}\ ,\label{dia}
\ee
pour n'importe quel choix de renormalisation de la cocha\^{\i}ne \^eta, avec $\ch_R(\rho)=[\d,\eta_R]\gamma$. A l'aide des transgressions $\tch_R^{2n+1}(\rho)=\etah_R^{2n+1}\gamma$, on peut donc \'ecrire (\ref{dia}) comme l'image d'une s\'erie formelle (non convergente) dans $\Rch$ sous l'application de bord $\nat\dd: \Rch\to \Om^1\Rch_{\nat}$
$$
\Delta(u)= \sqrt{2\pi i}\,\nat\dd \sum_{n=0}^{\infty} \tch_R^{2n+1}(\rho)\cdot\ch(u)\ ,
$$
\`a condition d'organiser la compensation d'une infinit\'e de termes de la m\^eme mani\`ere que dans le sch\'ema (\ref{ser}). Nous avons montr\'e dans \cite{P9} comment proc\'eder, en interpr\'etant la s\'erie formelle comme d\'eveloppement en puissances d'un \emph{potentiel de jauge}. Sa repr\'esentation graphique reproduit exactement le d\'eveloppement en graphes de Feynman d'une th\'eorie de jauge non-commutative, tandis que $\Delta(u)$ est l'anomalie chirale associ\'ee \`a la transformation de jauge $u$. Avant de d\'ecrire ce calcul plus en d\'etail voyons un exemple pratique de renormalisation par fonction z\^eta.

\section{Renormalisation z\^eta}

Afin d'illustrer la ``localit\'e'' du cycle $\chi_R=[\d,\eta_R]$, on d\'eveloppe ici un calcul pseudodiff\'erentiel abstrait construit sur un op\'erateur de Dirac \cite{P9}. Le r\'esultat obtenu en termes de r\'esidus de fonctions z\^eta est une g\'en\'eralisation bivariante de la formule locale de Connes et Moscovici pour le caract\`ere de Chern des triplets spectraux \cite{CM95}. \\

Soient $\Ac$ et $\Bc$ deux $m$-alg\`ebres de Fr\'echet et soit $(H,\rho,F)$ un $\Ac$-$\Bc$-bimodule born\'e mis sous la forme d'un quasihomomorphisme, avec $H$ espace de Hilbert $\zz_2$-gradu\'e. On se donne en plus un endomorphisme non born\'e $|D|$ de degr\'e pair sur $H_{\Bc}$ et commutant avec $F$. Afin de d\'evelopper un calcul pseudodiff\'erentiel adapt\'e \`a la situation bivariante, on suppose en plus la donn\'ee d'une \emph{alg\`ebre de symboles abstraits} (\cite{P9}), d\'efinie comme limite inductive 
\be
\Pc= \varinjlim_{\al\in\rr} \Pc_{\al}\ ,
\ee
o\`u les $\Pc_{\al}$ sont des espaces de Fr\'echet $\zz_2$-gradu\'es avec injections continues $\Pc_{\al}\to \Pc_{\beta}$ pour $\al\leq\beta$, v\'erifiant les propri\'et\'es suivantes: 
\begin{itemize}
\item Chaque $\Pc_{\al}$ est un espace d'op\'erateurs non-born\'es sur $H$ pour $\al>0$, born\'es pour $\al\leq 0$, et $F\in \Pc_0$.
\item $\Pc$ est une alg\`ebre filtr\'ee avec unit\'e, c'est-\`a-dire munie d'un produit associatif continu $\Pc_{\al}\times \Pc_{\beta} \to \Pc_{\al+\beta}$ et $1\in \Pc_0$. En particulier $\Pc_{\al}$ est une alg\`ebre de Fr\'echet pout $\al\leq 0$.
\item L'homomorphisme $\rho:\Ac\to \End(H_{\Bc})$ se factorise \`a travers le produit tensoriel projectif $\Pc_0\hotimes\Bc$.
\item Le commutateur $[F,\rho(a)]$ est dans l'id\'eal $\Pc_{-1}\hotimes\Bc\subset \Pc_0\hotimes\Bc$ pour tout $a\in\Ac$.
\item $|D|\in \Pc_1\hotimes\Bc^+$, et son spectre en tant qu'\'el\'ement de l'alg\`ebre $\Pc\hotimes\Bc^+=\varinjlim_{\al} \Pc_{\al}\hotimes\Bc^+$ est contenu dans un intervalle r\'eel $[\eps,+\infty)$, $\eps>0$.
\item Pour $\la\in\cc$ hors du spectre la r\'esolvante $(\la-|D|)^{-1}$ est dans $\Pc_{-1}\hotimes\Bc^+$ et l'application $\la \mapsto (\la-|D|)^{-1}$ est une fonction holomorphe et born\'ee sur tout demi-plan $\re(\la)\leq \eps'<\eps$ disjoint du spectre de $|D|$.
\end{itemize}
On notera bri\`evement $(H,\rho,F,|D|)$ un tel bimodule muni d'une alg\`ebre de symboles abstraits. Dans les exemples pratiques $\Pc$ est r\'ealis\'ee comme une alg\`ebre d'op\'erateurs born\'es entre espaces de Sobolev \cite{P9}; l'utilisation du mot ``symbole'' est donc largement arbitraire puisqu'il ne s'agit pas n\'ecessairement de symboles d'op\'erateurs pseudodiff\'erentiels. D'autre part, le fait d'imposer que $\rho(\Ac)$ et $|D|$ font partie d'une alg\`ebre produit tensoriel $\Pc\hotimes\Bc^+$ est commode mais \'evidemment restrictif. Il est certainement possible de consid\'erer des alg\`ebres plus larges \`a condition d'adapter la discussion ci-dessous. Fixons maintenant une extension $0\to \Jc\to \Rc\to \Bc\to 0$ de $m$-alg\`ebres de Fr\'echet, avec $\Rch$ la compl\'etion $\Jc$-adique de $\Rc$. Alors $|D|$ se rel\`eve en un op\'erateur $|\Dh|\in \Pc_1\hotimes\Rch^+$, et sa r\'esolvante est bien d\'efinie:
\be
(\la - |\Dh|)^{-1}\in \Pc_{-1}\hotimes\Rch^+\ .
\ee
De m\^eme $\rho$ se rel\`eve en un homomorphisme $\rho_*: \Th\Ac\to \Pc_0\hotimes\Rch$. Les puissances complexes $|\Dh|^z\in \Pc_{\re(z)}\hotimes\Rch^+$ sont alors obtenues par une int\'egrale de contour \cite{P9}. On veut ensuite contr\^oler les commutateurs embo\^{\i}t\'es de $|\Dh|$ et de l'op\'erateur de Dirac $\Dh=F|\Dh|$ avec l'image de $\Th\Ac$ dans $\Pc\hotimes\Rch$. Introduisons \`a cette fin la filtration de $\Pc\hotimes\Rch^+$ par les sous-espaces
$$
\Fc_{\al}^k= \Pc_{\al}\hotimes \Jch^k\ ,\qquad \Fc_{\al}= \sum_{k\geq 0} \Fc_{\al+k}^k \ ,\qquad \al\in\rr\ , \ k\in\nn\ ,
$$
o\`u $\Jch^k$ est la compl\'etion $\Jc$-adique de $\Jc^k$. On a $\Th\Ac\subset \Fc_0^0$, $\Jh\Ac\subset \Fc_0^1$ et pour tout $z\in\cc$, $|\Dh|^z\in \Fc_{\re(z)}^0$. La d\'efinition suivante est une g\'en\'eralisation bivariante de la condition de r\'egularit\'e introduite dans le cadre des triplets spectraux par Connes et Moscovici \cite{CM95}.

\begin{definition}[\cite{P9}]\label{dreg}
Un quasihomomorphisme muni d'une alg\`ebre de symboles abstraits est \emph{r\'egulier} relativement \`a l'extension $0\to \Jc\to \Rc\to \Bc\to 0$ si les puissances de la d\'erivation $\delta=[|\Dh|, \ ]$ sur l'alg\`ebre $\Pc\hotimes\Rch$ v\'erifient, pour tout $n\geq 0$ et $k\geq 0$,
\beq
\delta^n(\Th\Ac)+ \delta^n [\Dh,\Th\Ac] &\subset& \Fc_0^0 + \Fc_1^1 + \ldots + \Fc_n^n\ \subset \ \Fc_0 \ ,\non\\
\delta^n((\Jh\Ac)^k)+ \delta^n [\Dh,(\Jh\Ac)^k] &\subset& \Fc_0^k + \Fc_1^{k+1} + \ldots + \Fc_n^{k+n}\ \subset \ \Fc_{-k}\ . \non
\eeq
\end{definition}
La r\'egularit\'e implique que la d\'erivation $\delta$ ne peut augmenter le degr\'e symbolique $\al$ que si elle augmente simultan\'ement le degr\'e adique $k$. En ce sens l'image de $\Th\Ac$ est ``lisse'' dans $\Pc\hotimes\Rch$. Il faut noter que cette condition impose de fortes contraintes sur le choix de l'extension $\Rc$. En particulier l'alg\`ebre tensorielle $\Rc=T\Bc$ est rarement compatible avec la r\'egularit\'e. \\

On peut maintenant d\'efinir l'alg\`ebre des op\'erateurs pseudodiff\'erentiels abstraits $\Psi\Ac$ comme la sous-alg\`ebre (non compl\`ete) de $\Pc\hotimes\Rch^+$ engendr\'ee par l'image de $\Th\Ac$, l'op\'erateur $F$ et toutes les puissances complexes $|\Dh|^z$, $z\in\cc$. En utilisant la condition de r\'egularit\'e, on voit que tout op\'erateur pseudodiff\'erentiel est combinaison lin\'eaire d'\'el\'ements $x\in\Psi\Ac$ ayant un d\'eveloppement asymptotique
$$
x\simeq \sum_{k\geq 0}(a_k+b_k F)\, |\Dh|^{z-k}\ ,\quad z\in\cc\ ,
$$
o\`u $a_k$ et $b_k$ sont dans l'alg\`ebre engendr\'ee par les d\'eriv\'ees $\delta^n(\Th\Ac)$ et $\delta^n[\Dh,\Th\Ac]$, et l'\'egalit\'e $\simeq$ signifie que pour tout $N\geq 0$, la diff\'erence $x-\sum_{k= 0}^{N}(a_k+b_k F)\, |\Dh|^{z-k}$ est dans $\Fc_{\re(z)-N-1}$. On peut alors filtrer l'alg\`ebre $\Psi\Ac$ par les sous-espaces $(\Psi\Ac)^k_{\al}=\Psi\Ac\cap \Fc^k_{\al}$. De mani\`ere analogue le bimodule des 1-formes non-commutatives $\Om^1\Psi\Ac$ est filtr\'e par des sous-espaces $(\Om^1\Psi\Ac)^k_{\al}$ et l'on obtient une famille de $X$-complexes
$$
X(\Psi\Ac)_{\al}^k\ :\ (\Psi\Ac)_{\al}^k\ \rightleftarrows\ \nat(\Om^1\Psi\Ac)_{\al}^k
$$ 
index\'es par le degr\'e symbolique $\al\in\rr$ et le degr\'e adique $k\in\nn$ des op\'erateurs pseudodiff\'erentiels impliqu\'es. L'id\'ee est que pour $\al$ suffisamment n\'egatif, $X(\Psi\Ac)_{\al}^k$ soit dans le domaine de la supertrace d'op\'erateurs sur $H$. On peut alors tenter de saturer le facteur $\Pc_{\al}$ par la supertrace et obtenir ainsi un morphisme $X(\Psi\Ac)_{\al}^k \to X(\Rch)$, compatible avec la filtration de $X(\Rch)$ associ\'ee aux sous-complexes $F^{2k-1}\Xh(\Rc,\Jc): \Jch^k\rightleftarrows \nat(\Jch^k\dd\Rch+\Jch^{k-1}\dd\Jch)$. Cela conduit \`a la d\'efinition suivante, qui g\'en\'eralise la notion de spectre de dimension des triplets spectraux \cite{CM95}.

\begin{definition}[\cite{P9}]\label{dsd}
Un $\Ac$-$\Bc$-bimodule $(H,\rho,F,|D|)$ r\'egulier relativement \`a une extension $0\to \Jc\to \Rc\to \Bc\to 0$ a une \emph{dimension analytique finie} $p\in\rr$ si la supertrace d'op\'erateurs sur $H$ d\'efinit un morphisme de complexes 
$$
\tau: X(\Psi\Ac)_{\al}^k \to F^{2k-1}\Xh(\Rc,\Jc)
$$
pour tout $\al<-p$ et $k\in\nn$. De plus le quasihomomorphisme a un \emph{spectre de dimension discret} si pour tous op\'erateurs pseudodiff\'erentiels $x,y\in (\Psi\Ac)_0^k$, les fonctions z\^eta
$$
\tau(x|\Dh|^{-z}y)\in \Rch\ ,\qquad \left. \begin{array}{c}
\tau\nat(x|\Dh|^{-z}\dd y) \\
 \tau\nat(x|\Dh|^{-z-1}y\dd |\Dh|) \end{array} \right\} \in \Om^1\Rch_{\nat}
$$
sont holomorphes sur le demi-plan complexe $\re(z) > p$ et admettent un prolongement m\'eromorphe sur le compl\'ementaire d'un sous-ensemble discret de $\cc$. 
\end{definition}
Dans la suite on consid\`ere que le bimodule $(H,\rho,F,|D|)$ est $(p+1)$-sommable, de dimension analytique $p$ et de parit\'e $i\equiv p\mod 2$. Les puissances complexes $|\Dh|^{-z}$ avec $\re(z)\gg 0$ permettent alors de construire les composantes renormalis\'ees $\etah_R^{n+1}:\Omh\Th\Ac\to X(\Rch)$ de la cochaine \^eta en degr\'es $n<p$ de parit\'e $i$. Posons
\beq
\lefteqn{\etah^{n+1}_{R0}(x_0\dd x_1\ldots\dd x_{n+1}) =\frac{\Gamma(\frac{n}{2}+1)}{(n+2)!}\,  \frac{1}{2} \Pf_{z=0}\tau\Big( |\Dh|^{-z}F x_0[F,x_1]\ldots [F,x_{n+1}]+ } \non \\
&&\qquad\qquad    \sum_{i=1}^{n+1}(-)^{(n+1)i} [F,x_i]\ldots [F,x_{n+1}] |\Dh|^{-z}Fx_0 [F,x_1]\ldots [F,x_{i-1}] \Big) \non
\eeq
\beq
\lefteqn{\etah^{n+1}_{R1}(x_0\dd x_1\ldots\dd x_{n+2}) =\frac{\Gamma(\frac{n}{2}+1)}{(n+3)!}\frac{1}{2}\times   }\label{zeta}\\
&&\sum_{i=1}^{n+2}  \Pf_{z=0} \Big( \sum_{j=1}^i \tau\nat x_0F[F,x_1]\ldots [F,x_{j-1}]|\Dh|^{-z}[F,x_j] \ldots \dd x_i \ldots [F,x_{n+2}] \non\\
&&\qquad  +\sum_{j=i}^{n+2} \tau\nat Fx_0[F,x_1]\ldots \dd x_i \ldots [F,x_j]|\Dh|^{-z}[F,x_{j+1}] \ldots [F,x_{n+2}] \Big) \ ,\non
\eeq
o\`u $\Pf_{z=0}$ est la partie finie des fonctions z\^eta, c'est-\`a-dire le terme constant dans leur d\'eveloppement en s\'erie de Laurent autour de $z=0$. Lorsque $n\geq p$ les fonctions z\^eta n'ont pas de p\^ole en z\'ero et l'on retrouve les formules donn\'ees initialement pour $\etah^{n+1}$. Notons que la renormalisation (\ref{zeta}) est loin d'\^etre unique: en effet on pourrait changer de position l'op\'erateur $|\Dh|^{-z}$ ou bien multiplier les fonctions z\^eta par une fonction $h(z)$ holomorphe autour de z\'ero avec $h(0)=1$. Dans tous les cas les nouvelles composantes $\etah_R^{n+1}$ ne changeraient que par des \emph{sommes de r\'esidus de fonctions z\^eta} \cite{P9}, autrement dit des termes ``locaux''. Dans un langage de th\'eorie quantique des champs ces r\'esidus correspondent \`a des \emph{contre-termes} effectuant le passage d'une renormalisation \`a une autre. Le lien pr\'ecis entre th\'eorie locale de l'indice non-commutative et th\'eorie des champs sera tra\^{\i}t\'e dans le paragraphe suivant.\\
Supposons que l'alg\`ebre $\Rc$ est quasi-libre. D'apr\`es le lemme \ref{lan} on sait que le bord de la s\'erie $\eta_R=\sum_{n<p}\etah_R^{n+1}+ \sum_{n\geq p}\etah^{n+1}$ repr\'esente le caract\`ere de Chern en cohomologie cyclique bivariante p\'eriodique $\ch(\rho)\in HP^*(\Ac,\Bc)$. Nous avons montr\'e dans \cite{P9} qu'il est n\'ecessairement donn\'e par une somme de r\'esidus. De surcro\^{\i}t, la renormalisation sp\'ecifique (\ref{zeta}) est judicieusement choisie de fa\c{c}on \`a contr\^oler aussi sa classe de cohomologie cyclique bivariante non-p\'eriodique.

\begin{theorem}[\cite{P9}]\label{tres}
Soit $(H,\rho,F,|D|)$ un $\Ac$-$\Bc$-bimodule $(p+1)$-sommable de parit\'e $i\equiv p\mod 2$ muni d'une alg\`ebre de symboles abstraits. On le suppose r\'egulier relativement \`a une extension quasi-libre $0\to \Jc\to \Rc\to \Bc\to 0$, de dimension analytique $p$ et de spectre de dimension discret. Alors le bord de sa cocha\^{\i}ne \^eta renormalis\'ee par (\ref{zeta}) repr\'esente le caract\`ere de Chern en cohomologie cyclique bivariante non-p\'eriodique 
$$
\ch^p(\rho)\equiv [\d,\eta_R]\circ \gamma \ \in HC^p(\Ac,\Bc)\ .
$$
De plus $[\partial,\eta_R]$ est cohomologue, dans le complexe $\hom(\Omh\Th\Ac, X(\Rch))$, au cocycle $\chi_R$ dont les composantes $\chi^n_{R0}:\Omh^n\Th\Ac\to \Rch$ et $\chi^n_{R1}:\Omh^{n+1}\Th\Ac\to \Om^1\Rch_{\nat}$  sont d\'efinies en tout degr\'e $n\equiv i\mod 2$ par une somme de r\'esidus
\beq
\lefteqn{\chi_{R0}^{n}(x_0\dd x_1\ldots\dd x_{n}) =\sum_{k_0,\ldots,k_{n}\geq 0} (-)^{k+n}c(k_0,\ldots,k_{n})\, \res \Big(\frac{\Gamma(z+k+\frac{n}{2})}{\Gamma(z+1)}\times}\non\\ 
&&\qquad\qquad  \sum_{i=0}^{n} (-)^{i(n-1)}\, \tau(dx_{n-i+1}^{(k_0)}\ldots x_0^{(k_i)}dx_1^{(k_{i+1})}\ldots dx_{n-i}^{(k_{n})}|\Dh|^{-2(z+k)-n})\Big)\non
\eeq
\beq
\lefteqn{\chi_{R1}^{n}(x_0\dd x_1\ldots\dd x_{n+1}) =\sum_{k_0,\ldots,k_{n}\geq 0} (-)^{k+n}c(k_0,\ldots,k_{n})\, \res \Big(\frac{\Gamma(z+k+\frac{n}{2})}{\Gamma(z+1)}\times}\non\\
&&\qquad \sum_{i=0}^{n}(-)^{in} \tau\nat (dx_{n-i+2}^{(k_0)}\ldots x_0^{(k_i)}dx_1^{(k_{i+1})}\ldots dx_{n-i}^{(k_{n})}|\Dh|^{-2(z+k)-n}\dd x_{n-i+1})\Big)\non\\
&&\qquad\qquad - \sum_{k_0,\ldots,k_{n+1}\geq 0} (-)^{k+n}c(k_0,\ldots,k_{n+1})\, \res \Big(\frac{\Gamma(z+k+\frac{n}{2}+1)}{\Gamma(z+1)}\times\non\\
&&\qquad \sum_{i=0}^{n+1}(-)^{in} \tau\nat (dx_{n-i+2}^{(k_0)}\ldots x_0^{(k_i)}dx_1^{(k_{i+1})}\ldots dx_{n-i+1}^{(k_{n+1})}|\Dh|^{-2(z+k+1)-n}\dd \Dh)\Big)\non
\eeq
avec $k=\sum_i k_i$, $dx=[\Dh,x]$, $x^{(k_i)}$ est la $k_i$-i\`eme d\'eriv\'ee de $x$ par rapport au commutateur $[\Dh^2,\ ]$, et la constante $c(k_0,\ldots,k_{n})$ est donn\'ee par
$$
c(k_0,\ldots,k_{n})^{-1}=k_0!\ldots k_n! (k_0+1)(k_0+k_1+2)\ldots(k_0+\ldots+k_n+n+1)\ .
$$
En cons\'equence le cocycle $\ch_R=\chi_R\circ\gamma$ repr\'esente le caract\`ere de Chern en cohomologie cyclique bivariante p\'eriodique $\ch(\rho) \in HP^i(\Ac,\Bc)$.  \cqfd
\end{theorem}
A priori $\chi_R$ poss\`ede une infinit\'e de composantes $\chi_R^n$ non nulles. Cependant apr\`es projection de $X(\Rch)$ sur un complexe quotient $X(\Rc/\Jc^k)$ seul un nombre fini de composantes survivent \cite{P9}, et $\chi_R\in \hom(\Omh\Th\Ac,X(\Rch))$ donne r\'eellement un cocycle cyclique bivariant p\'eriodique. Dans le cas particulier $\Rc=\Bc=\cc$ on retrouve d'ailleurs la formule locale de Connes et Moscovici pour le caract\`ere de Chern des triplets spectraux r\'eguliers \cite{CM95}, avec $\chi_R^n=0$ d\`es que $n>p$.

\begin{remark}\textup{La formule de r\'esidus donn\'ee pour $\chi_R$ est en fait intimement li\'ee au cocycle (\ref{map}) construit autour du noyau de la chaleur. Plus pr\'ecis\'ement $\chi_R$ est extrait du d\'eveloppement asymptotique de ce cocycle \`a la limite $t\downarrow 0$. Pour cette raison le th\'eor\`eme \ref{tres} permet de red\'emontrer la formule de localisation entrant dans le th\'eor\`eme de l'indice \'equivariant du chapitre \ref{ceq}, au moins dans le cas d'un groupe $G$ discret (\cite{P9}). }
\end{remark}

\section{Triplets spectraux et anomalies}\label{sreg}

Nous allons maintenant utiliser la renormalisation z\^eta dans le cas particulier des triplets spectraux et constater que la formule locale $\chi_R=[\d,\eta_R]$ correspond exactement \`a l'anomalie chirale d'une th\'eorie de jauge non-commutative \cite{P7}. \\

On consid\`ere un triplet spectral $(p+1)$-sommable $(\Ac,H,D)$ de parit\'e $i=0$, avec $p$ entier pair et $\Ac$ une $m$-alg\`ebre de Fr\'echet involutive $*$-repr\'esent\'ee par op\'erateurs born\'es sur $H$. Sans perte de g\'en\'eralit\'e on peut supposer que l'op\'erateur de Dirac $D$ est inversible. Dans une d\'ecomposition de l'espace de Hilbert $H=H_+ \oplus H_-$ correspondant \`a sa $\zz_2$-graduation, la repr\'esentation de $\Ac$ et l'op\'erateur de Dirac s'\'ecrivent
$$
a=\left(\begin{matrix}
a_+ & 0 \\
0 & a_- \end{matrix} \right)\ ,\qquad
D=\left(\begin{matrix}
0 & Q^* \\
Q & 0 \end{matrix} \right)\ .
$$
L'objectif est de calculer l'indice de l'op\'erateur de Dirac \`a coefficients dans un projecteur $e=e^*=e^2\in M_{\infty}(\Ac)$ repr\'esentant une classe de $K$-th\'eorie $[e]\in \Kt_0(\Ac)$. En modifiant la repr\'esentation de $\Ac$ le triplet spectral se remet sous la forme d'un quasihomomorphisme $\rho:\Ac\to \Ec^s\triangleright \Ic^s$, avec $\Ec=\End(H_+)$ et $\Ic=\ell^{p+1}(H_+)$, muni d'un module de Dirac $|D|$:
$$
\rho(a)=\left(\begin{matrix}
a_+ & 0 \\
0 & Q^{-1}a_-Q \end{matrix} \right)\ ,\quad
F=\left(\begin{matrix}
0 & 1 \\
1 & 0 \end{matrix} \right)\ ,\quad
|D|= \Id(\cc^2)\otimes (Q^*Q)^{1/2}\ .
$$
Soit $S\Ac=\Ac\hotimes \cinf(0,1)$ la suspension lisse de $\Ac$ et soit $\Bc=\cinf(S^1)$ l'alg\`ebre du cercle. $\rho$ s'\'etend de mani\`ere \'evidente en un quasihomomorphisme de $S\Ac$ vers $\Bc$, et en vertu de l'isomorphisme de Bott $\Kt_0(\Ac)\cong\Kt_1(S\Ac)$ l'indice de l'op\'erateur de Dirac $eDe$ correspond \`a l'image de $[e]$ sous la diagonale du diagramme commutatif 
\be
\vcenter{\xymatrix{
\Kt_1(S\Ac) \ar[r] \ar[d] \ar@{.>}[rd]^{\Delta} & \Kt_1(\Ic\hotimes\Bc)\cong\zz \ar[d] \\
HP_1(S\Ac) \ar[r] & HP_1(\Bc)\cong\cc }} \label{di}
\ee
Nous avons montr\'e dans \cite{P7} comment relier $\Delta$ \`a l'anomalie d'une th\'eorie de jauge chirale non-commutative associ\'ee au triplet spectral. On d\'efinit d'abord un \emph{potentiel de jauge} comme un op\'erateur born\'e $A: H_+\to H_-$ form\'e de combinaisons lin\'eaires finies du type $a_-(Qb_+-b_-Q)$ avec $a,b\in \Ac$. On peut alors coupler $A$ \`a des champs chiraux $\psi\in H_+$ et $\psib\in H_-^*$ par l'interm\'ediaire d'une fonctionnelle d'action 
\be
S(\psi,\psib,A)= \langle \psib, (Q+A)\psi \rangle\ \in\cc\ . \label{action}
\ee
On veut maintenant quantifier les champs $\psi$ et $\psib$ en tant que fermions, en laissant le potentiel $A$ fix\'e \cite{IZ}. Cela revient \`a ajouter une fluctuation quantique $W(A)$ \`a l'action (\ref{action}), obtenue comme logarithme d'un d\'eterminant (voir \cite{P7}, on prend ici la constante de Planck $\hbar=1$)
\be
W(A) = \ln \Det(1+Q^{-1}A)\ .
\ee
Remarquons que l'op\'erateur $Q^{-1}A$ est dans la classe de Schatten $\ell^{p+1}(H_+)$, donc le d\'eterminant n'est pas bien d\'efini lorsque $p>0$ et n\'ecessite une renormalisation. Dans l'approche perturbative on ne s'int\'eresse qu'au d\'eveloppement de $W(A)$ en s\'erie formelle de puissances de $A$, obtenu par la relation na\"{\i}ve $\ln\Det = \Tr\ln$. La repr\'esentation graphique se fait au moyen de diagrammes de Feynman \`a une boucle,
\beq
W(A)  &=& \Tr \ln(1+Q^{-1}A)\ =\ \sum_{n\geq 1}\frac{(-)^{n+1}}{n}\, \Tr((Q^{-1}A)^n) \non\\
&=& \xymatrix{ *=0{\bullet} \ar@(ur,dr)@{-}[] |-{\SelectTips{cm}{}\object@{>}} }   
\ - \frac{1}{2}\  
\vcenter{\xymatrix@R=3pc{ *=0{\bullet} \ar@/^/@{-}[d] |-{\SelectTips{cm}{}\object@{>}} \\*=0{\bullet} \ar@/^/@{-}[u] |-{\SelectTips{cm}{}\object@{>}} } }
\ +\ \frac{1}{3}\ 
\vcenter{\xymatrix@C=1pc{ & *=0{\bullet} \ar@{-}[dr] |-{\SelectTips{cm}{}\object@{>}} &  \\
*=0{\bullet} \ar@{-}[ur] |-{\SelectTips{cm}{}\object@{>}} & & *=0{\bullet} \ar@{-}[ll] |-{\SelectTips{cm}{}\object@{>}} }}
\ -\frac{1}{4}\ 
\vcenter{\xymatrix@R=3pc@C=2.8pc{ *=0{\bullet} \ar@{-}[r] |-{\SelectTips{cm}{}\object@{>}} & *=0{\bullet} \ar@{-}[d] |-{\SelectTips{cm}{}\object@{>}} \\
*=0{\bullet} \ar@{-}[u] |-{\SelectTips{cm}{}\object@{>}} & *=0{\bullet} \ar@{-}[l] |-{\SelectTips{cm}{}\object@{>}} }}
\ +\frac{1}{5}\ 
\vcenter{\xymatrix@R=0.25pc@C=0.5pc{ & & *=0{\bullet} \ar@{-}[ddrr] |-{\SelectTips{cm}{}\object@{>}} & \\
 & & & & \\
*=0{\bullet} \ar@{-}[uurr] |-{\SelectTips{cm}{}\object@{>}} &  &  &  & *=0{\bullet} \ar@{-}[dddl] |-{\SelectTips{cm}{}\object@{>}} \\
 & & & &  \\
 & & & & \\
 & *=0{\bullet} \ar@{-}[uuul] |-{\SelectTips{cm}{}\object@{>}} & & *=0{\bullet} \ar@{-}[ll] |-{\SelectTips{cm}{}\object@{>}} & }}
\ +\ldots \non
\eeq
o\`u chaque sommet repr\'esente une insertion de potentiel $A$ et chaque ar\^ete repr\'esente le propagateur $Q^{-1}$. Puisque l'op\'erateur $Q^{-1}A$ est tra\c{c}able lorsque $n\geq p+1$, les termes $W^n(A)=\frac{(-)^{n+1}}{n}\, \Tr((Q^{-1}A)^n)$ sont bien d\'efinis dans ce cas. Par cons\'equent seuls les premiers termes du d\'eveloppement n\'ecessitent une renormalisation. Supposons que le triplet spectral soit \emph{r\'egulier} et de dimension analytique $p$. On peut alors choisir une renormalisation z\^eta et d\'efinir pour $n\leq p$ 
\be
W_R^n(A) = \frac{(-)^{n+1}}{n} \Pf_{z=0} \Tr((Q^{-1}A)^n |D|_+^{-2z})\label{ren}
\ee
avec $|D|_+=(Q^*Q)^{1/2}$. A partir de maintenant nous allons consid\'erer des familles de potentiels de jauge param\'etr\'ees par le cercle. Soit $u$ un \'el\'ement unitaire tel que $u-1$ soit dans le produit tensoriel \emph{alg\'ebrique} $\Ac\otimes \cinf(0,1)$, qui est une sous-alg\`ebre dense de la suspension $S\Ac$. Il d\'efinit une classe de $K$-th\'eorie topologique $[u]\in \Kt_1(S\Ac)$. Par exemple, on peut prendre l'unitaire $u=1+ e\otimes (\beta-1)$ qui correspond \`a un projecteur $e\in\Ac$ par p\'eriodicit\'e de Bott, avec $\beta$ le g\'en\'erateur de Bott du cercle. Regardons maintenant $u$ comme une boucle de transformations de jauge dans la th\'eorie des champs, \`a point-base l'identit\'e. La famille d'op\'erateurs
\be 
A= u_-^{-1}Qu_+ - Q \label{pot}
\ee
est donc une boucle de potentiels, obtenue en appliquant la transformation de jauge $u$ sur le potentiel trivial $A_0=0$. On note $\dd: \cinf(S^1)\to \Om^1(S^1)$ la diff\'erentielle de de Rham sur le cercle, et $\om=u^{-1}\dd u\in \Ac\otimes \Om^1(S^1)$ la forme de Maurer-Cartan associ\'ee \`a la boucle $u$. Alors la diff\'erentielle de $A$ s'exprime au moyen de $\om$ (on utilise la convention que $A$ et $Q$ sont de degr\'e impair, voir les \'equations BRS \cite{MSZ})
$$
- \dd A = (Q+A)\om_+ + \om_-(Q+A)\ .
$$
La diff\'erentielle $\dd W_R^n(A)$ est alors calculable en fonction de $\om$ et $A$ en tout degr\'e $n$. Elle d\'epend lin\'eairement de $\om$ et se scinde en deux termes de degr\'es respectifs $n-1$ et $n$ par rapport \`a $A$. En sommant la s\'erie formelle $\dd W_R(A)$ en puissances de $A$, on constate (voir \cite{P7}) que les termes de degr\'e $>p$ se compensent deux \`a deux, alors qu'en degr\'e inf\'erieur cette propri\'et\'e est bris\'ee par la renormalisation (\ref{ren}). Le bord $\Delta(\om,A):=\dd W_R(A)$ est donc une \emph{somme finie} de formes diff\'erentielles au-dessus du cercle, qui d\'epend polyn\^omialement de $A$. Par exemple pour $p=2$ on obtient graphiquement
\be
\vcenter{\xymatrix@!0@=2.5pc{
W_R(A) \ar[d]_{\dd} & = &  & W^1_R(A) \ar[dl] \ar[dr]  & + & W^2_R(A) \ar[dl] \ar[dr] & + & W^3(A) \ar[dl] \ar[dr] & + & W^4(A) \ar[dl] \ar[dr] & +  \ldots  \\
\Delta(\om,A) & = & \Delta^0(\om,A) & + & \Delta^1(\om,A) & + & \Delta^2(\om,A) & + & 0 & + & \ 0  \ldots }} \label{gr}
\ee
o\`u $\Delta^n(\om,A)$ est de degr\'e $n$ en $A$. La 1-forme $\Delta(\om,A)$ est l'anomalie chirale associ\'ee \`a la th\'eorie des champs et mesure la brisure d'invariance de jauge de l'action quantique renormalis\'ee. L'anomalie est n\'ecessairement donn\'ee par une formule locale; dans le cas de la renormalisation (\ref{ren}) c'est une somme finie de r\'esidus de fonctions z\^eta \cite{P7}. Notons qu'un changement de renormalisation affecte seulement les premiers termes de la s\'erie formelle $W_R(A)$, en ajoutant une \emph{fonctionnelle locale et polyn\^omiale} $P(A)$ de degr\'e $\leq p$ (contre-termes). L'anomalie correspondante $\Delta'(\om,A)= \Delta(\om,A) +\dd P(A)$ diff\`ere donc d'un cobord et sa classe de cohomologie dans $H^1(S^1)$ est ind\'ependante de la renormalisation choisie. Par exemple en utilisant une renormalisation z\^eta \`a la Ray-Singer \cite{RS, S}, nous avons donn\'e dans \cite{P7} la formule de r\'esidus suivante pour l'anomalie (ici nos conventions de notation et de signe diff\`erent l\'eg\`erement de \cite{P7}):
\beq
\lefteqn{\Delta'(\om,A)= \res\, \frac{1}{z}\tau(\omega |D|^{-2z})\ + } \label{an} \\
&& \sum_{\substack{n\geq 1 \\ k\geq 0}}(-1)^{n+k}c(k) \res \, \frac{\Gamma(z+n+k)}{\Gamma(z+1)}\textup{Tr} \big( q\omega A^{(k_1)}Q^*A^{(k_2)}\ldots Q^*A^{(k_n)}|D|_+^{-2(z+n+k)} \big) \non
\eeq
o\`u $\tau$ est la supertrace d'op\'erateurs sur $H=H_+\oplus H_-$, $k=(k_1,\ldots,k_n)$ est un multi-indice, $q\omega=\omega_+Q^*+Q^*\om_-$, $A^{(k_i)}$ est la $k_i$-\`eme d\'eriv\'ee de $A$ par rapport au commutateur $[Q^*Q,\ ]$, et $c(k)$ est la constante  
$$
c(k)^{-1}=(k_1!\ldots k_n!)(k_1+1)(k_1+k_2+2)\ldots(k_1+\ldots +k_n+n)\ .
$$
Puisque les r\'esidus s'annulent pour des op\'erateurs tra\c{c}ables la somme sur $n,k$ est finie. Le r\'esultat essentiel est que la s\'erie renormalis\'ee $W_R(A)$ correspond exactement \`a la cochaine \^eta du quasihomomorphisme $\rho$, \'evalu\'ee sur le caract\`ere de Chern $\ch(u)\in HP_1(S\Ac)$. En effet, puisque l'alg\`ebre $\Bc=\cinf(S^1)$ est quasi-libre on peut choisir l'extension triviale $\Rc=\Bc$ et son homologie cyclique est calcul\'ee par le $X$-complexe isomorphe au complexe de de Rham $X(\Bc): \cinf(S^1)\stackrel{\dd}{\to} \Om^1(S^1)$. Les composantes de la cocha\^{\i}ne \^eta renormalis\'ee $\etah_{R0}^{2n+1}: \Om^{2n+1}\Th(S\Ac)\to \Bc$ sont donn\'ees par les \'equations (\ref{zeta}). On peut alors \'evaluer la transgression $\tch^{2n+1}_R(\rho)=\etah^{2n+1}_{R}\gamma$ sur $\ch(u)$, et exprimer le r\'esultat en fonction du potentiel $A= u_-^{-1}Qu_+ - Q$:
\beq
\lefteqn{\tch^{2n+1}_R(\rho)\cdot \ch(u) = }\non\\
&& \qquad \frac{(-)^n}{\sqrt{2\pi i}} \frac{(n!)^2}{(2n+1)!} \Pf_{z=0} \Tr \left( \left(\frac{Q^{-1}A}{1+Q^{-1}A}\right)^{2n+1}(1+Q^{-1}A/2) |D|_+^{-2z} \right) \ . \non
\eeq
Le d\'eveloppement en s\'erie formelle de cette quantit\'e ne contient que des puissances de $A$ sup\'erieures \`a $2n+1$. Par cons\'equent la somme infinie $\sum_{n=0}^{\infty} \tch_R^{2n+1}(\rho)\cdot\ch(u)$ existe au sens des s\'eries formelles en puissances de $A$.

\begin{proposition}[\cite{P9}]
La fonction de partition $W_R(A)$ renormalis\'ee par (\ref{ren}) est proportionnelle, au sens des s\'eries formelles en puissances du potentiel $A= u_-^{-1}Qu_+ - Q$, \`a la somme des transgressions $\tch_R^{2n+1}(\rho)\cdot\ch(u)$ renormalis\'ees par (\ref{zeta}) modulo un contre-terme:
\be
W_R(A) + P(A)=\sqrt{2\pi i}\, \sum_{n=0}^{\infty}\tch_R^{2n+1}(\rho)\cdot \ch(u) \ ,\label{coin}
\ee
o\`u $P(A)$ est une fonctionnelle locale et polyn\^omiale en $A$ donn\'ee par une somme finie de r\'esidus de fonctions z\^eta. \cqfd
\end{proposition}
La somme $\sqrt{2\pi i}\, \sum_{n=0}^{\infty} \tch_R^{2n+1}(\rho)\cdot\ch(u)$ est donc simplement un autre choix de renormalisation pour la fonction de partition $W(A)$. Prenons ensuite la diff\'erentielle $\dd$ de chaque membre de l'\'equation (\ref{coin}). Le membre de gauche donne l'anomalie $\Delta(\om,A)$, tandis que d'apr\`es la discussion du paragraphe \ref{sprin} le membre de droite correspond \`a l'image de $u$ sous l'application diagonale du diagramme (\ref{di}). Ainsi leurs classes de cohomologie dans $HP_1(\Bc)=H^1(S^1)$ co\"{\i}ncident. Cette discussion se g\'en\'eralise de mani\`ere \'evidente au cas d'une boucle unitaire $u$ dans l'alg\`ebre des matrices $M_{\infty}(\Ac)^+$.

\begin{corollary}[\cite{P7, P9}]\label{cta}
Soit $(\Ac,H,D)$ un triplet spectral r\'egulier $p$-sommable de degr\'e pair, et $W_R(A)$ une renormalisation quelconque de la th\'eorie de jauge chirale qui lui est associ\'ee. Soit $e\in M_{\infty}(\Ac)$ un projecteur repr\'esentant une classe de $K$-th\'eorie $[e]\in \Kt_0(\Ac)$, et $u=1+e\otimes(\beta-1)$ la boucle unitaire qui lui correspond par p\'eriodicit\'e de Bott. Alors l'anomalie chirale $\Delta(\om,A)=\dd W_R(A)$ int\'egr\'ee le long de la boucle de potentiels $A=u_-^{-1}Qu_+-Q$ calcule l'indice
\be
\Ind(eDe) = \frac{1}{2\pi i}\oint \Delta(\omega,A) \ \in\zz\ .
\ee
L'anomalie est donn\'ee par une formule locale pour tout choix de renormalisation, par exemple la somme de r\'esidus (\ref{an}) dans le cas de la renormalisation z\^eta \`a la Ray-Singer. \cqfd
\end{corollary} 

Terminons en mentionnant une formule suggestive pour l'application r\'egulateur $\rho_!:MK_{p+1}(\Ac)\to\cc^{\times}$ associ\'ee au triplet spectral (voir l'exemple \ref{ereg}). Soit $(u,\te)$ un couple repr\'esentant une classe de $K$-th\'eorie multiplicative tel que $u\in \GL_{\infty}(\Ac)$ soit un \'el\'ement unitaire \emph{alg\'ebriquement} homotope \`a 1, dans le sens o\`u il existe un chemin unitaire $v\in \GL_{\infty}(\Ac\otimes\cinf[0,1])$ tel que $v(0)=1$ et $v(1)=u$. On peut d\'efinir un d\'eterminant renormalis\'e $\Det_R(u)\in \cc^{\times}$ en int\'egrant l'anomalie $\dd W_R(A)$ le long du chemin de potentiels $A=v_-^{-1}Qv_+ - Q$. Alors
\be
\rho_!(u,\te)= \exp(\sqrt{2\pi i}\,\ch_R(\rho)\cdot \te)\,\Det_R^{-1}(u) \ . 
\ee
Notons que dans cette formule on doit utiliser la m\^eme renormalisation pour le caract\`ere de Chern $\ch_R(\rho)$ et pour le d\'eterminant $\Det_R(u)$ (on s'arrange pour que le contre-terme $P(A)$ dans (\ref{coin}) soit nul). Ainsi tout changement dans le choix de renormalisation pour le d\'eterminant est automatiquement compens\'e par un facteur de phase dans l'exponentielle. C'\'etait attendu puisque l'application r\'egulateur est canoniquement d\'efinie. Une autre cons\'equence de cette formule est qu'elle permet de calculer l'anomalie multiplicative du d\'eterminant renormalis\'e, voir \cite{P9}.

\chapter{Groupo\"{\i}des conformes et localisation}\label{cgrou}

On expose dans ce chapitre un th\'eor\`eme de l'indice local bas\'e sur une renormalisation sans op\'erateur de Dirac. Consid\'erons un groupe discret $G$ op\'erant par transformations conformes sur le plan complexe not\'e $\Si$. Afin d'incorporer des exemples non triviaux on suppose que l'action de tout \'el\'ement $g\in G$ n'est d\'efinie que sur un domaine $\dom(g)\subset\Si$, \'eventuellement vide. L'op\'erateur de Dolbeault d\'efinit naturellement un quasihomomorphisme $p$-sommable entre des compl\'etions convenables du produit crois\'e $\Ac_0=\cinfc(\Si)\cp G$ et de l'alg\`ebre $\Bc_0$ du groupe. Notre objectif est d'utiliser la formule locale d'anomalie. Contrairement \`a la situation du chapitre \ref{ceq}, l'action de $G$ ne pr\'eserve aucune m\'etrique riemannienne sur $\Si$. Il n'est donc pas naturel d'introduire un op\'erateur de Dirac et la renormalization z\^eta n'est pas adapt\'ee. Un choix beaucoup plus judicieux est d'exploiter uniquement la structure complexe du plan. L'anomalie est alors plus facile \`a calculer, et comme pr\'evu se localise automatiquement aux points fixes de l'action de $G$. Le th\'eor\`eme de l'indice qui en r\'esulte s'exprime en fonction de deux cocycles cycliques sur $\Ac_0$. Le premier est une trace qui g\'en\'eralise la formule de Lefschetz aux points fixes d'ordre sup\'erieur. Le deuxi\`eme est un 2-cocycle cyclique construit \`a partir du groupe d'automorphismes modulaires du produit crois\'e. Ce travail fait l'objet de la pr\'epublication \\

\noindent \cite{P10} D. Perrot: Localization over complex-analytic groupoids and conformal renormalization, preprint arXiv:0804.3969.\\

Notons que les r\'esultats pr\'esent\'es ici g\'en\'eralisent ceux obtenus en utilisant l'alg\`ebre de Hopf de Connes et Moscovici et publi\'es dans \\

\noindent \cite{P2} D. Perrot: A Riemann-Roch theorem for one-dimensional complex groupoids, {\it Comm. Math. Phys.} {\bf 218} (2001) 373-391.

\section{Quasihomomorphisme de Dolbeault}

Dans tout le chapitre $\Si=\cc$ d\'esigne le plan complexe vu comme surface de Riemann, avec $z$ son syst\`eme de coordonn\'ees complexes. On note $\Omc^*(\Si)$ l'alg\`ebre des formes diff\'erentielles complexes \`a support compact sur $\Si$. La sous-alg\`ebre des 0-formes est donc isomorphe aux fonctions lisses \`a support compact $\cinfc(\Si)$, et l'espace des 1-formes se scinde en somme directe $\Omc^{1,0}(\Si)\oplus \Omc^{0,1}(\Si)$ des formes proportionnelles respectivement \`a $dz$ et $d\zb$. De m\^eme la diff\'erentielle de de Rham $d=\d+\deb$ se scinde en la somme de $\d=dz\d_z$ et de l'op\'erateur de Dolbeault $\deb=d\zb\d_{\zb}$. On notera
\be
Q= \deb\, :\, \cinfc(\Si)\to \Omc^{0,1}(\Si)
\ee
sa restriction \`a l'espace des 0-formes. D\'esignons improprement par $Q^{-1}:\Omc^{0,1}(\Si)\to \cinf(\Si)$ l'op\'erateur de Green associ\'e \`a $Q$. Son noyau distributionnel sur $\Si\times\Si$ est proportionnel au noyau de Cauchy \cite{P10}: en deux points quelconques $w,z$ du plan 
\be
Q^{-1}(z,w) = \frac{1}{\pi(z-w)}\ .
\ee
Puisque la multiplication des fonctions scalaires par une 1-forme $A=d\zb A_{\zb}\in \Omc^{0,1}(\Si)$ induit une application $\cinfc(\Si)\to \Omc^{0,1}(\Si)$, la compos\'ee $Q^{-1} A$ est un op\'erateur $\cinfc(\Si)\to\cinf(\Si)$. On veut l'\'etendre en un op\'erateur compact sur un espace de Hilbert. Pour tout poids $\al\in\rr$, d\'efinissons l'espace de Hilbert $H_{\al}$ comme la compl\'etion de $\cinfc(\Si)$ pour la norme 
$$
\|\xi\|_{\al}=\Big(\int_{\Si} d^2z\, (1+|z|)^{\al} |\xi(z)|^2 \Big)^{1/2}\qquad \forall \xi\in\cinfc(\Si)\ ,
$$
o\`u $d^2z=d\zb\wedge dz/2i$ est la forme volume euclidienne. Comme cons\'equence du lemme de Rellich \cite{Gi} on obtient l'estimation suivante.

\begin{lemma}[\cite{P10}]\label{lcom}
Pour toute 1-forme $A\in\Omc^{0,1}(\Si)$, la compos\'ee $Q^{-1}A$ s'\'etend en un op\'erateur compact sur $H_{\al}$ d\`es que $\al<-1$. Dans ce cas $Q^{-1}A$ est dans la classe de Schatten $\ell^p(H_{\al})$ pour tout $p>2$. \cqfd
\end{lemma}

Consid\'erons maintenant $G$ un groupe \emph{discret} agissant sur $\Si$ par transformations conformes de la mani\`ere suivante \cite{P10}. A tout \'el\'ement $g\in G$ on associe deux ouverts (\'eventuellement vides) $\dom(g)\subset \Si$ et $\ran(g)\subset \Si$ ainsi qu'une transformation conforme inversible $\dom(g)\to\ran(g)$, avec la condition $\dom(gh)\supset h^{-1}(\dom(g))\cap \dom(h)$ pour tous $g,h\in G$. Par exemple $G$ est un sous-groupe discret de $\SL(2,\cc)$ agissant sur le plan par homographies. Insistons sur le fait qu'en g\'en\'eral, la transformation conforme $\dom(g)\to \ran(g)$ ne caract\'erise pas $g$ de mani\`ere unique comme \'el\'ement du groupe $G$. Par exemple, $G$ est un groupe quelconque agissant trivialement sur $\Si$, avec $\dom(g)=\Si$ pour tout $g\in G$. \\

D\'esignons ensuite par $\Ac_0$ l'espace engendr\'e par les sommes finies de symboles $fU^*_g$ avec $f\in\cinfc(\Si)$ et $g\in G$ tels que $\supp f\subset\dom(g)$. Le produit de convolution $(f_1U^*_{g_1})(f_2U^*_{g_2})= f_1(f_2\circ g_1)U^*_{g_2g_1}$ d\'efinit une structure d'alg\`ebre associative sur $\Ac_0$ que l'on \'ecrira comme un produit crois\'e
\be
\Ac_0= \cinfc(\Si)\cp G\ .
\ee
L'alg\`ebre $\Ac_0$ est naturellement repr\'esent\'ee lin\'eairement sur l'espace $\Omc^*(\Si)$ et respecte le bidegr\'e des formes diff\'erentielles: pour tout $fU^*_g\in \Ac_0$ on notera $fr(g)_+$ sa repr\'esentation sur $\cinfc(\Si)$ et $fr(g)_-$ sa repr\'esentation sur $\Omc^{0,1}(\Si)$. Soit $\Bc_0$ l'alg\`ebre de convolution du groupe $G$: c'est l'espace engendr\'e par sommes finies de symboles $U^*_g$ muni du produit $U^*_{g_1}U^*_{g_2} = U^*_{g_2g_1}$. Choisissons une compl\'etion $\Bc\supset\Bc_0$ en $m$-alg\`ebre de Fr\'echet. Pour tout choix de poids $\al<-1$, il existe deux homomorphismes $\rho_+,\rho_-: \Ac_0\to\End(H_{\al})\otimes\Bc$ d\'efinis par
$$
\rho(fU^*_g)_+ = fr(g)_+\otimes U^*_g\ , \qquad \rho(fU^*_g)_- = Q^{-1}fr(g)_- Q \otimes U^*_g\ .
$$ 
Comme l'op\'erateur de Dolbeault est invariant conforme, on a $r(g)_-Q=Qr(g)_+$ et le lemme \ref{lcom} implique que la diff\'erence $\rho_+-\rho_-$ est \`a valeurs dans l'id\'eal $\ell^p(H_{\al})\otimes\Bc$, $p>2$. En posant $\Ic=\ell^p(H_{\al})$ on obtient donc un quasihomomorphisme $p$-sommable de degr\'e pair
\be
\rho: \Ac\to \Ec^s\triangleright\Ic^s\hotimes\Bc
\ee
pour une compl\'etion ad\'equate $\Ac\supset\Ac_0$ en $m$-alg\`ebre de Fr\'echet \cite{P10}. Puisque l'image de $\rho_{\pm}$ est contenue dans un produit tensoriel $\End(H_{\al})\otimes\Bc$, le quasihomomorphisme est automatiquement admissible relativement \`a toute extension quasi-libre $0\to \Jc \to \Rc\to \Bc\to 0$. On pourra prendre par exemple l'alg\`ebre tensorielle $\Rc=T\Bc$. Le th\'eor\`eme \ref{trr} appliqu\'e aux invariants primaires donne le diagramme commutatif (\ref{diag}):
\be
\vcenter{\xymatrix{
\Kt_j(\Ic\hotimes\Ac) \ar[r]^{\rho_!} \ar[d] \ar@{.>}[rd]^{\Delta} & \Kt_j(\Ic\hotimes\Bc) \ar[d] \\
HP_j(\Ac) \ar[r]_{\ch(\rho)} & HP_j(\Bc) }} \qquad j\in\zz_2 \label{didi}
\ee
Par p\'eriodicit\'e de Bott, il suffit de consid\'erer comme d'habitude la $K$-th\'eorie impaire ($j=1$). Nous avons donn\'e dans \cite{P10} une formule locale de l'incice, en calculant la diagonale $\Delta$ sous la forme d'une anomalie chirale associ\'ee \`a une th\'eorie des champs conforme, d'apr\`es la m\'ethode expos\'ee au chapitre \ref{cano}. La renormalisation effectu\'ee dans cette situation n\'ecessite cependant de se restreindre \`a la $K$-th\'eorie de la sous-alg\`ebre dense $\Ac_0$. Consid\'erons donc un \'el\'ement inversible $u\in (\Ac_0)^+$ repr\'esentant une classe $[u]\in \Kt_1(\Ic\hotimes\Ac)$, et choisissons un rel\`evement inversible arbitraire $\uh\in (\Th\Ac_0)^+$. L'image de $\uh$ sous les rel\`evements de $\rho_{\pm}$ en des homomorphismes $(\rho_*)_{\pm}:\Th\Ac_0\to \End(H_{\al})\otimes\Rch$ permet de d\'efinir deux inversibles $\uh_+$ et $\uh_-$ par 
$$
\rho_*(\uh)_+ = \uh_+\ ,\qquad \rho_*(\uh)_-=Q^{-1}\uh_-Q\ .
$$
On introduit ensuite un potentiel de jauge $A\in\hom(\cinfc(\Si),\Omc^{0,1}(\Si))\otimes\Rch$ par la formule (\ref{pot}),
\be
A= \uh_-^{-1}Q\uh_+ - Q\ ,
\ee
de sorte que l'op\'erateur $Q^{-1}A=Q^{-1}\uh_-^{-1} Q \uh_+ -1$ s'\'etende en un \'el\'ement de l'id\'eal $\Ic\otimes\Rch\subset \End(H_{\al})\otimes\Rch$. En combinant la trace $\Tr$ des op\'erateurs sur $H_{\al}$ avec la trace universelle $\nat:\Rch\to \Rch_{\nat}=\Rch/[\ ,\ ]$ la fonctionnelle d'action quantique est d\'efinie comme s\'erie formelle en puissances de $A$
\be
W(A)= \sum_{n\geq 1}W^n(A)\ ,\qquad W^n(A)= \frac{(-)^{n+1}}{n}\Tr\nat((Q^{-1}A)^n)\in \Rch_{\nat}\ .
\ee
Puisque $\Ic=\ell^p(H_{\al})$ pour $p>2$, la trace d'op\'erateurs n'a de sens que pour $n\geq 3$ et seuls les termes de plus bas degr\'e $W^1(A)$ et $W^2(A)$ n\'ecessitent une renormalisation. L'anomalie $\Delta(\om,A)$, qui correspond \`a l'image de la s\'erie formelle renormalis\'ee $W_R(A)$ sous l'application de bord $\dd:\Rch_{\nat}\to\Om^1\Rch_{\nat}$, est donc un polyn\^ome en $A$ de degr\'e au plus 2 (voir (\ref{gr})) et d\'efinit une classe d'homologie cyclique dans $HP_1(\Bc)$ ind\'ependante de la renormalisation choisie. Notons que la construction pr\'ec\'edente se g\'en\'eralise de mani\`ere \'evidente au cas d'un inversible dans l'alg\`ebre des matrices $u\in M_{\infty}(\Ac_0)^+\subset (\Ic\hotimes\Ac)^+$. La discussion du chapitre \ref{cano} implique donc la proposition suivante.

\begin{proposition}[\cite{P10}]\label{pano}
Soit $u\in M_{\infty}(\Ac_0)^+$ un \'el\'ement inversible repr\'esentant une classe dans $\Kt_1(\Ic\hotimes\Ac)$. Alors pour toute renormalisation de la s\'erie formelle $W_R(A)$ associ\'ee au potentiel de jauge $A=\uh_-^{-1}Q \uh_+ - Q$, l'anomalie
\be
\Delta(\om,A) =\dd W_R(A) \equiv \sqrt{2\pi i}\, \ch(\rho_!(u)) \ \in HP_1(\Bc)
\ee
est un polyn\^ome en $A$ de degr\'e au plus $2$ qui calcule la diagonale du diagramme commutatif (\ref{didi}). \cqfd
\end{proposition}

\section{Renormalisation conforme}

Nous allons maintenant renormaliser les deux premiers termes $W^1(A)$ et $W^2(A)$ de la fonctionnelle d'action quantique associ\'ee au potentiel de jauge
$$
A\in \hom(\cinfc(\Si),\Omc^{0,1}(\Si))\otimes\Rch\ ,
$$
en exploitant uniquement la structure complexe de $\Si$. Rappelons que $r(g)_+$ d\'esigne l'action de $g\in G$ sur les fonctions $f\in\cinfc(\Si)$ dont le support est contenu dans $\dom(g)$. On peut donc d\'ecomposer le potentiel en une somme
$$
A=\sum_{g\in G} A(g) r(g)_+\quad \mbox{avec} \quad A(g)\in \Omc^{0,1}(\Si)\otimes\Rch\ ,
$$
o\`u l'espace des 1-formes $\Omc^{0,1}(\Si)$ est vu dans $\hom(\cinfc(\Si),\Omc^{0,1}(\Si))$ par multiplication sur $\cinfc(\Si)$. Dans le syst\`eme de coordonn\'ees complexes $z$ sur $\Si$ \'ecrivons $A(g)=d\zb A_{\zb}(g)$ et regardons chaque composante $A_{\zb}(g)\in \cinfc(\Si)\otimes\Rch$ comme une fonction test sur $\Si$ \`a valeurs dans $\Rch$. Un calcul na\"{\i}f au moyen du noyau distributionnel de $Q^{-1}$ donne (\cite{P10})
$$
W^1(A) = \Tr\nat(Q^{-1}A) = \sum_{g\in G} \int_{\Si} d^2z\, \frac{\nat A_{\zb}(g,z)}{\pi(g(z)-z)}\ ,
$$
o\`u $g(z)$ est fonction holomorphe de $z$. A priori cette expression n'a pas de sens. \emph{Renormaliser} $W^1(A)$ revient \`a prolonger la fonction $z\mapsto 1/(g(z)-z)$ en une distribution sur $\Si$. La difficult\'e provient bien entendu de ses p\^oles, autrement dit des \emph{points fixes} de la transformation $g$. On dit qu'un point fixe $z_0\in\dom(g)$ est d'ordre $n\geq 1$ si $g(z)-z$ se comporte comme $(z-z_0)^n$ au voisinage de $z_0$. Le cas $n=\infty$ signifie que tous les points sont fixes au voisinage de $z_0$. Le prolongement distributionnel de la fonction $1/(g(z)-z)$ en un point fixe d\'epend de son ordre:\\

\noindent $\bullet$ Si $n=\infty$ alors $1/(g(z)-z)$ n'a aucun sens autour de $z_0$. Dans ce cas on assigne une valeur quelconque \`a cette fonction, par exemple z\'ero (le choix le plus simple).\\

\noindent $\bullet$ Si $n<\infty$, alors $z_0$ est n\'ecessairement un point fixe isol\'e. Remarquons que pour $z\neq z_0$ on a une \'egalit\'e de fonctions
$$
\frac{1}{g(z)-z} = \frac{1}{(z-z_0)^n}\, H^n_{g,z_0}(z)\quad \mbox{avec} \quad H^n_{g,z_0}(z) := \frac{(z-z_0)^n}{g(z)-z}\ ,
$$
et $H^n_{g,z_0}$ est holomorphe sur un voisinage de $z_0$. Il suffit donc de construire un prolongement distributionnel de la fonction m\'eromorphe $1/(z-z_0)^n$. On peut \'ecrire
\be
\frac{1}{(z-z_0)^n} = \frac{(-)^{n-1}}{(n-1)!}\, \d_z^{n-1}\left(\frac{1}{z-z_0}\right)\ , \label{conf}
\ee
et le membre de droite d\'efinit bien une distribution sur un voisinage de $z_0$. En proc\'edant de la sorte en tous les points fixes on obtient le terme renormalis\'e $W_R^1(A)$. D'autres renormalisations sont possibles, mais celle d\'ecrite ici est la seule \emph{invariante conforme}, c'est-\`a-dire ind\'ependante du syst\`eme de coordonn\'ees complexes choisi. \\

La renormalisation du terme $W_R^2(A)$ est analogue. Dans ce cas il faut prolonger la fonction de deux variables complexes $(z,w)\mapsto 1/(h(w)-z)(g(z)-w)$ en une distribution pour tous $g,h\in G$. Les formules sont aussi bas\'ees sur le prolongement (\ref{conf}) et nous renvoyons \`a \cite{P10} pour plus de d\'etails. La s\'erie formelle $W_R(A)$ est alors bien d\'efinie et on peut calculer l'anomalie $\Delta(\om,A)=\dd W_R(A)$ au moyen des \'equations BRS (chapitre \ref{cano}) $-\dd A = (Q+A)\om_+ + \om_-(Q+A)$. Ici la forme de Maurer-Cartan se d\'ecompose de mani\`ere analogue au potentiel $A$, 
$$
\om_{\pm}=\sum_{g\in G} \om(g)r(g)_{\pm}\quad \mbox{avec} \quad \om(g)\in \cinfc(\Si)\otimes \Om^1\Rch\ .
$$
Puisque les probl\`emes de renormalisation sont concentr\'es aux points fixes de l'action de $G$ sur $\Si$, il n'est pas \'etonnant de constater que l'anomalie est aussi localis\'ee aux points fixes:

\begin{proposition}[\cite{P10}]\label{pan}
L'anomalie associ\'ee \`a la renormalisation conforme est un polyn\^ome de degr\'e $1$ par rapport \`a $A$. Sa composante de degr\'e z\'ero $\Delta^0(\om,A)$ est une somme sur les points fixes isol\'es:
\be
\Delta^0(\om,A)=\sum_{g\in G}\sum_{\substack{z_0=g(z_0)\\ \textup{isol\'e}}} \frac{-1}{(n-1)!}\, \d_z^{n-1}\big(H_{g,z_0}^n(z) \nat\om(g,z)\big)_{z=z_0}\ , \label{d0}
\ee
o\`u $n\in \nn^*$ d\'enote l'ordre de $z_0$. La composante de degr\'e un $\Delta^1(\om,A)$ est une int\'egrale sur la vari\'et\'e complexe des points fixes non isol\'es ( = d'ordre infini):
\be
\Delta^1(\om,A)=\frac{1}{\pi} \sum_{g,h\in G}\int_{z=hg(z)}d^2z\, \nat \big(\d_z- \frac{1}{2}\d_z\ln g'(z)\big)A_{\zb}(g,z)\,\,\om(h,g(z))\ , \label{d1}
\ee
o\`u $g'(z)$ est la d\'eriv\'ee de la fonction holomorphe $g(z)$. \cqfd
\end{proposition}
La preuve est un calcul direct. Il convient de remarquer que la contribution d'un point fixe d'ordre $n$ dans $\Delta^0(\om,A)$ d\'epend uniquement des d\'eriv\'ees de $g$ d'ordre $\leq 2n-1$.  Par exemple aux plus bas ordres le calcul donne
\beq
\lefteqn{\frac{-1}{(n-1)!}\, \d_z^{n-1}\big(H_{g,z_0}^n(z) \nat\om(g,z)\big)_{z=z_0} = } \non\\
&&\qquad (n=1):\qquad \frac{1}{1-g'(z_0)}\, \nat\om(z_0)    \non\\
&&\qquad (n=2):\qquad \frac{2}{g''(z_0)} \left( \frac{1}{3} \frac{g'''(z_0)}{g''(z_0)} \, \nat\om(z_0) -  \nat\d_z\om(z_0) \right)   \non\\
&&\qquad (n=3):\qquad \frac{3}{2g'''(z_0)} \left( \frac{1}{10} \frac{g^{(5)}(z_0)}{g'''(z_0)}\right.  \nat\om(z_0) - \frac{1}{8} \left(\frac{g^{(4)}(z_0)}{g'''(z_0)}\right)^2 \, \nat\om(z_0)  \non\\
&& \qquad \qquad \qquad \qquad \qquad \qquad \qquad \left. + \frac{1}{2} \frac{g^{(4)}(z_0)}{g'''(z_0)} \, \nat \d_z\om(z_0) - \nat \d_z^2\om(z_0) \right)   \non
\eeq
On retrouve donc le nombre de Lefschetz bien connu dans le cas $n=1$, tandis que pour $n>1$ les jets d'ordre sup\'erieur de la transformation $g$ interviennent.

\section{Th\'eor\`eme de l'indice}

Nous pouvons maintenant calculer la diagonale du diagramme (\ref{didi}) restreinte \`a la sous-alg\`ebre $\Ac_0$. D\'efinissons l'ensemble
\be
\Gamma = \coprod_{g\in G} \dom(g) =  \{(g,z)\in G\times \Si\ |\ z\in \dom(g)\}\ .
\ee
$\Gamma$ muni de la loi de composition partielle $(g,z)\cdot (h,w)=(gh,w)$ pour $z=h(w)\in \dom(g)$ est un groupo\"{\i}de \'etale. Remarquons que $\Ac_0=\cinfc(\Si)\cp G$ s'identifie \`a l'alg\`ebre de convolution des fonctions lisses \`a support compact sur $\Gamma$.\\
Les automorphismes de $\Gamma$ correspondent aux couples $\gamma_0=(g,z_0)$ form\'es d'un \'el\'ement $g\in G$ agissant par transformation conforme sur $\Si$ et d'un point fixe $z_0\in \dom(g)$. L'ensemble de tous les automorphismes de $\Gamma$ est la r\'eunion du sous-ensemble discret $\Gamma_f$ des automorphismes isol\'es (ordre $n<\infty$), et de la vari\'et\'e complexe de dimension un $\Gamma_{\infty}$ des automorphismes non isol\'es (ordre $n=\infty$). Un examen attentif des formules (\ref{d0}) et (\ref{d1}) permet de ``deviner'' certains cocycles cycliques sur l'alg\`ebre de convolution du groupo\"{\i}de $\Gamma$. 

\begin{lemma}[\cite{P10}]\label{linv}
Soit $\gamma_0=(g,z_0)\in\Gamma_f$ un automorphisme isol\'e d'ordre $n\in \nn^*$. La fonctionnelle lin\'eaire $\Ac_0\to\cc$ donn\'ee par
$$
a\mapsto \d_z^{n-1} \big( H^n_{g,z_0}(z)\,a(g,z) \big)_{z=z_0}
$$
est ind\'ependante du choix de syst\`eme de coordonn\'ees complexes $z$ choisi. Ecrivons $\gamma=(g,z)\in \Gamma$ dans un voisinage de $\gamma_0$ muni de sa structure complexe, $H^n_{g,z_0}(z)=H^n_{\gamma_0}(\gamma)$ et identifions $\d_z$ et $\d_{\gamma}$. Alors en sommant sur tous les automorphismes isol\'es, la fonctionnelle $\Phi(\Gamma):\Ac_0 \to \cc$ d\'efinie par
$$
\Phi(\Gamma)(a)= \sum_{\gamma_0 \in\Gamma_f} \frac{-1}{(n-1)!}\, \d_{\gamma}^{n-1} \big( H^n_{\gamma_0}(\gamma)\,a(\gamma) \big)_{\gamma=\gamma_0}
$$
est une trace sur l'alg\`ebre $\Ac_0$. \cqfd
\end{lemma}
Ainsi la composante de degr\'e z\'ero de l'anomalie $\Delta^0(\om,A)$ est essentiellement une trace \'evalu\'ee sur $\om$, d\'efinie uniquement \`a partir de la structure complexe de $\Gamma$. La composante de degr\'e un $\Delta^1(\om,A)$ fait quant \`a elle intervenir un analogue non-commutatif de la classe de Todd \cite{P2}. Remarquons d'abord que les trois diff\'erentielles $\d$, $\deb$ et $d=\d+\deb$ sur l'alg\`ebre $\Omc^*(\Si)$ commutent avec les transformations conformes, donc s'\'etendent en des diff\'erentielles sur le produit crois\'e $\Omc^*(\Si)\cp G$. Il existe une quatri\`eme diff\'erentielle provenant du \emph{groupe d'automorphismes modulaires}: son g\'en\'erateur est une d\'erivation $D$ sur $\Omc^*(\Si)\cp G$,
$$
D(fU^*_g)= \ln|g'|^2 fU^*_g\ ,\qquad \forall\ f\in \Omc^*(\Si)\ ,\ g\in G\ ,
$$
o\`u la fonction scalaire $z\mapsto |g'(z)|^2$ mesure la dilatation du volume euclidien sur $\Si$ induite par la transformation $g$. Le commutateur de la d\'erivation $D$ avec la diff\'erentielle $\d$ d\'efinit donc une nouvelle diff\'erentielle 
\be
\delta = [\d,D]\ ,\qquad \delta^2=0\ ,
\ee
qui anticommute avec $d$, $\d$, et $\deb$. On a explicitement $\delta(fU^*_g)=(\d\ln g')fU^*_g$. Alors l'int\'egration des 2-formes sur la vari\'et\'e $\Gamma_{\infty}$ orient\'ee par sa structure complexe permet de construire des $2$-cocycles cycliques sur la sous-alg\`ebre $\Ac_0\subset \Omc^*(\Si)\cp G$:

\begin{lemma}[\cite{P2, P10}]
La classe fondamentale du groupo\"{\i}de $\Gamma$ est le 2-cocycle cyclique sur $\Ac_0$ d\'efini par la fonctionnelle  
$$
[\Gamma] (a_0\dd a_1 \dd a_2) = \int_{\Gamma_{\infty}} a_0da_1da_2\qquad \forall a_i\in \Ac_0\ .
$$
La classe de Chern du groupo\"{\i}de est le 2-cocycle cyclique 
$$
c_1(\Gamma) (a_0 \dd a_1 \dd a_2) = \int_{\Gamma_{\infty}} a_0(da_1\delta a_2 + \delta a_1 da_2)\qquad \forall a_i\in \Ac_0\ .
$$
La somme $\Td(\Gamma) :=  [\Gamma] -\frac{1}{2} c_1(\Gamma)$ est appel\'ee \emph{classe de Todd} du groupo\"{\i}de. \cqfd
\end{lemma}
La classe de Todd admet une expression simple en fonction de la diff\'erentielle $\nabla=d-\frac{1}{2}\delta$. En effet $\delta a_1\delta a_2=0$ pour des raisons dimensionnelles et 
\be
\Td(\Gamma)(a_0 \dd a_1 \dd a_2) = \int_{\Gamma_{\infty}} a_0\nabla a_1\nabla a_2\ . \label{tod}
\ee
Ce 2-cocycle cyclique n'est cependant pas enti\`erement canonique car le g\'en\'erateur $D$ du groupe modulaire d\'epend du choix de la mesure euclidienne en plus de la structure complexe sur $\Si$. Comme il n'y a aucune raison de pr\'ef\'erer la mesure euclidienne on peut aussi bien repr\'esenter le groupe modulaire au moyen d'une forme volume quelconque sur $\Si$. Nous avons montr\'e dans \cite{P2} comment modifier en cons\'equence le cocycle (\ref{tod}) sans changer sa classe de cohomologie: un terme proportionnel \`a la courbure de la m\'etrique de K\"ahler appara\^{\i}t. On retrouve ainsi l'expression de la classe de Todd d'une surface de Riemann.\\

Soit maintenant $\Bc$ l'alg\`ebre de convolution du groupe discret $G$ compl\'et\'ee en $m$-alg\`ebre de Fr\'echet. On d\'efinit un homomorphisme $\tilde{\rho}: \Ac_0\to \Ac_0\otimes \Bc$ en posant $\tilde{\rho}(f U^*_g)=f U^*_g\otimes U^*_g$. Si $e\in M_{\infty}(\Ac_0)$ est un idempotent et $u\in M_{\infty}(\Ac_0)^+$ un inversible, les caract\`eres de Chern de leurs images sous $\tilde{\rho}$ sont des classes d'homologie cyclique p\'eriodique
$$
\ch(\tilde{\rho}(e)) \in HP_0(\Ac_0\otimes\Bc)\ ,\qquad \ch(\tilde{\rho}(u)) \in HP_1(\Ac_0\otimes\Bc)\ ,
$$
o\`u l'alg\`ebre $\Ac_0\otimes\Bc$ est consid\'er\'ee comme discr\`ete. D'autre part, toute classe de cohomologie cyclique p\'eriodique $\varphi$ sur $\Ac_0$ induit une application cap-produit
$$
\varphi\cap : HP_*(\Ac_0 \otimes \Bc)\to HP_*(\Bc)\ .
$$
En utilisant la formule locale d'anomalie, nous avons montr\'e dans \cite{P10} que la diagonale du diagramme (\ref{didi}) restreinte \`a la sous-alg\`ebre $\Ac_0\subset\Ac$ se ram\`ene \`a un cap-produit avec la classe de cohomologie cyclique $\varphi=\Phi(\Gamma)+ \Td(\Gamma)$:

\begin{theorem}[\cite{P10}]\label{tcup}
Soit $e\in M_{\infty}(\Ac_0)$ un idempotent et $u\in M_{\infty}(\Ac_0)^+$ un inversible. Alors les caract\`eres de Chern de leurs images directes $\rho_!(e)\in \Kt_0(\Ic\hotimes\Bc)$ et $\rho_!(u)\in \Kt_1(\Ic\hotimes\Bc)$ sont donn\'es par les cap-produits
\beq
\ch(\rho_!(e)) &=& (\Phi(\Gamma)+ \Td(\Gamma))\cap \ch(\tilde{\rho}(e)) \ \in HP_0(\Bc)\ ,\non\\
\ch(\rho_!(u)) &=& (\Phi(\Gamma)+ \Td(\Gamma))\cap \ch(\tilde{\rho}(u)) \ \in HP_1(\Bc)\ ,\non
\eeq
o\`u $\tilde{\rho}: \Ac_0\to\Ac_0\otimes\Bc$ est l'homomorphisme canonique. \cqfd
\end{theorem}
On peut donner des formules explicites. L'homomorphisme $\tilde{\rho}$ se rel\`eve en un homomorphisme de pro-alg\`ebres $\tilde{\rho}_*:\Th\Ac_0\to \Ac_0\otimes\Rch$. En oubliant l'alg\`ebre des matrices $M_{\infty}$ pour simplifier l'\'ecriture, l'idempotent $\tilde{\rho}(e)\in \Ac_0\otimes\Bc$ se rel\`eve donc en un idempotent $\tilde{e}=\tilde{\rho}_*(\eh)\in \Ac_0\otimes\Rch$ et de m\^eme l'inversible $\tilde{\rho}(u)\in (\Ac_0\otimes\Bc)^+$ se rel\`eve en un inversible $\tilde{u}=\tilde{\rho}_*(\uh)\in (\Ac_0\otimes\Rch)^+$. On a alors
\beq
\ch(\rho_!(e)) &=& \Phi(\Gamma) \nat (\et) -  \int_{\Gamma_{\infty}} \nat\, \frac{\et\nabla\et\nabla\et}{2\pi i}  \ \in \Rch_{\nat}\ , \\
\ch(\rho_!(u)) &=& \Phi(\Gamma) \nat \Big(\frac{\ut^{-1}\dd\ut}{\sqrt{2\pi i}}\Big) -  \int_{\Gamma_{\infty}} \nat\, \frac{\ut^{-1}\nabla\ut\nabla\ut^{-1}\dd\ut}{2(2\pi i)^{3/2}} \ \in \Om^1\Rch_{\nat}\ . \label{expl}
\eeq
Le th\'eor\`eme \ref{tcup} se d\'emontre en remarquant que pour $A=\uh_-^{-1}Q\uh_+-Q$, l'anomalie $\Delta(\om,A)$ donn\'ee par la proposition \ref{pan} co\"{\i}ncide avec le membre de droite dans (\ref{expl}) modulo un bord (et un facteur $\sqrt{2\pi i}$). Le cas pair s'en d\'eduit par p\'eriodicit\'e de Bott. 

\begin{remark}\textup{Dans le cas o\`u $\Gamma_f=\varnothing$ et $\Gamma_{\infty}=\Si$, le cocycle cyclique $\varphi=\Phi(\Gamma)+ \Td(\Gamma)$ se r\'eduit \`a la classe de Todd
$$
\Td(\Gamma)(a_0 \dd a_1 \dd a_2) = \int_{\Si} a_0 \nabla a_1\nabla a_2 \ .
$$
Nous avions d\'ej\`a obtenu cette formule dans \cite{P2}, en utilisant l'approche de Connes et Moscovici par l'alg\`ebre de Hopf des diff\'eomorphismes \cite{CM98}. Dans cette situation, $G$ est un pseudogroupe de transformations conformes dont l'action est relev\'ee au fibr\'e $P$ des m\'etriques de K\"ahler au-dessus de $\Si$, et l'on obtient un $K$-cycle sur l'alg\`ebre $\cinfc(P)\cp G$ en combinant l'op\'erateur de Dolbeault horizontal et l'op\'erateur de signature vertical. Pour cette raison un facteur 2 global appara\^{\i}t dans la formule de \cite{P2}. Notons que la diff\'erentielle $\delta=[\d,D]$ est l'un des g\'en\'erateurs de l'alg\`ebre de Hopf de Connes et Moscovici. }
\end{remark}

\end{document}